\documentclass[reqno]{amsart}
\usepackage{amssymb,enumerate,mathrsfs}

\numberwithin{equation}{section}

\newtheorem{theorem}[equation]{Theorem}
\newtheorem{proposition}[equation]{Proposition}
\newtheorem{lemma}[equation]{Lemma}
\newtheorem{corollary}[equation]{Corollary}

\theoremstyle{definition}
\newtheorem{definition}[equation]{Definition}
\newtheorem{remark}[equation]{Remark}
\newtheorem{example}[equation]{Example}

\DeclareMathOperator{\sym}{ \sigma\!\!\!\sigma}

\def\C{\mathbb C}
\def\Dom{\mathcal D}
\def\K{\mathcal K}
\def\KK{\mathcal K}
\def\L{\mathscr L}
\def\N{\mathbb N}
\def\R{\mathbb R}
\def\RR{\mathcal R}
\def\Sing{\mathcal E}
\def\Vee{\mathcal V}
\def\W{\mathcal W}

\def\Left{\mathrm L}
\def\Right{\mathrm R}
\def\Gr{\mathrm{Gr}}
\def\Hol{\mathfrak{H}}
\def\Mero{\mathfrak{M}}

\def\bev{\,{}^b\hspace{-0.5pt}{\ev}}
\def\bT{\,{}^b\hspace{-0.5pt}T}
\def\bpi{\,{}^b\!\pi}
\def\bsym{\,{}^b\!\!\sym}

\def\cev{\,{}^c {\ev}}
\def\csym{\,{}^c\!\sym}
\def\cT{\,{}^c T}
\def\cpi{\,{}^c\hspace{-1.5pt}\pi}

\def\cn{\mathrm{cn}}

\def\e{\mathrm e}
\def\embed{\hookrightarrow}
\def\eps{\varepsilon}

\def\id{I}
\def\ie{i.e.}
\def\im{i}
\def\m{\mathfrak m}

\def\open#1{\smash[t]{\overset{{}_{\circ}}{#1}{}}}
\def\open#1{\smash[t]{\overset{{}_{\,\,\circ}}{#1}{}}}

\def\set#1{\{#1\}}

\def\minus{\backslash}
\def\Wedge{\raise2ex\hbox{$\mathchar"0356$}}

\def\contr{{\mathchar"0263}}

\def\neutral#1{}  %this is used only to balance parentheses

\DeclareMathOperator{\bgres}{bg-res}
\DeclareMathOperator{\bgspec}{bg-spec}

\DeclareMathOperator{\ev}{ev}
\DeclareMathOperator{\Ind}{ind}

\DeclareMathOperator{\LinSpan}{span}
\DeclareMathOperator{\Diff}{Diff}
\DeclareMathOperator{\rg}{rg}
\DeclareMathOperator{\resolv}{res}
\DeclareMathOperator{\spec}{spec}

\begin{document}
\title[Geometry and spectra of cone operators]{Geometry and spectra of closed extensions of elliptic cone operators}
\author{Juan B. Gil}
\address{Penn State Altoona\\ 3000 Ivyside Park \\
 Altoona, PA 16601-3760}
\author{Thomas Krainer}
\address{Institut f\"ur Mathematik\\ Universit\"at Potsdam\\
 14415 Potsdam, Germany}
\author{Gerardo A. Mendoza}
\address{Department of Mathematics\\ Temple University\\
 Philadelphia, PA 19122}
\begin{abstract}
We study the geometry of the set of closed extensions of index $0$ of an elliptic differential cone operator and its model operator in connection with the spectra of the extensions, and give a necessary and sufficient condition for the existence of rays of minimal growth for such operators.
\end{abstract}

\maketitle

\section{Introduction}

The purpose of this paper is to study the spectra and resolvents of the closed extensions of an elliptic differential cone operator $A$ on a compact manifold $M$ with boundary, and of its model operator $A_\wedge$. It is well known that the closed extensions of $A$ are in one to one correspondence with the subspaces of a finite dimensional space, $\Dom_{\max}/\Dom_{\min}$, the spaces $\Dom_{\max}$ and $\Dom_{\min}$ being certain subpaces determined by $A$ of an $L^2$ space on $M$, cf. Lesch \cite{Le97}. It is thus natural to view the extensions as corresponding to points in the various Grassmannians associated with $\Dom_{\max}/\Dom_{\min}$. Extending this, we develop a viewpoint in which issues pertaining spectra and resolvents, both for the closed extensions of $A$ and of $A_\wedge$, are expressed and examined in (finite-dimensional) geometric terms.

Cone differential operators are generalizations of the operators that arise when standard differential operators are written using polar coordinates. Their study is therefore of interest in the context of manifolds with conical singularities, both in themselves and as guiding examples in a general theory of analysis of differential operators on manifolds with other kinds of singularities, cf. Schulze \cite{SzNorthHolland}.

Our motivation for undertaking this study comes from the desirability of executing Seeley's program \cite{Seeley} in the case of elliptic cone operators. This requires a detailed understanding of the resolvent in terms of the symbol of $A$ and the domain of the extension. From the pseudodifferential point of view, the symbol of $A$ is a pair consisting of $A_\wedge$ and its $b$-symbol, or more invariantly, its cone-symbol $\csym(A)$ as defined in Section \ref{sec-TheSymbols}. As in the standard theory of elliptic operators on a manifold without boundary, the statement that a given sector $\Lambda\subset \C $ is a sector of minimal growth for $\csym(A)$ is domain-insensitive. The operator $A_\wedge$, however, is a differential operator, so the analogous statement for $A_\wedge$ requires that a domain be specified. In Section \ref{sec-Domains} we shall construct a natural bijection from the set of domains of closed extensions of $A$ to that of $A_\wedge$, and in \cite{GiKrMe-b} we use pseudodifferential techniques to show that if $\Dom_\wedge$ is the domain of $A_\wedge$ associated with the domain $\Dom$ of $A$, and if $\Lambda$ is a sector of minimal growth for $\csym(A)$ and for $A_{\wedge}$ with domain $\Dom_\wedge$, then it is also a sector of minimal growth for $A$ with domain $\Dom$. This of course brings up the question of how to determine whether a given sector, or even a ray, is a sector of minimal growth for $A_\wedge$ with domain $\Dom_\wedge$. In connection with this we give, in Theorem \ref{NecAndSuffWedge}, a necessary and sufficient condition for $\Lambda$ to be a sector of minimal growth for $A_\wedge$ with domain $\Dom_\wedge$. The condition \eqref{NecAndSuffWedgeCond} of Theorem \ref{NecAndSuffWedge} is in principle verifiable.

Resolvents for cone-elliptic operators written as pseudodifferential operators have been constructed by other authors in special cases, e.g. Br\"uning-Seeley \cite{BrSe88}, Mooers \cite{Mooers}, Gil \cite{GiHeat01}, and Schrohe-Seiler \cite{SchrSe03}, the last mentioned article being the one closest to our own aims in \cite{GiKrMe-b}. Also of interest is Loya \cite{LoRes01} in the context of $b$-operators. Our goal, here and in \cite{GiKrMe-b}, is to study the problem with minimal assumptions.

\medskip
A description of the paper follows.

We shall be working with a fixed elliptic cone operator $A$ acting on sections of a Hermitian vector bundle $E$ over a manifold $M$; the latter is assumed to be compact of dimension $n$ with nonempty boundary. The definition of cone operators is recalled in Section \ref{sec-DefinitionsAndConventions}, where we also recall the definitions of the spaces on which cone operators act. In this section we introduce certain strongly continuous one-parameter groups of isometries $\kappa_\varrho$, one associated with $M$ and one with the interior pointing part of the normal bundle of $\partial M$ in $M$ (where $A_\wedge$ lives). These actions generally play an important role in the analysis of degenerate elliptic operators, see Schulze \cite{SzNorthHolland}, and they do so here as well.

The $c$-cotangent bundle, $\cT^*M$, is defined in Section \ref{sec-TheSymbols}. Its definition is analogous to that of the $b$-tangent bundle of Melrose \cite{RBM1,RBM2}. It is a vector bundle over $M$ which is canonically isomorphic to $T^*M$ over the interior of $M$. Cone operators have invariantly defined symbols, $\csym(A)$, defined on $\cT^*M$. We also recall in this section the definition of $A_\wedge$, and discuss some properties inherited by $A_\wedge$ from $A$. We also briefly recall the definition of the conormal symbol.

In Section \ref{sec-Domains} we first recall known facts about the closed extensions of cone-elliptic operators on compact manifolds, such as $M$, and sketch proofs of analogous results for the operator $A_\wedge$. Proofs are needed since $A_\wedge$, though elliptic in the proper sense, is not a Fredholm operator on the spaces naturally associated with it. For $A$, as is well known, there is a minimal closed extension with domain $\Dom_{\min}$, and there is a maximal extension with domain $\Dom_{\max}$. Likewise, for $A_\wedge$ there is the domain of the minimal extension, $\Dom_{\wedge,\min}$, and the maximal domain $\Dom_{\wedge,\max}$. In both cases, the minimal domain has finite codimension in the maximal domain (in fact the same codimension). The set of domains of closed extensions can be viewed as a Grassmannian variety, and there is a natural map $\Theta$, cf. \eqref{ThetaOperator}, one can use to pass from one variety to the other. This is most relevant in \cite{GiKrMe-b}; indeed, the meaning of the condition that $\csym(A)$ admits a ray of minimal growth is clear, but to express the analogous condition for $A_\wedge$ requires the specification of a domain for $A_\wedge$. This domain is the one associated by $\Theta$ with the given domain for $A$.

The analysis of the spectrum of a given closed extension of $A$ is taken up in Section \ref{sec-DomainsAndSpectra}. It is natural to classify the set of extensions of $A$ by the index. The ones with index $0$ being the only relevant in the problem of studying the spectrum, we let $\mathfrak G$ be the set of domains $\Dom$ such that $\Ind A_\Dom=0$; here and elsewhere $A_\Dom$ means $A$ with domain $\Dom$. The simple condition that both numbers $d''=-\Ind A_{\Dom_{\min}}$ and $d'=\Ind A_{\Dom_{\max}}$ be nonnegative is necessary and sufficient for $\mathfrak G$ to be nonempty, see Lemma \ref{Spec1}, and if this is the case, then $\mathfrak G$ can be viewed as a (complex) Grassmannian variety (based on $\Dom_{\max}/\Dom_{\min})$. An at first surprising fact is that if $\dim\mathfrak G>0$, then for every $\lambda\in\C$ there is $\Dom\in \mathfrak G$ such that $\lambda\in\spec A_{\Dom}$, see Proposition~\ref{Spec4}.

Letting
\begin{equation*}
\bgspec A= \bigcap _{\Dom\in \mathfrak G} \spec A_\Dom,\quad \bgres A=\C\minus \bgspec A,
\end{equation*}
we get
\begin{equation*}
\spec A_\Dom= \bgspec A \cup (\spec A_\Dom \cap \bgres A),
\end{equation*}
a disjoint union. It is the part of $\spec A_\Dom$ in $\bgres A$ that is most amenable to study. For $\lambda\in \bgres A$, the dimension of $\KK_{\lambda}=\ker (A_{\Dom_{\max}} - \lambda)$ is constant, equal to $d'$, and
\begin{equation*}
\lambda\in \resolv A_\Dom  \iff \lambda\in \bgres A \text{ and } \KK_\lambda\cap \Dom = 0,
\end{equation*}
cf. Lemma \ref{ResAndKernelBundle}; by the same lemma, if $\KK_\lambda\cap \Dom = 0$ then $\Dom_{\max}=\KK_\lambda\oplus \Dom$. Let then $\pi_{\KK_\lambda,\Dom}$ be the projection on $\KK_\lambda$ according to this decomposition.

If $\lambda\in \bgres A$, then $A_{\Dom_{\min}}-\lambda$ is injective and $A_{\Dom_{\max}}-\lambda$ is surjective (this property characterizes $\bgres A$). For such $\lambda$ let $B_{\max}(\lambda)$ be the right inverse of $A_{\Dom_{\max}}-\lambda$ whose range is orthogonal to $\KK_\lambda$ with respect to the inner product
\begin{equation*}
(u,v)_A=(Au,Av)_{x^{-m/2}L^2_b}+(u,v)_{x^{-m/2}L^2_b},
\end{equation*}
and let $B_{\min}(\lambda)$ be the left inverse of $A_{\Dom_{\min}}-\lambda$ with kernel the orthogonal of $\rg (A_{\Dom_{\min}}-\lambda)$ in $x^{-m/2}L^2_b(M;E)$. Then, if $\lambda \in \resolv A_\Dom$, one has the formula
\begin{equation*}
B_\Dom(\lambda)=B_{\max}(\lambda) - \big(\id-B_{\min}(\lambda)(A-\lambda)\big) \pi_{\KK_\lambda,\Dom} B_{\max}(\lambda)
\end{equation*}
for the resolvent $B_\Dom(\lambda)=(A_\Dom-\lambda)^{-1}$ of $A_\Dom$, cf. \eqref{ResolventOnDBis}. This formula is evident if one notes that the factor in front of $\pi_{\KK_\lambda,\Dom}$ is the identity on $\KK_\lambda$. In principle both $B_{\min}(\lambda)$ and $B_{\max}(\lambda)$ can be written as pseudodifferential operators, a purely analytic problem, so inverting $A_\Dom-\lambda$ is reduced to an algebraic problem, indeed, a problem in a finite dimensional space, as follows.

Let $\Sing_{\max}$ be the orthogonal of $\Dom_{\min}$ in $\Dom_{\max}$ with respect to the inner product defined above; this is a finite dimensional space. Let $\pi_{\max}:\Dom_{\max}\to \Dom_{\max}$ be the orthogonal projection on $\Sing_{\max}$. Both $I-B_{\min}(\lambda)(A-\lambda)$ and $\pi_{\KK_\lambda,\Dom}$ vanish on $\Dom_{\min}$, so
\begin{equation*}
B_\Dom(\lambda)=B_{\max}(\lambda) -  \big(\id-B_{\min}(\lambda)(A-\lambda)\big) \pi_{\max}\,\pi_{\KK_\lambda,\Dom}\,\pi_{\max} B_{\max}(\lambda).
\end{equation*}
On the other hand,
\begin{equation*}
\lambda\in \resolv A_\Dom \iff \lambda\in \bgres A \text{ and } \pi_{\max}\KK_\lambda\cap \pi_{\max}\Dom= 0,
\end{equation*}
and for such $\lambda$, $\Sing_{\max}=\pi_{\max}\KK_\lambda\oplus \pi_{\max}\Dom$, cf. Lemma \ref{ResAndKernelBundle}. The map $\pi_{\max}\,\pi_{\KK_\lambda,\Dom}\big|_{\Sing_{\max}}$ is just the projection on $\pi_{\max}\KK_\lambda$ according to this decomposition of $\Sing_{\max}$, cf. Lemma~\ref{ReductionOfPiToEmax}.

Organizing the information in terms of Grassmannians turns out to be quite useful. The set $\mathfrak G$ can be viewed as the Grassmannian $\Gr_{d''}(\Sing_{\max})$ of $d''$-dimensional subspaces of $\Sing_{\max}$, and the spaces $\K_\lambda$ (which are the fibers of a holomorphic vector bundle over $\bgres A$) give a holomorphic map $\lambda\mapsto \pi_{\max}\KK_\lambda\in \Gr_{d'}(\Sing_{\max})$. The condition that $\lambda\in \bgres A\cap \spec A_\Dom$ is that $\pi_{\max}\KK_\lambda$ belongs to the set
\begin{equation*}
\mathfrak V_\Dom = \set{\Vee\in \Gr_{d'}(\Sing_{\max}):\Vee\cap \pi_{\max}(\Dom)\ne 0}.
\end{equation*}
This is a complex analytic variety in $\Gr_{d'}(\Sing_{\max})$ of codimension $1$. The condition that for some nonzero $\lambda_0\in \bgres A$, the ray $\set{r\lambda_0:r>R}$ contains no point of $\spec A_\Dom$ is that the curve in $\Gr_{d'}(\Sing_{\max})$ given by $r\mapsto\pi_{\max}\KK_{r\lambda_0}$ has no point in $\mathfrak V_\Dom$ when $r>R$. And if $\Vee\in \Gr_{d'}(\Sing_{\max})\minus \mathfrak V_\Dom$, then the norm of the projection on $\Vee$ using $\Sing_{\max}=\Vee\oplus \pi_{\max}\Dom$ can be estimated in simple terms. This can be useful for estimating the norm of the resolvent of $A_\Dom$ near a point in $\spec A_\Dom\cap \bgres A$.

In Section \ref{sec-SelfAdjointness} we discuss some aspects of symmetric cone operators from the geometric perspective developed in Section \ref{sec-DomainsAndSpectra}. Among other things we show that for such operators, the set of domains of selfadjoint extensions is a real-analytic submanifold of $\mathfrak G$, and that if $\dim \mathfrak G>0$, then for every real $\lambda$ there is a selfadjoint extension of $A$ with $\lambda$ in its spectrum. This is so even if the operator with minimal domain is bounded below (or above). A more detailed study of geometric aspects of the spectrum of selfadjoint extensions will be taken up elsewhere.

In Section \ref{sec-ModelOperator} we analyze $A_\wedge$, also from the perspective of Section \ref{sec-DomainsAndSpectra}. While $A_\wedge$ is not a Fredholm operator, the fact that it is homogeneous under the action of the one-parameter group $\kappa_\varrho$ permits a rather complete analysis of the operator, its background spectrum and the resolvents of the various extensions with index $0$. Theorem \ref{NecAndSuffWedge} gives a necessary and sufficient condition for a given extension of $A_\wedge$ to admit a sector of minimal growth.

We finish the paper proving Theorem \ref{NecAndSuff}, an analogue of Theorem \ref{NecAndSuffWedge} giving a necessary and sufficient condition for an extension of $A$ to admit a sector of minimal growth. While the proofs of these theorems are quite similar, some assumptions in Theorem \ref{NecAndSuff} are automatically satisfied in the case of Theorem \ref{NecAndSuffWedge}.

Most of the nonstandard notation used in this paper, and not mentioned in this introduction, is presented in Sections \ref{sec-Domains} and \ref{sec-DomainsAndSpectra}. In general, objects associated with $A_\wedge$ have the symbol $\wedge$ as part of the notation. For example, $\Sing_{\wedge,\max}$ is the orthogonal of $\Dom_{\wedge,\min}$ in $\Dom_{\wedge,\max}$, and $\pi_{\wedge,\max}$ is the corresponding orthogonal projection. All the other projections will usually indicate the space on which they map: If $H=E\oplus F$, then $\pi_{E,F}:H\to H$ will denote the projection on $E$ according to this decomposition, and $\pi_E$ is the orthogonal projection on $E$.

%%%%%%%%%%%%%%%%%%%%%%%%%%%%%%%%%%%%%%%%%%%%%%%%%%%%%%%%%%%
%%%%%%%%%%%%%%%%%%%%%%%%%%%%%%%%%%%%%%%%%%%%%%%%%%%%%%%%%%%

\section{Definitions and conventions}
\label{sec-DefinitionsAndConventions}

Throughout the paper $M$ is a compact $n$-manifold with boundary, $\m$ is a smooth $b$-measure, $E\to M$ is a Hermitian vector bundle, and $\nabla$ a Hermitian connection on $E$. The boundary of $M$ will be denoted by $Y$. By $x$ we shall mean a smooth defining function of $Y$, positive in the interior $\open M$ of $M$. This function will be fixed shortly so as to have certain properties that simplify the calculations.

The $b$-tangent bundle of Melrose, $\bT M$, is the vector bundle over $M$ whose space of sections is
\begin{equation}\label{GammaBTan}
C^\infty_{\tan}(M;TM)=\set{X\in C^\infty(M,TM): X\text { is tangent to }\partial M},
\end{equation}
see \cite{RBM1,RBM2}. The space $C^\infty_{\tan}(M;\C TM)$ is
a Lie algebra over $\C$ under the usual Lie bracket, and the
collection of elements of order $\leq m$ in its enveloping algebra is
the space $\Diff^m_b(M)$ of totally characteristic differential
operators of order $\leq m$. If $E\to M$ is a complex vector bundle
and $\Diff^m(M;E)$ is the space of differential operators on
$C^\infty(M;E)$ of order $m$, then $\Diff^m_b(M;E)$ denotes the
subspace consisting of totally characteristic differential operators
on $C^\infty(M;E)$ of order $m$, cf. Melrose \cite{RBM2}.

The elements of $x^{-m}\Diff^m_b(M;E)$, that is, differential
operators of the form $A=x^{-m}P$ with $P\in \Diff^m_b(M;E)$, are the
cone operators of order $m$.

The Hilbert space $L^2_b(M;E)$ is the $L^2$ space of
sections of $E$ with respect to the Hermitian form on $E$ and the
density $\m$. Thus the inner product is
\begin{equation*}
(u,v)_{L^2_b}=\int (u,v)_E\, \m
\quad\text{ if }u,\ v \in L^2_b(M;E).
\end{equation*}
For non-negative integers $s$ the Sobolev spaces $H^s_b(M;E)$ are defined as
\begin{equation*}
 H^s_b(M;E)=\set{u\in L^2_b(M;E): Pu\in L^2_b(M;E)\
 \forall P\in \Diff^s_b(M;E)}.
\end{equation*}
The Hilbert space structure is defined using the vector fields in
$C^\infty_{\tan}(M;TM)$ with the aid of the connection on $E$ and a
partition of unity. The spaces $H^s_b(M;E)$ for general $s\in\R$ are
defined by interpolation and duality, and we set
\begin{equation*}
H^\infty_b(M;E)=\bigcap_s H^s_b(M;E),\quad
H^{-\infty}_b(M;E)=\bigcup_s H^s_b(M;E).
\end{equation*}
The weighted spaces
\begin{equation*}
x^\mu H^s_b(M;E)=\set{x^\mu u: u\in H^s_b(M;E)}
\end{equation*}
are Hilbert spaces with the inner product for which $H^s_b(M;E)\ni u\mapsto x^\mu u\in x^\mu
H^s_b(M;E)$ is an isometry. In the case of $s=0$ one has
\begin{equation*}
x^\mu H^0_b(M;E)=x^\mu L^2_b(M;E)=L^2(M,x^{-2\mu}\m;E),
\end{equation*}
and the Sobolev spaces based on $L^2(M,x^{-2\mu}\m;E)$ and $\Diff^s_b(M;E)$ are isomorphic to $x^\mu H^s_b(M;E)$. The topological structure of these spaces is independent of the particular $b$-density on $M$, Hermitian structure and connection of $E$, and defining function $x$.

To simplify a number of computations and constructions it is convenient to introduce additional structure. Let $\pi:NY\to Y$ be the normal bundle of $Y$ in $M$, $NY=T_YM/TY$. Let $x:M\to\R$ be any defining function for $Y$, positive in $\open M$. Since $dx$ vanishes on $TY$, $dx$ defines a function $x_\wedge=dx:NY\to\R$. Define
\begin{equation*}
 Y^\wedge=\set{v\in NY: x_\wedge v \geq 0},
\end{equation*}
and let $\pi_+:Y^\wedge\to Y$ be the restriction of $\pi$.

Let $x\partial_x$ denote the canonical section of $\bT M$ along $Y$.

\begin{lemma}
Let $\m_Y=x\partial_x\contr \m$ be the contraction of $\m$ by $x\partial_x$ along $Y$; $\m_Y$ is a smooth positive density on $Y$. There is a tubular neighborhood map
\begin{equation}\label{TubularNbdMap}
\Phi:V\subset Y^\wedge\to U\subset M
\end{equation}
and a defining function $x$ for $Y$ in $M$ such that \begin{equation}\label{ProductMeasure}
\Phi^*\m=\frac{dx}{x}\otimes \pi_+^*\m_Y\text{ in } V.
\end{equation}
\end{lemma}
\begin{proof} Pick some smooth but otherwise arbitrary tubular neighborhood map $\tilde \Phi$ and a defining function $\tilde x$. Trivialize $N_+Y$ as $[0\neutral],\neutral(\infty)\times Y$ by choosing some smooth vector field $\partial_{\tilde x}$ in $M$ along $Y$ such that $\partial_{\tilde x} \tilde x=1$. Trivialized in this manner, $\tilde x_\wedge: [0\neutral],\neutral(\infty) \times Y\to [0,\neutral(\infty)\neutral]$ is the canonical projection. The $b$-density $\frac{d\tilde x_\wedge}{\tilde x_\wedge}\otimes  \pi^*_+\m_Y$ is smooth, positive, and globally defined on $Y^\wedge$. Therefore, near $\tilde x_\wedge=0$,
\begin{equation*}
\tilde \Phi^*\m=f \frac{d\tilde x_\wedge}{\tilde x_\wedge} \otimes \pi^*_+\m_Y
\end{equation*}
with some smooth function $f$. From the fact that $\tilde \Phi$ is a tubular neighborhood map it follows that $f=1$ when $\tilde x_\wedge=0$. There is $g$ smooth, defined near $\tilde x_\wedge=0$, and equal to $1$ at $\tilde x_\wedge=0$, such that if
\begin{equation*}
F(\tilde x_\wedge,y)=(\tilde x_\wedge g(\tilde x_\wedge, y), y),
\end{equation*}
then
\begin{equation*}
F^*\big(f\frac{d\tilde x_\wedge}{\tilde x_\wedge}\otimes  \pi^*_+\m_Y\big) = \frac{d\tilde x_\wedge}{\tilde x_\wedge}\otimes  \pi^*_+\m_Y.
\end{equation*}
Indeed, this holds if $g$ solves
\begin{equation*}
\partial_{\tilde x_\wedge}g = \frac{g}{f(\tilde x_\wedge g , y)}
\frac{1-f(\tilde x_\wedge g, y)}{\tilde x_\wedge}.
\end{equation*}
Since $f(0,y)=1$, There is a smooth solution with initial condition $g(0,y)=1$. Define $\Phi = \tilde\Phi\circ F$. Then $\Phi$ is a tubular neighborhood map satisfying \eqref{ProductMeasure}. Let $x$ be a smooth function on $M$, positive in $\open M$, that agrees with  $\tilde x_\wedge \circ \tilde \Phi^{-1}$ near $Y$. Then $\Phi$ and $x$ are as required.
\end{proof}

We fix a tubular neighborhood map \eqref{TubularNbdMap} and defining function $x$ for $Y$ such that \eqref{ProductMeasure} holds, and take
\begin{equation}\label{ProductMeasureWedge}
\m_\wedge=\frac{dx_\wedge}{x_\wedge}\otimes\pi_+^*\m_Y
\end{equation}
as density on $Y^\wedge$. We also fix $x_\wedge$ as defining function for $Y$ in $N_+Y$. Both $U$ and $V$ contain $Y$.

Let $X_\wedge=\partial_{x_\wedge}$ be the canonical vertical vector field. Fix a smooth real vector field $X$  on $M$ which coincides with $d\Phi (X_\wedge)$ near $Y$. Shrinking $V$ and $U$ we assume that this holds in $U$.

\begin{definition}\label{ConstantCoeff}
An operator $P\in \Diff^m_b(M;E)$ is said to have coefficients independent of $x$ near $Y$ if $[P,\nabla_{xX}] = 0$ near $Y$.
\end{definition}

The operators on $M$ we are concerned with need not have coefficients independent of $x$. They appear, however, in the form of Taylor coefficients. Namely, if $P\in \Diff^m_b(M;E)$, then for any $N$ there are operators $P_k$, $\tilde P_N\in \Diff^m_b(M;E)$ such that
\begin{equation}\label{bOperatorTaylor}
 P = \sum_{k=0}^{N-1} P_k x^k + \tilde P_N x^N
\end{equation}
where each $P_k$ has coefficients independent of $x$ near $Y$. The operators $P_k$ are uniquely determined near $Y$ by $P$ and our choices of connection on $E$, defining function $x$, and vector field $X$. These Taylor expansions will be used in the course of the construction of the map $\theta$ in Theorem \ref{def:thetaWedge}.

If $P$ has coefficients independent of $x$ near $Y$ then so does its formal adjoint $P^\star$ in $L^2_b(M;E)$. This follows immediately from
\begin{equation}\label{SAixNablaX}
 (\nabla_{xX}u,v)_{L^2_b(M;E)} = -
 (u,\nabla_{xX}v)_{L^2_b(M;E)},\quad u,\ v\in C_0^\infty(\open U;E).
\end{equation}
To see that this formula holds we note that
\begin{equation*}
xX (u,v)_E = (\nabla_{xX}u,v)_E + (u,\nabla_{xX}v)_E
\end{equation*}
because the connection is Hermitian. Near $Y$, the Lie derivative $\mathcal L_{xX}\m$ vanishes because of \eqref{ProductMeasure} and the choice of $X$. So if $u$ and $v$ are supported in $U$ and $h=(u,v)_E$, then $xX h\m=\mathcal L_{xX} h\m=d(h\,xX\contr\m)$. Therefore, by Stokes' theorem, $\int xX\,(u,v)_E\,\m=0$ if $u$, $v\in C_0^\infty(\open U;E)$. This gives \eqref{SAixNablaX}.

Let $E^\wedge\to Y^\wedge$ be the vector bundle $\pi^*_+(E|_Y)$ and give it the canonical Hermitian metric and connection.
An operator $P\in \Diff^m_b(Y^\wedge;E^\wedge)$ is said to have coefficients independent of $x_\wedge$ if it commutes with $\nabla_{x_\wedge X_\wedge}$. The spaces
\begin{equation*}
x^\mu_\wedge H^s_b(Y^\wedge;E^\wedge)
\end{equation*}
are defined in a manner completely analogous to those associated with $M$, using operators with coefficients independent of $x_\wedge$; for nonnegative integers $s$ they may be defined using smooth vector fields in $\bT Y^\wedge$ that commute with $x^\wedge X^\wedge$. Since $Y^\wedge$ is non-compact, $x^\mu_\wedge L^2_b(Y^\wedge;E^\wedge)$ literally means the $L^2$-space corresponding to the measure $x_\wedge^{-2\mu} \m_\wedge$.

Using the tubular neighborhood map $\Phi$, define
\begin{equation*}
 \Phi_*:E^\wedge|_V\to E|_U
\end{equation*}
as follows: For $\eta=(p,\eta')\in E^\wedge_p$ with $p\in V$ and $\eta'\in E_{\pi_+(p)}$, let $\Phi_*\eta\in E_{\Phi(p)}$ be the element obtained by parallel transport of $\eta'$ along the curve $t\mapsto \Phi(tp)$, $t\in [0,1]$. The map $\Phi_*$ is a smooth vector bundle isomorphism covering $\Phi$, an isometry because $\nabla$ is Hermitian. For this reason, and because of \eqref{ProductMeasure}, the induced map
\begin{equation}\label{PhiStar}
\Phi_*:x^{-m/2}_\wedge L^2_b(V;E^\wedge|_V)\to x^{-m/2}L^2_b(U;E|_U)
\end{equation}
is an isometry.

Let $\chi_t$ be the one parameter group of diffeomorphisms of $M$ generated by $xX$. If $u$ is a section of $E$, let $(\kappa_\varrho^{\smash[t]{{}_\parallel}} u)(p)\in E_p$ be the result of parallel transport of $u(\chi_{\log\varrho} p)\in E_{\chi_{\log\varrho} p}$ along the curve
\begin{equation*}
[0,1]\ni s\mapsto \chi_{(1-s)\log\varrho}(p)\in M.
\end{equation*}
There is a unique smooth positive function $f_\varrho: M\to\R$ with the property that
\begin{equation*}
f_\varrho^2 x^m\m = \chi_{\log \varrho}^* (x^m \m).
\end{equation*}

\begin{definition}\label{DilationGroup}
Let $\kappa_\varrho$ act on $C_0^\infty(\open M;E)$ as  $\kappa_\varrho u = f_\varrho \kappa_\varrho^{\smash[t]{{}_\parallel}}$. Denote also by $\kappa_\varrho$ the analogously defined family of maps on $C_0^\infty(\open Y^\wedge;E^\wedge)$ obtained using $\m_\wedge$ and $x_\wedge\partial_{x_\wedge}$.
\end{definition}

The context will indicate whether an instance of $\kappa_\varrho$ means the operator on sections of $E$ over $M$ or sections  of $E^\wedge$ over $Y^\wedge$. In the case of $Y^\wedge$, the function $f_\varrho$ is $\varrho^{m/2}$. Because of the following lemma, the function $f_\varrho$, in the case of $M$, is equal to $\varrho^{m/2}$ near $Y$.

\begin{lemma}
Let $u\in C_0^\infty(V;E^\wedge|_V)$. Then $\kappa_\varrho\Phi_* u=\Phi_*\kappa_\varrho u$ for all $\varrho\geq 1-\eps$ for some $\eps>0$ depending on $u$.
\end{lemma}

This follows from the definitions of $\Phi_*$ and $\kappa_\varrho$, using that near $Y$, $\Phi^*\m=\m_\wedge$, $x^\wedge=x\circ \Phi$, and $\Phi_*\partial_{x_\wedge}=X$. The number $\eps$ serves only to ensure that the support of $\kappa_\varrho u$ is contained in $V$.

\begin{lemma}
The family $\varrho\mapsto \kappa_\varrho$, initially defined on $C_0^\infty(\open M;E)$, extends to $x^{-m/2}L^2_b(M;E)$ as a strongly continuous one-parameter group of isometries.
\end{lemma}

\begin{proof}
Let $h$ denote the Hermitian metric on $E$. If $u$, $v\in C_0^\infty(\open M;E)$, then
\begin{multline*}
h(\kappa_\varrho u,\kappa_\varrho v)x^m \m
= h(f_\varrho \kappa_\varrho^{{}_\parallel} u,f_\varrho \kappa_\varrho^{{}_\parallel} v)x^m \m
\\= h(\kappa_\varrho^{{}_\parallel} u,\kappa_\varrho^{{}_\parallel} v)f_\varrho^2 x^m \m
=\chi_{\log\varrho}^* (h(u, v) x^m \m)
\end{multline*}
so $\kappa_\varrho$ extends to $x^{-m/2}L^2_b(M;E)$ as an isometry. Next we note that
\begin{equation*}
f_{\varrho'\varrho} = f_{\varrho'} \chi_{\log\varrho'}^*f_\varrho.
\end{equation*}
Indeed,
\begin{multline*}
f_{\varrho'\varrho}^2 x^m \m
=\chi_{\log \varrho\varrho'}^*x^m\m
=\chi_{\log\varrho'}^*\chi_{\log\varrho'}^*x^m\m
\\=\chi_{\log\varrho'}^*(f_\varrho^2 x^m\m)
=(\chi_{\log\varrho'}^*f_\varrho^2) f_{\varrho'}^2 x^m\m.
\end{multline*}
Thus
\begin{equation*}
\kappa_{\varrho'\varrho}
= f_{\varrho'\varrho}\kappa_{\smash[t]{\varrho'\varrho}}^{{}_\parallel}
= f_{\varrho'} (\chi_{\log\varrho'}^*f_\varrho) \kappa_{\smash[t]{\varrho'}}^{{}_\parallel} \kappa_{\varrho}^{{}_\parallel}
= f_{\varrho'}  \kappa_{\smash[t]{\varrho'}}^{{}_\parallel} f_\varrho \kappa_{\varrho}^{{}_\parallel}
= \kappa_{\varrho'}\kappa_{\varrho}.
\end{equation*}
That $\varrho\mapsto\kappa_\varrho u$ is continuous follows from the fact that this holds when $u$ belongs to the dense subspace $C_0^\infty(\open M;E)$ of $x^{-m/2}L^2_b(M;E)$ and the continuity of each $\kappa_\varrho$.
\end{proof}

\medskip

We end the section with a brief comment on what we mean by the Mellin transform of an element of $x^{-m/2}L^2_b(M;E)$. Fix $\omega\in C_0^\infty(-1,1)$ real valued, nonnegative and such that $\omega=1$ in a neighborhood of $0$. Also fix a Hermitian connection $\nabla$ on $E$. The Mellin transform of an element $u\in C_0^\infty(\open M;E)$ is defined to be the entire function $\hat u:\C\to C^\infty(Y;E|_Y)$ such that for any $v\in C^\infty(Y;E|_Y)$
\begin{equation*}
 (x^{-\im\sigma}\omega u,\pi_Y^*v)_{L^2_b(M;E)}=
 (\hat u(\sigma),v)_{L^2(Y;E|_Y)}.
\end{equation*}
By $\pi_Y^*v$ we mean the section of $E$ over $U$ obtained by parallel transport of $v$ along the fibers of $\pi_Y$. As is well known, the Mellin transform extends to the spaces $x^\mu L^2_b(M;E)$ in such a way that if $u \in x^\mu L^2_b(M;E)$ then $\hat u(\sigma)$ is holomorphic in $\set {\Im \sigma > -\mu}$ and in $L^2(\set{\Im\sigma=-\mu}\times Y)$ with respect to $d\sigma\otimes \m_Y$.

\medskip

The density $\m$, the map $\Phi$, the function $x$
and the Hermitian connection are fixed throughout the paper. For the sake of  some notational simplification we will henceforth write $x$, $\m$,  and $E$ instead of $x_\wedge$, $\m_\wedge$, and $E^\wedge$. Fixing a defining function $x$ for $Y$ in $M$, as we have done, is equivalent to fixing a trivialization of $Y^\wedge$, a diffeomorphism $Y^\wedge\to [0\neutral],\neutral(\infty)\times Y$.

%%%%%%%%%%%%%%%%%%%%%%%%%%%%%%%%%%%%%%%%%%%%%%%%%%%%%%%%%%%
%%%%%%%%%%%%%%%%%%%%%%%%%%%%%%%%%%%%%%%%%%%%%%%%%%%%%%%%%%%

\section{The symbols of a cone operator}
\label{sec-TheSymbols}

Let $E$, $F\to M$ be complex vector bundles over $M$. An operator
\begin{equation*}
A\in x^{-m}\Diff^m_b(M;E,F)
\end{equation*}
 is called $c$-elliptic if $P=x^mA$ is $b$-elliptic, which means that
its $b$-symbol,
\begin{equation*}
 \bsym(x^mA)\in C^\infty(\bT^*M\minus 0;\mathrm{Hom}(\bpi^*E,\bpi^*F))
\end{equation*}
(cf. Melrose, op. cit.), is invertible. Here $\bpi:\bT^*M\minus 0 \to M$ is the projection map. This definition depends in a mild way on the choice of defining function: if $\tilde x$ is another defining function for $\partial M$, then
\begin{equation}\label{bSymbolsCompared}
\bsym(\tilde x^mA)=(\tilde x/x)^m\,\bsym(x^mA).
\end{equation}

Alternatively, consider the following construction of the $c$-cotangent bundle of $M$,  $\cT^*M$, motivated by Melrose's definition of $\bT M$, and definition of an invariant replacement of the $b$-symbol. Let $\iota:\partial M\to M$ be the inclusion map and define
\begin{equation*}
C^\infty_{\cn}(M;T^*M)=\set{\eta\in C^\infty(M,T^*M): \iota^*\eta=0},
\end{equation*}
the space of smooth $1$-forms on $M$ which are, over $\partial M$, sections of the conormal bundle of $\partial M$ in $M$. Just as with the $b$-tangent bundle, there is the $c$-cotangent bundle, $\cT^* M$, whose space of smooth sections is $C^\infty_{\cn}(M;T^*M)$, and a homomorphism
\begin{equation*}
\cev:\cT^* M\to T^*M
\end{equation*}
which is an isomorphism over the interior. The fiber over $p$ is
\begin{equation*}
 \cT^*_pM = C^\infty_{\cn}(M;T^*M)/\big(\mathcal I_p(M)\cdot
 C^\infty_{\cn}(M;T^*M)\big)
\end{equation*}
where $\mathcal I_p(M)$ is the ideal in $C^\infty(M)$ of functions vanishing at $p$, and the homomorphism $\cev$ is the one induced by
\begin{equation*}
C^\infty_{\cn}(M;T^*M)\ni \eta\mapsto \eta(p)\in T^*_p M.
\end{equation*}
Since the latter map has $\mathcal I_p(M)\cdot C^\infty_{\cn}(M;T^*M)$
in its kernel, it induces a map $\cev_p:\cT^*_pM\to T^*_p M$. Let $\cT
M$ be the dual bundle and let $\cpi:\cT^* M\to M$ be the projection map.

At this point it is convenient to recall that the $b$-tangent bundle of $M$ is defined in a completely analogous manner using $C^\infty_{\tan}(M;TM)$, cf. \eqref{GammaBTan}, so that
\begin{equation*}
 \bT_pM = C^\infty_{\tan}(M;TM)/\big(\mathcal I_p(M)
 \cdot C^\infty_{\tan}(M;TM)\big).
\end{equation*}
Thus we have a map $\bev:\bT M\to TM$.

Now let $A\in x^{-m}\Diff^m_b(M;E,F)$. Since $A$ is a differential operator in the interior of $M$, it has a principal symbol there, given by the standard  formula
\begin{equation*}
 \sym(A)(\xi)(\phi(p))=\lim_{\tau\to\infty} \tau^{-m}e^{-\im \tau f(p)}
 A(e^{\im\tau f}\phi)(p)
\end{equation*}
with $f$ a real-valued smooth function such that $df(p)=\xi$ and with $\phi$ a smooth section of $E$. Suppose now that $f$ is defined in a neighborhood of a point $p_0\in \partial M$ and vanishes on $\partial M$, so that $df$ is conormal to $\partial M$ and therefore represents a local section of $\cT^*M$. If, with local coordinates $x,y_1,\dotsc,y_{n-1}$ and with respect to some frame $\phi_1\dotsc,\phi_r$ of $E$ and frame $\psi_1,\dotsc,\psi_{s}$ of $F$, near $p_0$, we have
\begin{equation*}
A(\sum_\mu h^\mu \phi_\mu)=x ^{-m}\sum_{\mu,\nu}\sum_{k+|\alpha|\leq m} a^\nu_{k\alpha\mu}(x,y)D_y^\alpha(xD_x)^k h^\mu\, \psi_\nu,
\end{equation*}
then, away from the boundary,
\begin{equation*}
\sym(A)(df)(\sum_\mu h^\mu \phi_\mu)=x^{-m}\sum_{\mu,\nu} \sum_{k+|\alpha|=m} a^\nu_{k\alpha\mu}(x,y)(\partial_y f)^\alpha(x\partial_x f)^k h^\mu \,\psi_\nu
\end{equation*}
where by $\partial_y f$ we mean the gradient of $f$ in the $y$ variables. Since $f=xg$ with smooth $g$, this is equal to
\begin{equation*}
\sum_{\mu,\nu}\sum_{k+|\alpha|=m} a^\nu_{k\alpha\mu}(x,y)(\partial_y g)^\alpha(g+x\partial_x g)^k h^\mu \,\psi_\nu,
\end{equation*}
which is smooth up to the boundary. Suppose that $\tilde f$ is another smooth function defined near $p_0$ and vanishing on the boundary, so that $\tilde f=x\tilde g$ for some $\tilde g$. Then the statement that
\begin{equation*}
df-d\tilde f\in \mathcal I_{p_0}(M)\cdot C^\infty_{\cn}(M;T^*M)
\end{equation*}
is equivalent to the statement that
\begin{equation*}
g(p_0)=\tilde g(p_0)\text{ and }\partial_{y_j} g(p_0)=\partial_{y_j} \tilde g(p_0)\text{ for }j=1,\dotsc,n-1;
\end{equation*}
recall that $p_0\in \partial M$. Thus if $df$ and $d\tilde f$ represent the same element of $\cT_{p_0}M$, then
\begin{equation*}
\lim_{p\to p_0} \sym(A)(df(p))(\phi(p))=\lim_{p\to p_0} \sym(A)(d\tilde f(p))(\phi(p))
\end{equation*}
for any smooth section $\phi$ of $E$ defined near $p_0$. It follows that the function
\begin{equation*}
\cT^*_{\,\open M}M\ni \eta \mapsto \csym(A)(\eta)=\sym(A)(\cev(\eta))
\end{equation*}
extends by continuity to a function
\begin{equation*}
\cT^*M\ni \eta \mapsto \csym(A)(\eta),
\end{equation*}
a section of $\mathrm{Hom}(\cpi^*E,\cpi^*F)$ over $\cT^*M\minus 0$. It is easy to see that $\csym(A)$ is smooth.
\begin{definition}
The section $\csym(A)$ is the $c$-symbol of $A$.
\end{definition}
For example, with the notation above, taking $g=\xi +\sum \gamma_j y_j$ with $\xi$ and $\gamma_j$ real constants, and $f=xg$ we get
\begin{equation*}
\sym(A)(df)(\sum_\mu h^\mu \phi_\mu) = \sum_{\mu,\nu}\sum_{k+|\alpha|=m} a^\nu_{k\alpha\mu}(x,y)\gamma^\alpha\xi^k h^\mu \,\psi_\nu,
\end{equation*}
so if $\eta$ is the element of $\cT^*M$ represented by $df$, then the right hand side of this formula is $\csym(A)(\eta)(\sum_\mu h^\mu\phi_\mu)$.

By the definition,
\begin{equation*}
\csym(A)(\eta)=\sym(A)(\cev(\eta)), \quad \eta\in \cT^*_{\,\open M}M.
\end{equation*}
From the fact that $\cev$ is an isomorphism over the interior of $M$, invertibility of $\csym(A)$ over $\open M$ is equivalent to ellipticity of $A$ in that set.

To relate the $c$-symbol of $A$ and the $b$-symbol of $x^mA$ recall first that if $P\in \Diff^m_b(M;E,F)$, then
\begin{equation*}
\bsym(P)(\bev^*\eta)=\sym(P)(\eta), \quad\eta\in T^*M\minus 0;
\end{equation*}
here $\bev^*:T^*M\to \bT^*M$ is the dual of $\bev:\bT M\to TM$. Thus, if $\eta\in \bT M$ projects on an interior point of $M$, then
\begin{equation*}
\bsym(x^mA)(\tilde \eta)=\sym(P)((\bev^*)^{-1}(\tilde\eta)), \quad\tilde \eta\in \bT^*M\minus 0.
\end{equation*}
The fact that $x^mA$ is totally characteristic implies that $\tilde \eta\mapsto \sym(P)((\bev^*)^{-1}(\tilde\eta))$ extends by continuity to the boundary. Let $\eta\in \cT^* M$ project over an interior point. Then
\begin{multline*}
\csym(A)(\eta) = \sym(A)(\cev(\eta)) = x^{-m}\sym(x^mA)(\cev(\eta)) \\= \sym(x^mA)(x^{-1}\cev(\eta))=\bsym(x^mA)(\bev^*(x^{-1}\cev(\eta))).
\end{multline*}
Writing the map $\eta\mapsto \bev^*(x^{-1}\cev(\eta))$ in coordinates one sees that it extends as a smooth isomorphism $\mathbf x^{-1}:\cT^*M\to\bT^*M$, so
\begin{equation*}
\csym(A)(\eta) =\bsym(x^mA)(\mathbf x^{-1}(\eta)).
\end{equation*}
In particular, invertibility of the $c$-symbol of $A$ is equivalent to invertibility of the $b$-symbol of $x^mA$.

The isomorphism $\mathbf x^{-1}:\cT^*M\to\bT^*M$ is determined by the defining function $x$, so is not natural. Write $\mathbf x$ for its inverse. If $\tilde x$ is another defining function for $\partial M$ then $\mathbf x^{-1}\mathbf {\tilde x}$ is multiplication by $\tilde x/x$; this is the reason for \eqref{bSymbolsCompared}.

\begin{definition}\label{ConeParameterElliptic}
The operator $A \in x^{-m}\Diff^m_b(M;E,F)$ is called $c$-elliptic if
\begin{equation*}
\csym(A) \in C^\infty(\cT^*M\minus 0;\mathrm{Hom}(\cpi^*E,\cpi^*F))
\end{equation*}
is an isomorphism. If $F=E$, the family $\lambda\mapsto A-\lambda$ is called $c$-elliptic with parameter in a set $\Lambda\subset \C$ if
\begin{equation*}
\csym(A)-\lambda \in C^\infty((\cT^*M\times \Lambda)\minus 0;\mathrm{End}((\cpi\times\mathit{id})^*E))
\end{equation*}
is an isomorphism. Here $\cpi\times\mathit{id}:(\cT^*M\times\Lambda)\minus 0\to M\times\Lambda$ is the canonical map.
\end{definition}

Let $\chi_t$ be the one-parameter group of diffeomorphisms generated by the vector field $xX$, cf. Section \ref{sec-DefinitionsAndConventions}. Fix $t$ and let $\eta\in C^\infty_\cn(M;T^*M)$. Then $\chi_t^*\eta\in C^\infty_\cn(M;T^*M)$, since $\chi_t\circ \iota=\iota$. Since also $\chi_t^*\mathcal I_{\chi_t(p)}(M)=\mathcal I_p(M)$, we get a map
\begin{equation}
\chi_t^*:\cT M\to \cT M,
\end{equation}
a vector bundle morphism covering $\chi_{-t}$. It is not hard to see that this map is smooth. If $A\in x^{-m}\Diff^m_b(M;E,F)$, let $A_\varrho=\varrho^{-m}\kappa_\varrho^{-1}A\kappa_\varrho$.
Then $A_\varrho\in x^{-m}\Diff^m_b(M;E)$ and
\begin{equation}\label{csymArho}
\csym(A_\varrho) = \varrho^{-m}\kappa_\varrho^{-1}(\csym(A)\circ \chi_{\log\varrho}^*)\kappa_\varrho.
\end{equation}
Thus $A_\varrho$ is $c$-elliptic if $A$ is.

\medskip
We now recall the definitions of conormal and wedge symbols, and of boundary spectrum.

If $P\in \Diff^m_b(M;E)$ and if $u$ is a smooth section of $E$ that vanishes on $Y=\partial M$, then $Pu$ also vanishes on $Y$. Therefore, if $v$ is a section of $E$ over $Y$ and $u$ is an extension of $v$, then $(Pu)|_Y$ does not depend on the extension. Thus, associated with $P$ there is a differential operator
\begin{equation*}
\hat P(0):C^\infty(Y;E|_Y)\to C^\infty(Y;E|_Y)
\end{equation*}
of order $m$. Fix $\sigma\in \C$. Since $u\mapsto x^{-\im\sigma}P (x^{\im\sigma}u)$ is an operator in $\Diff^m_b(M;E)$, there is, for each $\sigma\in \C$, an operator $\hat P(\sigma)\in \Diff^m(Y;E|_Y)$. The conormal symbol of $P$ is defined to be the operator-valued polynomial
\begin{equation}\label{ConormalSymbol}
\C\ni \sigma\mapsto \hat P(\sigma) \in \Diff^m(Y;E|_Y).
\end{equation}
It is easy to verify that $\hat P(\sigma)$ is elliptic for every $\sigma$ if $P$ is $b$-elliptic. The boundary spectrum of $P$, cf. Melrose \cite{RBM2}, Melrose-Mendoza \cite{MM}, is
\begin{equation*}
\spec_b(P)=\set{\sigma\in \C: \hat P(\sigma)\text{ is not invertible}}.
\end{equation*}

The definition of $\hat P(\sigma)$ depends on the choice of defining function $x$ but different choices of defining functions give operators related by conjugation with multiplication by $e^{\im\sigma g}$ for some smooth real-valued function $g$, so the particular choice of defining function to define the conormal symbol is not critical. The conormal symbol of $A\in x^{-m}\Diff^m_b(M;E)$ is defined to be that of the totally characteristic operator $x^mA$, and the boundary spectrum of $A$ is defined to be that of $x^mA$.

If $A\in x^{-m}\Diff^m_b(M;E,F)$, then $\Phi_*^{-1}A\Phi_*$ is a cone operator defined in $V$, cf. \eqref{PhiStar}, and if $u\in C_0^\infty(\open Y^\wedge;E^\wedge)$, then the limit in the following definition exists in $C_0^\infty(\open Y^\wedge;F^\wedge)$.
\begin{definition}
The wedge symbol of $A$ is the operator $A_\wedge\in x^{-m}\Diff^m_b(Y^\wedge;E,F)$ defined by
\begin{equation}\label{InvariantAWedge}
A_\wedge u = \lim_{\varrho\to 0}
\varrho^m\kappa_{\varrho}(\Phi_*^{-1}A\Phi_*)\kappa_\varrho^{-1} u
\end{equation}
\end{definition}
Writing $A$ as $x^{-m}P$ and expanding $P$ as in \eqref{bOperatorTaylor} with $N=1$ we get directly from the definition of $A_\wedge$ that
\begin{equation*}
A_\wedge = x^{-m}P_{\wedge,0}
\end{equation*}
where $P_{\wedge,0}$ is the operator on $Y^\wedge$ that coincides with $\Phi_*^{-1}P_0\Phi_*$ near $Y$ and satisfies 
\begin{equation*}
\kappa_{\varrho}(\Phi_*^{-1}P_{\wedge,0}\Phi_*)\kappa_\varrho^{-1}=P_{\wedge,0}
\end{equation*}
for all $\varrho>0$. Near  $Y$ (in $M$), with the notation and conventions of Section \ref{sec-DefinitionsAndConventions}, $P_0=\sum_{|\alpha|+k\leq m} a_{k \alpha}(y)D_y^\alpha (xD_x)^k$ near $P_0$ where $a_{k \alpha}(y)$ is independent of $x$. The operator $P_{\wedge,0}$ is the ``same'' operator (the pull-back by $\Phi$), but on $Y^\wedge$, and so
\begin{equation*}
A_\wedge = x^{-m}\sum_{|\alpha|+k\leq m} a_{k \alpha}(y)D_y^\alpha (xD_x)^k
\end{equation*}
for all $x$. The conormal symbols of $A$ and $A_\wedge$ are equal to each other, in coordinates the family
\begin{equation*}
x^{-m}\sum_{|\alpha|+k\leq m} a_{k \alpha}(y)\sigma^k D_y^\alpha ,\quad \sigma\in \C.
\end{equation*}

The wedge symbol $A_\wedge$ of $A$ inherits properties of $A$.  Using the tubular neighborhood map $\Phi$ we also get a bundle isomorphism
\begin{equation*}
\cT^*_V Y^\wedge\to\cT^*_{U}M
\end{equation*}
covering $\Phi$. It is not hard to verify that, over $\partial Y^\wedge=Y=\partial M$ we have
\begin{equation*}
\csym(A_\wedge)|_{\cT^*_Y Y^\wedge}=\csym(A)|_{\cT^*_Y M}.
\end{equation*}
Thus $c$-ellipticity is preserved.

Let $A^\star$ be the formal adjoint of $A$ acting on $x^{-m/2}L^2_b(M;E)$. Since $\Phi_*$ and $\kappa_\varrho$ are isometries (the former near $Y$),
\begin{equation*}
(\varrho^m\kappa_{\varrho}(\Phi_*^{-1}A\Phi_*)\kappa_\varrho^{-1} u,v)_{x^{-m/2}L^2_b}=
(u,\varrho^m\kappa_{\varrho}(\Phi_*^{-1}A^\star\Phi_*)\kappa_\varrho^{-1} v)_{x^{-m/2}L^2_b}
\end{equation*}
if $u,v\in C_0^\infty(\open Y^\wedge;E^\wedge)$ and $\varrho$ is small. Thus, taking the limit as $\varrho\to 0$ we get
\begin{equation}\label{AwedgeStarIsAstarWedge}
(A_\wedge)^\star=(A^\star)_\wedge.
\end{equation}
\begin{lemma}
Suppose that $A$ is symmetric on $C_0^\infty(\open M;E)$. Then $A_\wedge$ is symmetric on $C_0^\infty(\open Y^\wedge;E)$. If in addition $A$ is  semibounded from below, then $A_\wedge$ is semibounded from below by $0$.
\end{lemma}
The first assertion follows immediately from \eqref{AwedgeStarIsAstarWedge} and the hypothesis that $A^\star=A$. For the second, let $C\in \R$ be such that
\begin{equation*}
(A u,u)_{x^{-m/2}L^2_b}\geq C\|u\|_{x^{-m/2}L^2_b}^2,\quad u \in C_0^\infty(\open M;E).
\end{equation*}
Suppose that $u \in C_0^\infty(\open Y;E^\wedge)$. Then
\begin{multline*}
(\varrho^m \kappa_\varrho^{-1} (\Phi_*^{-1}A\Phi_*) \kappa_\varrho u,u)_{x^{-m/2}L^2_b}
= \varrho^m( A\Phi_* \kappa_\varrho u,\Phi_*\kappa_\varrho u)_{x^{-m/2}L^2_b}\\
\geq C \varrho^m\|\Phi_*\kappa_\varrho u\|_{x^{-m/2}L^2_b}^2
= C\varrho^m\|u\|_{x^{-m/2}L^2_b}^2
\end{multline*}
Passing to the limit as $\varrho\to 0$ we thus get the second assertion of the lemma.

\medskip

It also follows from \eqref{InvariantAWedge} that the family $\lambda\mapsto A_\wedge -\lambda$ satisfies the homogeneity relation
\begin{equation}\label{HomogeneityOfA}
 A_\wedge - \varrho^m\lambda =
 \varrho^m \kappa_\varrho (A_\wedge -\lambda)\kappa_\varrho^{-1}
 \;\text{ for every } \varrho>0.
\end{equation}
\begin{definition}\label{kappaHomogeneous}
A family of operators $A(\lambda)$
acting on a $\kappa$-invariant space of distributions on $Y^\wedge$ will be called $\kappa$-homogeneous of degree $\nu$ if
\begin{equation*}
 A(\varrho^m\lambda) =
 \varrho^\nu \kappa_\varrho A(\lambda)\kappa_\varrho^{-1}
\end{equation*}
for every $\varrho>0$.
\end{definition}

%%%%%%%%%%%%%%%%%%%%%%%%%%%%%%%%%%%%%%%%%%%%%%%%%%%%%%%%%%%
%%%%%%%%%%%%%%%%%%%%%%%%%%%%%%%%%%%%%%%%%%%%%%%%%%%%%%%%%%%
\section{Closed extensions}
\label{sec-Domains}

In this section we recall some known facts about the closed extensions of a $c$-elliptic cone operator $A$ on a compact manifold and, where needed, sketch proofs of analogous results for the closed extensions of its model operator $A_\wedge$. Theorem \ref{def:thetaWedge} gives a natural isomorphism between $\Dom_{\max}(A)/\Dom_{\min}(A)$ and $\Dom_{\max}(A_\wedge)/\Dom_{\min}(A_\wedge)$ that will play an important role in \cite{GiKrMe-b} but not in the remainder of the present paper. The rest of the material in this section will be used at various points in all later sections.

\medskip
Suppose that $A\in x^{-m}\Diff^m_b(M;E)$, fix $\mu\in \R$ and consider $A$ first as an unbounded operator
\begin{equation}\label{AonWeightMu}
A:C_0^\infty(\open M;E)\subset x^\mu L^2_b(M;E)\to x^\mu L^2_b(M;E)
\end{equation}
Write $\Dom_{\min}(A)$ for the domain of the closure of this operator; with this domain, $A$ is referred to as the minimal extension of $A$. The structure of $\Dom_{\min}(A)$ when $A$ is $c$-elliptic was characterized in Gil-Mendoza \cite[Proposition 3.6]{GiMe01}. Define also
\begin{equation*}
\Dom_{\max}(A)=\set {u\in x^\mu L^2_b(M;E): Au\in x^\mu L^2_b(M;E)}
\end{equation*}
The maximal extension of $A$ is
\begin{equation*}
A:\Dom_{\max}(A)\subset x^\mu L^2_b(M;E)\to x^\mu L^2_b(M;E),
\end{equation*}
also a closed operator. The space $\Dom_{\min}(A)$ is a closed subspace of $\Dom_{\max}(A)$ in the graph norm defined by $A$, and all closed extensions of \eqref{AonWeightMu} have as domain a subspace of $\Dom_{\max}(A)$ containing $\Dom_{\min}(A)$.

It is well known, see Lesch \cite{Le97}, that if $A$ is $c$-elliptic, then $A$ with domain $\Dom_{\max}(A)$ is Fredholm, $\Dom_{\min}(A)$ has finite codimension in $\Dom_{\max}(A)$, and if $\Dom$ is a subspace of $\Dom_{\max}(A)$ containing $\Dom_{\min}(A)$, then
\begin{equation}\label{WeakRelSInd}
\Ind A_\Dom =\Ind A_{\Dom_{\min}}+\dim \Dom/\Dom_{\min}.
\end{equation}
Here $A_\Dom$ means the operator
\begin{equation*}
A:\Dom\subset x^\mu L^2_b(M;E)\to x^\mu L^2_b(M;E).
\end{equation*}
The problem we wish to consider is the nature of the spectrum and structure of the resolvent of the closed extensions of \eqref{AonWeightMu} of index zero (if any).

Since multiplication by $x^\nu$ is an isomorphism (in fact an isometry)
\begin{equation*}
x^\nu:x^\mu L^2_b(M;E)\to x^{\mu+\nu} L^2_b(M;E)
\end{equation*}
we may conjugate $A$ with such operators with no essential change of the problem. For convenience we will work with the operator $x^{-\mu-m/2}Ax^{\mu+m/2}$ so as to base all the analysis on $x^{-m/2}L^2_b(M;E)$. Clearly, $x^{-\mu-m/2}Ax^{\mu+m/2}\in x^{-m}\Diff^m_b(M;E)$. Since
\begin{equation*}
\csym(x^{-\mu-m/2}Ax^{\mu+m/2})=\csym(A),
\end{equation*}
$c$-ellipticity is preserved by such conjugations. We thus assume that $\mu=-m/2$. The standing assumption, unless otherwise indicated, will be that $A$ is $c$-elliptic.

We will usually abbreviate $\Dom_{\min}(A)$ to $\Dom_{\min}$ and $\Dom_{\max}(A)$ to $\Dom_{\max}$ when the operator is clear from the context. As already indicated, the operator $A$ with domain $\Dom$ will be denoted by $A_{\Dom}$.

The inner product
\begin{equation}\label{AInner}
(u,v)_A=(u,v)_{x^{-m/2}L^2_b}+(Au,Av)_{x^{-m/2}L^2_b}
\end{equation}
on $\Dom_{\max}$ makes this space into a Hilbert space.

\begin{definition}\label{defEmax}
The orthogonal of $\Dom_{\min}(A)$ in $\Dom_{\max}(A)$ with respect to this inner product will be denoted $\Sing_{\max}(A)$, or $\Sing_{\max}$ if $A$ is clear from the context. We denote by $\pi_{\max}:\Dom_{\max}(A)\to \Dom_{\max}(A)$ the orthogonal projection on $\Sing_{\max}(A)$.
\end{definition}

Since $\Dom_{\min}$ is closed in $\Dom_{\max}$,
\begin{equation*}
\Dom_{\max}=\Dom_{\min} \oplus \Sing_{\max}
\end{equation*}
and since $\Dom_{\min}$ has finite codimension in $\Dom_{\max}$, $\Sing_{\max}$ is a finite-dimensional space.

\begin{lemma}\label{EmaxAsKernel}
The space $\Sing_{\max}(A)$ is equal to $\Dom_{\max}(A)\cap \ker (A^\star A +I)$, where the kernel is computed in the space of extendable distributions.
\end{lemma}

Here $A^\star$ is the formal adjoint of $A$ in the space $x^{-m/2}L^2_b(M;E)$, that is,
\begin{equation*}
(Au,v)=(u,A^\star v) \quad \forall u,\ v\in C_0^\infty(\open M;E).
\end{equation*}
It is immediate from the definitions of minimal and maximal domains that the Hilbert space adjoint of
\begin{equation}\label{AonDmin}
A:\Dom_{\min}(A)\subset x^{-m/2}L^2_b(M;E) \to x^{-m/2}L^2_b(M;E)
\end{equation}
is
\begin{equation*}
A^\star:\Dom_{\max}(A^\star)\subset x^{-m/2}L^2_b(M;E) \to x^{-m/2}L^2_b(M;E).
\end{equation*}

\begin{proof}
We first show that $\Sing_{\max}\subset \Dom_{\max}(A)\cap \ker (A^\star A +I)$. If $u\in \Sing_{\max}$, then $u\in \Dom_{\max}(A)$ and
\begin{equation*}
(Au,Av)_{x^{-m/2}L^2_b} = -(u,v)_{x^{-m/2}L^2_b} \quad \forall v\in \Dom_{\min}.
\end{equation*}
Therefore the map
\begin{equation*}
\Dom_{\min}(A)\ni v\mapsto (Au,Av)_{x^{-m/2}L^2_b} \in \C
\end{equation*}
is continuous in the norm of $x^{-m/2}L^2_b(M;E)$, and consequently, $Au$ belongs to the domain of the Hilbert space adjoint of \eqref{AonDmin}. Thus $u\in \Dom_{\max}(A^\star A)$ and the identity $(u,v)_A=0$, $v\in \Dom_{\min}(A)$ gives
\begin{equation*}
(A^\star A u, v)+(u,v)=0\quad \forall v\in C_0^\infty(\open M;E)
\end{equation*}
which gives that $u\in \ker (A^\star A+\id)$. Thus $\Sing_{\max}\subset \Dom_{\max}(A)\cap \ker (A^\star A +I)$.

To prove the opposite inclusion, suppose that $u\in \Dom_{\max}(A)\cap \ker (A^\star A +I)$. Then $Au\in x^{-m/2}L^2_b(M;E)$ (since $u\in \Dom_{\max}(A)$) and $A^\star (Au)=-u\in x^{-m/2}L^2_b(M;E)$, so $Au\in \Dom_{\max}(A^\star)$. Thus
\begin{equation*}
(A^\star A u, v)=(Au,Av)\quad \forall v\in C_0^\infty(\open M;E)
\end{equation*}
and it follows that $(u,v)_A=0$ for all $v\in C_0^\infty(\open M;E)$. Since the latter space is dense in $\Dom_{\min}$, we get that $u\in \Sing_{\max}$.
\end{proof}

In the course of the proof we also showed:
\begin{lemma}\label{EmaxASubsetAstarA}
$\Sing_{\max}(A)\subset \Dom_{\max}(A^\star A)$.
\end{lemma}

Since $A$ is $c$-elliptic, so are $A^\star$ and $A^\star A+\id$. It follows that the Mellin transform of any $u\in \Sing_{\max}$ is a meromorphic function defined on all of $\C$.

\medskip
We now discuss analogous aspects for the operator $A_\wedge$. The space $\Dom_{\min}(A_\wedge)$ is the domain of the closure of
\begin{equation*}
A_\wedge:C_0^\infty(\open Y^\wedge;E)\subset x^{-m/2}L^2_b(Y^\wedge;E) \to x^{-m/2}L^2_b(Y^\wedge;E)
\end{equation*}
and
\begin{equation*}
\Dom_{\max}(A_\wedge)=\set{u\in  x^{-m/2}L^2_b(Y^\wedge;E):A_\wedge u \in x^{-m/2}L^2_b(Y^\wedge;E)}.
\end{equation*}
Since $A_\wedge$ need not be Fredholm with either of these domains, we discuss these in some detail. We will usually write $\Dom_{\wedge,\min}$ and $\Dom_{\wedge,\max}$ for the minimal and maximal domains of $A_\wedge$.

\begin{lemma}\label{DomAtInfty}
If $u\in \Dom_{\max}(A_\wedge)$, then $(1-\omega)u\in \Dom_{\min}(A_\wedge)$ for every cut-off function $\omega$ with $\omega=1$ near $x=0$.
\end{lemma}

In other words, as far as closed extensions are concerned, there is no essential structure at infinity.

\begin{proof}
Let $j:Y^\wedge\to Y^\wedge$ be the involution $(x,y)\mapsto (\xi,y)=(1/x,y)$. Under this map, $P_0=x^mA_\wedge$ goes over to a certain other totally characteristic $b$-elliptic operator $\check P_0$, $A_\wedge$ goes to $\check A_\wedge=\xi^m \check P_0$, and
\begin{equation*}
j^*:x^\nu H^\mu_b(Y^\wedge,E)\to\xi^{-\nu}H^\mu_b(Y^\wedge;E)
\end{equation*}
is an isomorphism. We'll write $\check u$ for $j^*u$. Since $\check P_0$ is $b$-elliptic, there are properly supported operators $\check Q$ and $\check R$ defined on extendable distributions such that for every $\nu$ the operators
\begin{equation*}
\check Q:\xi^\nu L^2_b(Y^\wedge;E)\to \xi^\nu H^m_b(Y^\wedge;E),\quad
\check R:\xi^\nu L^2_b(Y^\wedge;E)\to \xi^\nu H^\infty_b(Y^\wedge;E)
\end{equation*}
are continuous and
\begin{equation*}
\check Q\check P_0=I-\check R.
\end{equation*}
If $u\in \Dom_{\max}(A_\wedge)$ then $\check u\in \xi^{m/2} L^2_b(Y^\wedge;E)$ and $\xi^m \check P_0\check u=f\in \xi^{m/2} L^2_b(Y^\wedge;E)$. From
\begin{equation*}
\check Q \xi^{-m}\check f = \check Q \xi^{-m} \xi^m \check P_0 \check u=\check u-\check R\check u,
\end{equation*}
we get
\begin{equation*}
\check u=\check Q \xi^{-m}\check f+\check R\check u
\end{equation*}
with $\check Q \xi^{-m}\check f\in \xi^{-m/2}H^m_b(Y^\wedge;E)$ and $\check R\check u\in \xi^{m/2}H^\infty_b(Y^\wedge;E)$. If $\omega$ is as in the statement of the lemma, then $(1-\check \omega)$ is supported near $\xi=0$, so $(1-\check \omega)\check R\check u\in \xi^{m/2}H^m_b(Y^\wedge;E)$. Thus $(1-\check \omega)\check u\in \xi^{m/2}L^2_b(Y^\wedge;E)\cap \xi^{-m/2}H^m_b(Y^\wedge;E)$. Let $\chi\in C_0^\infty(\R)$ be such that $\chi(\xi)=1$ near $0$ and let $\chi_\ell(\xi)=\chi(\ell\xi)$. Define
\begin{equation*}
\check v_\ell=(1-\chi_\ell)(1-\check \omega)\check u.
\end{equation*}
Then $\check v_\ell\in \xi^{m/2}L^2_b(Y^\wedge;E)\cap \xi^{-m/2}H^m_b(Y^\wedge;E)$ and $\check v_\ell\to (1-\check \omega)\check u$ as $\ell\to\infty$, in  $\xi^{m/2}L^2_b(Y^\wedge;E)$ as well as in $\xi^{-m/2}H^m_b(Y^\wedge;E)$. From the latter we get that $\check P_0\check v_\ell$ converges in $\xi^{-m/2}L^2_b(Y^\wedge;E)$ to $\check P_0(1-\check\omega)\check u$ as $\ell\to\infty$, and consequently, that
$\check A_\wedge \check v_\ell$ converges in $\xi^{m/2}L^2_b(Y^\wedge;E)$ to $\check A_\wedge (1-\check \omega)\check u$. This proves that $(1-\omega)u\in \Dom_{\min}(A_\wedge)$, since the $\check v_\ell$ are compactly supported.
\end{proof}

The structure of $\Dom_{\wedge,\min}$ near $Y$ is described in the first two items of the following proposition, which can be proved using the same arguments as in the proof of \cite[Proposition 3.6]{GiMe01}. The third follows from an analysis of Mellin transforms that takes advantage of the fact that the conormal symbols of $A_\wedge$ and $A$ are the same. An explicit, simple but fundamental isomorphism between the spaces $\Dom_{\wedge,\max}/\Dom_{\wedge,\min}$ and $\Dom_{\max}(A)/\Dom_{\min}(A)$ is given in Theorem \ref{def:thetaWedge}.

\begin{proposition}\label{OnWedgeDomains}
Let $A\in x^{-m}\Diff^m_b(M;E)$ be $c$-elliptic. Then
\begin{enumerate}[$(i)$]
\item $\Dom_{\wedge,\min}=\set{u\in \Dom_{\wedge,\max}:\omega u\in x^{m/2-\eps}H^m_b(Y^\wedge;E)\ \forall \eps >0}$.
\item $\Dom_{\wedge,\min} = x^{m/2}H^m_b(Y^\wedge;E)\cap x^{-m/2}L^2_b(Y^\wedge;E)$ if and only if $\spec_b(A)\cap \set{\Im \sigma=-m/2}=\varnothing$.
\item $\dim \Dom_{\wedge,\max}/\Dom_{\wedge,\min} = \dim \Dom_{\max}(A)/\Dom_{\min}(A)$.
\end{enumerate}
\end{proposition}

On $\Dom_{\wedge,\max}$ we take, naturally,
\begin{equation}\label{AWedgeInner}
(u,v)_{A_\wedge}
= (A_\wedge u,A_\wedge v)_{x^{-m/2}L^2_b} + (u, v)_{x^{-m/2}L^2_b},
\end{equation}
as inner product.

\begin{definition}\label{defWedgeEmax}
The orthogonal of $\Dom_{\wedge,\min}$ in $\Dom_{\wedge,\max}$ with respect to this inner product will be denoted $\Sing_{\max}(A_\wedge)$, or $\Sing_{\wedge,\max}$ if $A_\wedge$ is clear from the context. We denote by $\pi_{\wedge,\max}:\Dom_{\max}(A_\wedge)\to \Dom_{\max}(A_\wedge)$ the orthogonal projection on $\Sing_{\max}(A_\wedge)$.
\end{definition}

The proof of Lemma \ref{EmaxAsKernel} gives that $\Sing_{\wedge,\max}=\Dom_{\wedge,\max}\cap \ker (A_\wedge^\star A_\wedge+\id)$.

The following result, although elementary in nature, is of fundamental importance in expressing the relation between the domains of $A$ and the domains of $A_\wedge$. Let
\begin{equation*}
S=\set{\sigma\in \C: -m/2<\Im\sigma<m/2},
\end{equation*}
and for each $\sigma\in S$ let $N(\sigma)$ be the largest integer $N$ such that $\Im\sigma-N > -m/2$. Let $\sigma_j$, $j=1,\dotsc,\nu$ be an enumeration of the elements of $\Sigma = \spec_b(A)\cap S$.

\begin{theorem}\label{def:thetaWedge}
Let $A$ be an arbitrary $c$-elliptic cone differential operator. There are canonical decompositions
\begin{equation}\label{CanonicalDecompositionWedge}
\Sing_{\max}(A)=\bigoplus_{j=1}^\nu \Sing_{\sigma_j}(A),\quad
\Sing_{\max}(A_\wedge)=\bigoplus_{j=1}^\nu \Sing_{\sigma_j}(A_\wedge)
\end{equation}
such that
\begin{enumerate}[$(i)$]
\item if $u\in \Sing_{\sigma_j}(A)$, then $\hat u|_{\{\Im\sigma>-m/2\}}$ has poles at most at $\sigma_j-\im\vartheta$ for $\vartheta = 0,\dots, N(\sigma_j)$;
\item if $u\in \Sing_{\sigma_j}(A_\wedge)$, then $\hat u|_{\{\Im\sigma>-m/2\}}$ has a pole at most at $\sigma_j$;
\item if $u\in \Sing_{\sigma_j}(A)$ or $u\in \Sing_{\sigma_j}(A_\wedge)$ and $\hat u$ is holomorphic at $\sigma_j$, then $u =0$.
\end{enumerate}
There is a natural isomorphism
\begin{equation}\label{thetaWedge}
\theta:\Sing_{\max}(A)\to\Sing_{\max}(A_\wedge)
\end{equation}
 such that for each $j$,
\begin{equation}\label{FineThetaWedge}
\theta\big|_{\Sing_{\sigma_j}(A)}:
\Sing_{\sigma_j}(A)\to \Sing_{\sigma_j}(A_\wedge),
\end{equation}
and for each $j$ and for all $u\in \Sing_{\sigma_j}(A)$,
\begin{equation}\label{DecayThetaWedge}
(\Phi_*^{-1}\omega u-\theta u)\widehat{\,} \quad \text {is holomorphic near }\sigma_j
\end{equation}
where $\omega\in C_0^\infty(U;E)$ is such that $\omega=1$ near $Y$.
\end{theorem}

\begin{proof} For any open set $U\subset \C$ let $\Mero(U)$ be the space of meromorphic functions on $U$ with values in $C^\infty(\partial M;E|_{\partial M})$. For $\sigma_0\in U$ let $\Mero_{\sigma_0}(U)$ be the subspace of $\Mero(U)$ consisting of elements with pole only at $\sigma_0\in U$. Finally let $\Hol(U)$ be the subspace of holomorphic elements. We let
\begin{equation*}
\mathrm s_{\sigma_0} : \Mero_{\sigma_0}(U) \to \Mero_{\sigma_0}(\C)
\end{equation*}
be the map that sends an element in $\Mero_{\sigma_0}(U)$ to its singular part at $\sigma_0$.

If $A=x^{-m}P$ with $P\in \Diff^m_b(M;E)$, then near $Y=\partial M=\partial Y^\wedge$ we have
\begin{equation}\label{TaylorExpOfA}
x^m A= P =\sum_{k=0}^{m-1} P_k x^k + \tilde P_m x^m
\end{equation}
where each $P_k$, $k<m$, has coefficients independent of $x$, cf. Definition \ref{ConstantCoeff}. Then
\begin{equation*}
x^m A_\wedge=\Phi_*^{-1}P_0\Phi_*
\end{equation*}
near $Y$ in $Y^\wedge$. Let $\hat P_k$ be the conormal symbol of $P_k$. The operator $\hat P_0$ is the conormal symbol of both $A$ and $A_\wedge$.

Let $\sigma_0\in \Sigma$, and let $U\subset S$ be a neighborhood of $\sigma_0$ such that $U\cap \spec_b(A)=\set{\sigma_0}$. Then $\hat P_0$ gives an operator
\begin{equation*}
\mathcal P:\Mero(U)/\Hol(U)\to\Mero(U)/\Hol(U)
\end{equation*}
whose kernel is finite-dimensional. Since $\sigma_0$ is the only point of $\spec_b(A)$ in $U$, the elements in the kernel of ${\mathcal P}$ are represented by meromorphic functions on $U$ with pole only at $\sigma_0$. By taking the singular part of such functions we get a space $\hat \Sing_{\sigma_0}(A_\wedge)\subset \Mero_{\sigma_0}(\C)$  with the property that $h\in \Mero_{\sigma_0}(\C)$ and $\hat P_0 h\in \Hol(U)$ imply that
there is a unique element $\psi\in \hat \Sing_{\sigma_0}(A_\wedge)$ such that
$h-\psi\in \Hol(U)$:
\begin{equation}\label{HatSsigmaWedge}
\hat \Sing_{\sigma_0}(A_\wedge)=\set{\mathrm s_{\sigma_0} \big(\hat P_0(\sigma)^{-1}f(\sigma)\big): f\in \Hol(U)}.
\end{equation}

If $\psi\in \hat \Sing_{\sigma_0}(A_\wedge)$, then there is $u\in x^{-m/2}H^\infty_b(M;E)$ supported in $U$ such that $\hat u -\psi$ is holomorphic in $\Im\sigma>-m/2$. Such $u$ belongs to $x^{-m/2}L^2_b(M;E)$, and since $\hat P_0 \hat u=\hat P_0(\hat u-\psi)$ is holomorphic in $\Im\sigma>-m/2$, we get that $\Phi_*^{-1}u\in \Dom_{\max}(A_\wedge)$. If $v\in \Dom_{\max}(A_\wedge)$ also has the property that $\hat v-\psi$ is holomorphic in $\Im\sigma>-m/2$, then $\Phi_*^{-1}u-v\in\Dom_{\min}(A_\wedge)$, and consequently $\pi_{\wedge,\max}\Phi_*^{-1}u =\pi_{\wedge,\max} v$. Thus there is a well defined operator $F_{\wedge,\sigma_0}:\hat \Sing_{\sigma_0}(A_\wedge)\to\Sing_{\max}(A_\wedge)$, characterized by the property that
\begin{equation*}
\psi - [F_{\wedge,\sigma_0}\psi]\widehat{\,}\quad \text{ is holomorphic in }\Im\sigma>-m/2.
\end{equation*}
From this property one obtains that
\begin{equation*}
\mathrm s_{\sigma_0}([F_{\wedge,\sigma_0}\psi]\widehat{\,})=\psi
\end{equation*}
so the operator $F_{\wedge,\sigma_0}$ is an isomorphism onto its image. Define
\begin{equation}\label{SsigmaWedge}
\Sing_{\sigma_0}(A_\wedge)=F_{\wedge,\sigma_0}\hat \Sing_{\sigma_0}(A_\wedge).
\end{equation}
Clearly, if $\sigma_i$, $\sigma_j\in \Sigma$ and $\sigma_i\ne \sigma_j$, then $\Sing_{\sigma_i}(A_\wedge)\cap \Sing_{\sigma_j}(A_\wedge)=0$.

If $u\in \Sing_{\max}(A_\wedge)$, then $\hat u$ is meromorphic in $\Im\sigma>-m/2$ with poles in $\Sigma$, since $\hat P_0\hat u$ is holomorphic in $\Im\sigma>-m/2$. Therefore $\mathrm s_{\sigma_j}(\hat u)\in \hat \Sing_{\sigma_j}(A_\wedge)$ and $\hat u-\sum\limits_{\sigma_j \in \Sigma} \mathrm s_{\sigma_j}(\hat u)$ is holomorphic in $\Im \sigma>-m/2$. Thus the Mellin transform of
\begin{equation*}
v=u-\sum_{j=1}^\nu F_{\wedge,\sigma_j}\mathrm s_{\sigma_j}(\hat u)
\end{equation*}
is holomorphic in $\Im\sigma>-m/2$, and therefore $v\in \Dom_{\min}(A_\wedge)$. But since $v$ also belongs to $\Sing_{\max}(A_\wedge)$, $v=0$. Thus we have \eqref{CanonicalDecompositionWedge} for the operator $A_\wedge$.

We now construct the spaces $\Sing_{\sigma_j}(A)$ for $A$. Pick $\sigma_0\in \Sigma$ and let $\psi\in \hat \Sing_{\sigma_0}(A_\wedge)$. Thus
$\hat P_0\psi$ is entire. Define $\hat{\e}_{\sigma_0,0}$ as the identity map on $\hat \Sing_{\sigma_0}(A_\wedge)$ and inductively define
\begin{equation*}
\hat{\e}_{\sigma_0,\vartheta}:\hat \Sing_{\sigma_0}(A_\wedge)\to \Mero_{\sigma_0-\im\vartheta}(\C),\quad \vartheta=1,\dotsc,N(\sigma_0)
\end{equation*}
by
\begin{equation*}
\hat{\e}_{\sigma_0,\vartheta}(\psi)= -\mathrm s_{\sigma_0-\im\vartheta}\bigg(\hat P_0(\sigma)^{-1}\sum_{\ell=0}^{\vartheta-1} \hat P_{\vartheta-\ell}(\sigma) \big(\hat{\e}_{\sigma_0,\ell} (\psi) (\sigma + \im(\vartheta-\ell))\big)\bigg).
\end{equation*}
Then
\begin{equation}\label{AtPole}
\sum_{\ell=0}^{\vartheta} \hat P_{\vartheta-\ell}(\sigma) \big(\hat{\e}_{\sigma_0,\ell}(\psi)(\sigma+\im(\vartheta-\ell))\big)
\end{equation}
is entire, and
\begin{equation*}
\sum_{k=0}^{m-1} \hat P_k(\sigma) \sum_{\vartheta=0}^{N(\sigma_0)} \hat{\e}_{\sigma_0,\vartheta}(\psi)(\sigma+\im k)
\end{equation*}
is holomorphic in $\Im\sigma>-m/2$. Define
\begin{equation}\label{HatSsigma}
\hat \Sing_{\sigma_0}(A)=\set{\sum_{\vartheta=0}^{N(\sigma_0)} {\hat{\e}}_{\sigma_0,\vartheta}(\psi): \psi\in \hat \Sing_{\sigma_0}(A_\wedge)}.
\end{equation}
Given $\psi\in \hat\Sing_{\sigma_0}(A_\wedge)$, choose for each $\vartheta$ and element $u_\vartheta\in x^{-m/2}H^\infty_b(M;E)$ such that $\hat u_\vartheta - \hat{\e}_{\sigma_0,\vartheta}(\psi)$ is entire. Then $u=\sum u_{\vartheta} \in \Dom_{\max}(A)$.  If the $v_\vartheta\in x^{-m/2}H^\infty_b(M;E)$ also satisfy the condition that $\hat v_\vartheta - \hat{\e}_{\sigma_0,\vartheta}(\psi)$ is entire, then $\pi_{\max}(u)=\pi_{\max}(v)$, so again we have a well defined operator $F_{\sigma_0}:\hat \Sing_{\sigma_0}(A_\wedge)\to \Sing_{\max}(A)$. This operator is injective; we let
$\Sing_{\sigma_0}(A)$ be its image. It is more tedious than hard to verify that \eqref{CanonicalDecompositionWedge} holds.

Define $\theta$ so that \eqref{FineThetaWedge} holds, and on each
$\Sing_{\sigma_j}(A)$, $\theta=F_{\wedge,\sigma_j}\circ F_{\sigma_j}^{-1}$. Then \eqref{DecayThetaWedge} also holds.
\end{proof}

Let
\begin{equation*}
\mathfrak D(A)= \set{\Dom\subset \Dom_{\max}(A): \Dom \text { is a vector space and } \Dom_{\min}(A)\subset \Dom}.
\end{equation*}
The elements of $\mathfrak D(A)$ are in one-to-one correspondence with the subspaces of $\Sing_{\max}(A)$ via the map
\begin{equation*}
\mathfrak D\ni \Dom\mapsto \Dom\cap \Dom_{\max}(A)=\pi_{\max}(\Dom)\subset \Sing_{\max}(A),
\end{equation*}
so $\mathfrak D(A)$ can be viewed as the union of the Grassmannian varieties of various dimensions associated with $\Sing_{\max}(A)$. Likewise let
\begin{equation*}
\mathfrak D_\wedge= \set{\Dom\subset \Dom_{\wedge,\max}: \Dom \text { is a vector space and } \Dom_{\wedge,\min}\subset \Dom}.
\end{equation*}
With the map $\theta$ of Theorem \ref{def:thetaWedge} we then get a map
\begin{equation}\label{ThetaOperator}
\Theta:\mathfrak D(A)\to\mathfrak D_\wedge.
\end{equation}

%%%%%%%%%%%%%%%%%%%%%%%%%%%%%%%%%%%%%%%%%%%%%%%%%%%%%%%%%%%
%%%%%%%%%%%%%%%%%%%%%%%%%%%%%%%%%%%%%%%%%%%%%%%%%%%%%%%%%%%
\section{Domains and spectra}
\label{sec-DomainsAndSpectra}

We discuss here the spectra and resolvents of the closed extensions of a cone operator $A\in x^{-m}\Diff^m_b(M;E)$ in geometric terms. We continue to assume that $A$ is $c$-elliptic and that $M$ is compact. The results in this section will be relevant mostly in Sections \ref{sec-SelfAdjointness} and \ref{sec-ResolventsOnM}. The conceptual point of view developed here will be taken up in Section \ref{sec-ModelOperator} in the context of the model operator.

We begin with the elementary observation that only those extensions of $A$ that have index zero may have nonempty
resolvent set.

\begin{lemma}\label{Spec1}
There is $\Dom\in \mathfrak D$ such that $\Ind A_{\Dom}=0$ if and only if $\Ind A_{\Dom_{\min}}\leq 0$ and $\Ind A_{\Dom_{\max}}\geq 0$.
\end{lemma}
\begin{proof}
If there is a domain $\Dom\in \mathfrak D$ such that $\Ind
A_{\Dom}=0$, then the relative index formula \eqref{WeakRelSInd} gives
\begin{multline*}
\Ind A_{\Dom_{\min}}\leq\Ind A_{\Dom_{\min}}+\dim\Dom/\Dom_{\max}=0 \\
\leq \Ind A_{\Dom_{\min}}+\dim\Dom_{\max}/\Dom_{\min} =\Ind A_{\Dom_{\max}}.
\end{multline*}
Conversely, suppose that $0\leq -\Ind A_{\Dom_{\min}}$ and $\Ind
A_{\Dom_{\max}}\geq 0$. Using  \eqref{WeakRelSInd} again we get
\begin{equation*}
d=\Ind A_{\Dom_{\max}}-\Ind A_{\Dom_{\min}},
\end{equation*}
so $-\Ind A_{\Dom_{\min}}\leq d$, and there is a subspace of
$\Dom_{\max}/\Dom_{\min}$ of dimension $-\Ind A_{\Dom_{\min}}$.
This subspace corresponds to an element $\Dom\in \mathfrak D$ for
which \eqref{WeakRelSInd} gives
\begin{equation*}
\Ind A_{\Dom}=\Ind A_{\Dom_{\min}}-\Ind A_{\Dom_{\min}}=0,
\end{equation*}
which proves the lemma.
\end{proof}

The domains $\Dom\in \mathfrak D$ on which $A$ has index $0$ are those in
\begin{equation*}
\mathfrak G=\set{\Dom\in \mathfrak D: \dim \Dom/\Dom_{\min}=-\Ind A_{\Dom_{\min}}}.
\end{equation*}
By the lemma, $\mathfrak G$ is empty unless
\begin{equation}\label{ExistenceIndex0}
\Ind A_{\Dom_{\min}}\leq 0\text{ and }\Ind A_{\Dom_{\max}}\geq 0.
\end{equation}
Assuming this, let $d''=-\Ind A_{\Dom_{\min}}$. Then
\begin{equation}\label{GrassmannMap}
\mathfrak G\ni \Dom \mapsto \Dom\cap\Sing_{\max}=\pi_{\max}\Dom \in \Gr_{d''}(\Sing_{\max})
\end{equation}
is a bijection between $\mathfrak G$ and the Grassmannian of $d''$-dimensional subspaces of $\Sing_{\max}$ which we use to give $\mathfrak G$ the structure of a complex manifold. Let $d'=\Ind A_{\Dom_{\max}}$. Then  $d=d'+d''=\dim\Sing_{\max}$.

\medskip

An initial classification of points in the spectrum of a closed extension of $A$ begins with the notion of background spectrum.

\begin{definition}\label{bgSpecAndRes}
The background spectrum of $A$ is the set
\begin{equation*}
\bgspec A=\set{\lambda\in \C:\lambda\in \spec A_\Dom\ \forall \Dom\in \mathfrak D}.
\end{equation*}
The complement of this set, $\bgres A$, is the background resolvent set.
\end{definition}
Thus, if $\Dom\in \mathfrak G$, then $A_\Dom$ has as spectrum the (disjoint) union of $\bgspec A$ and a subset of $\bgres A$. Note that the resolvent set $\resolv A_\Dom$ of $A_\Dom$, $\Dom\in \mathfrak G$, is contained in $\bgres A$. As we shall see, the part of the spectrum of $A_\Dom$ in $\bgres A$ is amenable to detailed study. The set $\bgspec  A$ has the same generic structure as a spectrum:

\begin{lemma}\label{Spec3}
The set $\bgspec  A$ is either $\C$, or closed and discrete.
\end{lemma}

Indeed, $\textup{bg-spec} A$ is an intersection of closed sets, so itself closed, and either all spectra are ${\mathbb C}$ or there is one extension with discrete spectrum.

Thus $\bgres A$ is open. A useful description of $\bgres A$ is as follows.

\begin{lemma}\label{Spec2}
\begin{equation*}
\bgres A=\set{\lambda\in \C:
A_{\Dom_{\min}}-\lambda \text{ is injective and }A_{\Dom_{\max}}-\lambda\text{  is surjective}}.
\end{equation*}
\end{lemma}

\begin{proof}
If $\lambda\in \bgres A$, let $\Dom\in \mathfrak D$ be such that $\lambda\notin\spec A_\Dom$. Since $\Dom_{\min}\subset \Dom$,
$A_{\Dom_{\min}}-\lambda$ is injective, and since $\Dom\subset \Dom_{\max}$,
$A_{\Dom_{\max}}-\lambda$ is surjective. Thus $\lambda\in \bgres A$.

Conversely, suppose that $\lambda$ belongs to the set on the right in the statement of the lemma. Let $R\subset x^{-m/2}L^2_b(M;E)$ be the range of $A_{\Dom_{\min}}-\lambda$, and let $R^\perp$ be its orthogonal. Since $A_{\Dom_{\min}}-\lambda$ is injective, $\dim R^\perp=-\Ind A_{\Dom_{\min}}=d''$. Choose a basis $f_1,\dotsc,f_{d''}$ of $R^\perp$. Since $A_{\Dom_{\max}}-\lambda$ is surjective, we may choose $u_1,\dotsc, u_{d''}\in \Dom_{\max}$ such that $(A -\lambda)u_j=f_j$ for all $j$. The $u_j$ are independent modulo $\Dom_{\min}$. Let
\begin{equation*}
\Dom = \Dom_{\min}\oplus \LinSpan\{u_1,\dots, u_{d''}\}.
\end{equation*}
Then $\Dom\in \mathfrak D$ and $A_\Dom-\lambda$ is invertible, since $R$ is closed.
\end{proof}

\begin{proposition}\label{Spec4}
Suppose that \eqref{ExistenceIndex0} holds and that $\dim\mathfrak G>0$. Then, for every $\lambda\in \C$, there is $\Dom \in \mathfrak G$ such that $\lambda\in \spec A_{\Dom}$.
\end{proposition}

\begin{proof}
Let $\lambda\in \C$. If $\lambda\in \bgspec A$ then in fact $\lambda\in \spec A_\Dom$ for any $\Dom\in \mathfrak G$. Suppose then that $\lambda\notin \bgspec  A$, so by Lemma \ref{Spec2} there is $\Dom_0\in \mathfrak G$ such that $A_{\Dom_0}-\lambda$ is invertible. The hypothesis that $\dim\mathfrak G>0$ is equivalent to the statement that the two numbers $d''=-\Ind A_{\Dom_{\min}}$ and $d'=\Ind A_{\Dom_{\max}}$ are strictly positive; recall that their sum is $d$, the dimension of $\Dom_{\max}/\Dom_{\min}$. Let $w\in \Dom_{\max}\minus\Dom_0$. Such $w$ exists because $d''<d$. Let $f=(A-\lambda)w$, and let $v\in \Dom_0$ be such that $(A-\lambda)v=f$. Then $w-v\ne 0$ modulo $\Dom_{\min}$, and thus is an eigenvector of $A$. Let  $\Dom\in \mathfrak G$ contain $w-v$; such $\Dom$ exists because $d''>0$. Then $A_{\Dom}-\lambda$ has nontrivial kernel.
\end{proof}

We will write $\KK_\lambda$ for the kernel of $A_{\Dom_{\max}}-\lambda$, $\lambda \in \bgres A$. For such $\lambda$,
\begin{equation*}
\dim \KK_\lambda=\Ind A_{\Dom_{\max}},
\end{equation*}
since $A_{\Dom_{\max}}-\lambda$ is surjective and its index is independent of $\lambda$.

\begin{proposition}\label{BundleOfKernels}
Let $\KK=\bigsqcup_{\lambda\in\bgres A} \KK_\lambda$ and let $\rho:\KK\to\bgres A$ be the natural map. Then $\KK\to\bgres A$ is a locally trivial Hermitian holomorphic vector bundle.
\end{proposition}
\begin{proof}
Let $\lambda_0\in \bgres A$, let $\KK_{\lambda_0}^\perp$ be the orthogonal of $\KK_{\lambda_0}$ in $\Dom_{\max}$. The operators
\begin{equation*}
A_{1}(\lambda)=(A-\lambda)|_{\KK_{\lambda_0}} \quad
A_{2}(\lambda)=(A-\lambda)|_{\KK_{\lambda_0}^\perp}
\end{equation*}
are continuous as operators into $x^{-m/2}L^2_b(M;E)$ when the domains are given the graph norm of $A$, and depend holomorphically on $\lambda$. Since $A_{2}(\lambda_0)$ is invertible, the inverse $A_{2}(\lambda)^{-1}$ exists for $\lambda$ close to $\lambda_0$. It is easy to verify that if $u_0\in \KK_{\lambda_0}$, then
\begin{equation*}
u(\lambda)=u_0-A_{2}(\lambda)^{-1}A_{1}(\lambda)u_0 \in \KK_\lambda
\end{equation*}
for $\lambda$ close to $\lambda_0$. These are, by definition, holomorphic local sections of $\KK$. The statement that $\KK\to\bgspec A$ is a locally trivial holomorphic vector bundle follows by taking local frames near $\lambda_0$ of the form $u_j(\lambda)$ where the $u_j$ form a basis of $\KK_{\lambda_0}$. The Hermitian form in $\KK$ is the one whose restriction to $\KK_\lambda$ is the restriction of the inner product of $\Dom_{\max}$ to $\KK_\lambda$.
\end{proof}

Note that if $u$, $v\in \KK_\lambda$, then
\begin{equation}\label{HermitianForm}
(u,v)_A=(Au,Av)+(u,v)=(1+|\lambda|^2)(u,v).
\end{equation}

\begin{lemma}\label{ResAndKernelBundle}
Let $\Dom\in \mathfrak G$. The following are equivalent:
\begin{enumerate}[$(i)$]
\item $\lambda \in \resolv A_\Dom$;
\item $\lambda\in \bgres A$  and $\KK_\lambda\cap \Dom=0$;
\item $\lambda\in \bgres A$  and $\pi_{\max}\KK_\lambda\cap \pi_{\max}\Dom=0$.
\end{enumerate}
Moreover, if $\lambda\in \resolv A_\Dom$, then
\begin{equation}\label{KoplusDom}
\KK_\lambda\oplus\Dom=\Dom_{\max}\text{ and } \pi_{\max}\KK_\lambda\oplus\pi_{\max}\Dom=\Sing_{\max}.
\end{equation}
\end{lemma}
\begin{proof}
To prove the equivalence of (i) and (ii), we recall first that $\resolv A_\Dom\subset \bgres A$. A point $\lambda\in\C$ belongs to $\resolv A_\Dom$ if and only if $\ker (A_\Dom-\lambda)=0$, because $A$ is Fredholm of index $0$. But for $\lambda\in \bgres A$, $\ker (A_\Dom-\lambda)=\KK_\lambda\cap \Dom$. Thus (i) and (ii) are equivalent.

Suppose that $\lambda\in \bgres A$. If $u\in \pi_{\max}\KK_\lambda\cap \pi_{\max}\Dom$, then $u=\phi-v$ with $\phi\in \KK_{\lambda}$ and $v\in \Dom_{\min}$. Thus $\phi=u+v\in \pi_{\max}\Dom+\Dom_{\min}=\Dom$, and so $\phi\in \KK_\lambda\cap \Dom$. If $\pi_{\max}\KK_\lambda\cap \pi_{\max}\Dom\ne0$, pick $u\ne 0$. Then $\phi\ne 0$, so $\KK_\lambda\cap \Dom\ne 0$. Thus (ii) implies (iii).

Again suppose that $\lambda\in \bgres A$. To prove that (iii) implies (ii) we will first show that $\pi_{\max}|_{\KK_\lambda}:\KK_\lambda\to \Sing_{\max}$ is injective. Let $\phi\in \KK_\lambda$. If $\pi_{\max}\phi= 0$ then $\phi\in\Dom_{\min}$. But $A-\lambda$ is injective on $\Dom_{\min}$, so $\phi=0$. Thus if $\KK_\lambda\cap \Dom\ne 0$, then $\pi_{\max}\KK_\lambda\cap \pi_{\max}\Dom\ne 0$.

To prove the last statement we first observe that
\begin{equation*}
\dim \KK_\lambda+\dim \Dom/\Dom_{\min}=\Ind A_{\Dom_{\max}}-\Ind A_{\Dom_{\min}}=\dim \Sing_{\max}.
\end{equation*}
This gives, in view of (iii), that $\pi_{\max}\KK_\lambda\oplus\pi_{\max}\Dom=\Sing_{\max}$. Adding $\Dom_{\min}$ to both sides of this formula gives $\KK_\lambda+\Dom=\Dom_{\max}$, but this sum is direct in view of (ii).
\end{proof}

The lemma gives
\begin{equation}\label{TheSpectrum}
\spec A_\Dom=\bgspec A\cup\set{\lambda\in \bgres A: \KK_\lambda\cap \Dom\ne 0}
\end{equation}
for any $\Dom\in \mathfrak G$. Since $\KK_\lambda\cap \Dom=0$ if and only if $\pi_{\max}\KK_\lambda \cap \pi_{\max}\Dom=0$, the presence of spectrum in $\bgres A$ for a given extension $A_\Dom$ is a purely finite dimensional phenomenon. We will exploit this in Section \ref{sec-ResolventsOnM} to give estimates for the resolvent $B_\Dom(\lambda)$ of $A_\Dom-\lambda$ in terms of a canonical right inverse of $A_{\Dom_{\max}}-\lambda$, a canonical left inverse for $A_{\Dom_{\min}}-\lambda$, $\lambda\in \bgres A$, and a finite dimensional projection.

If $\lambda\in \bgres A$, more generally, if $A_{\Dom_{\max}}-\lambda$ is surjective, then $A_{\Dom_{\max}}-\lambda$ admits a continuous right inverse $B_{\max}(\lambda)$, namely, let $\KK_\lambda^\perp\subset \Dom_{\max}$ be the orthogonal of $\KK_\lambda$ with respect to the inner product \eqref{AInner} ($\KK_\lambda^\perp$ may not be, and does not need to be, an element of $\mathfrak G$). The operator
\begin{equation*}
(A-\lambda)|_{\KK_\lambda^\perp}:\KK_\lambda^\perp\to x^{-m/2}L^2_b(M;E)
\end{equation*}
is continuous and bijective. Then the inverse, $B_{\max}(\lambda)$, of $(A-\lambda)|_{\KK_\lambda^\perp}$ is continuous. For each $\lambda\in \bgres A$, the operator $B_{\max}(\lambda)$ has the smallest norm among all continuous right inverses of $A_{\Dom_{\max}}-\lambda$.

The operators $B_{\max}(\lambda)$ can be obtained from any family
\begin{equation*}
B'_{\max}(\lambda):x^{-m/2}L^2_b(M;E)\to\Dom_{\max}
\end{equation*}
of continuous right inverses for $A_{\Dom_{\max}}-\lambda$ by means of the formula
\begin{equation}\label{CorrectionFormulaMax}
B_{\max}(\lambda)=B'_{\max}(\lambda)-\pi_{\KK_\lambda} B'_{\max}(\lambda)
\end{equation}
in which $\pi_{\KK_\lambda}:\Dom_{\max}\to \KK_\lambda$ is the orthogonal projection on $\KK_\lambda$ (with respect to \eqref{AInner}).

The $B_{\max}(\lambda)$, as operators $x^{-m/2}L^2_b(M;E)\to \Dom_{\max}$, depend continuously, even smoothly, on $\lambda$. To see this, let $\lambda$, $\lambda_0\in \bgres A$. Then
\begin{equation*}
(A-\lambda)B_{\max}(\lambda_0)=\big( (A -\lambda_0)+(\lambda_0-\lambda)\big) B_{\max}(\lambda_0)= \id +(\lambda_0-\lambda)B_{\max}(\lambda_0).
\end{equation*}
Since both $B_{\max}(\lambda_0):x^{-m/2}L^2_b(M;E)\to\Dom_{\max}$ and the inclusion $\iota:\Dom_{\max}\embed x^{-m/2}L^2_b(M;E)$ are continuous,
\begin{equation*}
\iota B_{\max}(\lambda_0):x^{-m/2}L^2_b(M;E) \to x^{-m/2}L^2_b(M;E)
\end{equation*}
is continuous. So if $\lambda$ is close enough to $\lambda_0$, then
\begin{equation*}
B'_{\max}(\lambda) = B_{\max}(\lambda_0)\big(\id +(\lambda_0-\lambda)\iota B_{\max}(\lambda_0)\big)^{-1}
\end{equation*}
is a right inverse for $A_{\Dom_{\max}}-\lambda$ depending smoothly on $\lambda$. Since the $\pi_{\KK_\lambda}$, as operators $\Dom_{\max}\to\Dom_{\max}$, also depend smoothly on $\lambda$, the correction \eqref{CorrectionFormulaMax} gives the smoothness of $\lambda\mapsto B_{\max}(\lambda)$.

The operators $B_{\max}(\lambda)$ can be used to construct the resolvent of $A_\Dom-\lambda$ for any $\Dom\in \mathfrak G$, as follows. For each $\lambda\in \bgres A$ such that $\KK_\lambda\cap \Dom=0$ let
\begin{equation*}
\pi_{\KK_\lambda,\Dom}:\Dom_{\max}\to \KK_\lambda
\end{equation*}
be the projection according to the decomposition \eqref{KoplusDom}; this is a continuous operator. Noting that $\lambda \in \resolv A$ if and only if $\lambda\in \bgres A$ and $\KK_\lambda\cap \Dom=0$, define
\begin{equation}\label{ResolventOnDviaMax}
B_\Dom(\lambda)=B_{\max}(\lambda)-\pi_{\KK_\lambda,\Dom} B_{\max}(\lambda).
\end{equation}
Then $\pi_{\KK_\lambda,\Dom}B_\Dom(\lambda)=0$, so $B_\Dom(\lambda)$ maps into $\Dom$. Since $(A-\lambda) \pi_{\KK_\lambda,\Dom}=0$, $B_\Dom(\lambda)$ is a right inverse for $A_\Dom-\lambda$, which must also be the left inverse because $A_\Dom-\lambda$ is invertible.

The canonical left inverse for $A_{\Dom_{\min}}-\lambda$ is constructed in an analogous manner. Let $\RR_\lambda$ be the range of $A_{\Dom_{\min}}-\lambda$, $\lambda\in \bgres A$ (more generally, one can let $\lambda$ belong to the set where $A_{\Dom_{\min}}-\lambda$ is injective). Since $A_{\Dom_{\min}}-\lambda$ is injective if $\lambda\in \bgres A$, $A_{\Dom_{\min}}-\lambda:\Dom_{\min}\to \RR_\lambda$ has a continuous left inverse $B^0_{\min}(\lambda)$. The orthogonal $\RR_\lambda^\perp$ has dimension $-\Ind A_{\Dom_{\min}}$. Let $B_{\min}(\lambda)$ be the composition of the orthogonal  projection on $\RR_\lambda$ followed by $B^0_{\min}(\lambda)$. Viewing $\RR_{\lambda}^\perp$ as the kernel of $A^\star-\overline \lambda$ on $\Dom_{\max}(A^\star)$, we see that $\bigsqcup_{\lambda \in \bgspec A}\RR_{\lambda}$ is a smooth (anti-holomorphic) vector bundle over $\bgres A$. An analysis similar to that done for $B_{\max}(\lambda)$ gives that $B_{\min}(\lambda)$ depends smoothly on $\lambda\in \bgres A$.

If $B_{\min}'(\lambda)$ is a left inverse for $A_{\Dom_{\min}}-\lambda$, $\lambda\in \bgres A$, and $\pi_{\RR_\lambda}$ is the orthogonal projection on $\RR_\lambda$ (in $x^{-m/2}L^2_b(M;E)$), then
\begin{equation}\label{CorrectionFormulaMin}
B_{\min}(\lambda)=B_{\min}'(\lambda)\pi_{\RR_\lambda},
\end{equation}
and so
\begin{equation*}
\|B_{\min}(\lambda)\|_{\L(x^{-m/2}L^2_b,\Dom_{\max})} \leq \|B_{\min}'(\lambda)\|_{\L(x^{-m/2}L^2_b,\Dom_{\max})}.
\end{equation*}

Let $\Dom\in \mathfrak G$ and let $B_\Dom(\lambda)$ be the resolvent of $A_\Dom-\lambda$. It is immediate that the formula
\begin{equation}\label{ResolventOnDviaMin}
B_\Dom(\lambda)=B_{\min}(\lambda) + \big(\id-B_{\min}(\lambda)(A-\lambda)\big)B_{\Dom}(\lambda)
\end{equation}
holds. Replacing \eqref{ResolventOnDviaMax} in this formula we get
\begin{equation}\label{ResolventOnDviaMinMaxTmp}
B_\Dom(\lambda)=B_{\max}(\lambda) - \big(\id-B_{\min}(\lambda)(A-\lambda)\big)\pi_{\KK_\lambda,\Dom}B_{\max}(\lambda).
\end{equation}
Letting $\pi_{\min}=\id-\pi_{\max}$ we see that that
\begin{equation}\label{PiKLambdaDFactors}
\pi_{\KK_\lambda,\Dom}=\pi_{\KK_\lambda,\Dom}(\pi_{\max}+\pi_{\min})=\pi_{\KK_\lambda,\Dom}\pi_{\max}.
\end{equation}
The operator $\id-B_{\min}(\lambda)(A-\lambda)$ is a projection with kernel $\Dom_{\min}$ so
\begin{equation*}
\id-B_{\min}(\lambda)(A-\lambda) = \big(\id-B_{\min}(\lambda)(A-\lambda)\big)\pi_{\max}.
\end{equation*}
Thus we arrive at
\begin{equation}\label{ResolventOnDBis}
B_\Dom(\lambda) = B_{\max}(\lambda) - \big(\id -B_{\min}(\lambda)(A-\lambda)\big) \pi_{\max}\pi_{\KK_\lambda,\Dom}\pi_{\max}B_{\max}(\lambda),
\end{equation}
a formula which will prove to be very useful.

\begin{remark}
The range of the projector $\id-B_{\min}(\lambda)(A-\lambda)$ contains $\KK_{\lambda}$ so there is no difference between \eqref{ResolventOnDviaMinMaxTmp} and \eqref{ResolventOnDviaMax}. Writing $B_\Dom(\lambda)$ in the form \eqref{ResolventOnDBis} separates the geometric information, in $\pi_{\max}\pi_{\KK_\lambda,\Dom}\pi_{\max}$, from analytic contributions.
\end{remark}

We now focus on  $\pi_{\max} \pi_{\KK_\lambda,\Dom} \pi_{\max}$, in particular its norm as a map $\Sing_{\max}\to\Sing_{\max}$.

\begin{lemma}\label{ReductionOfPiToEmax}
Let $\Dom\in \mathfrak G$. Suppose $\lambda\in \bgres A$ and $\KK_\lambda\cap \Dom=0$. Then $\pi_{\max}\pi_{\KK_\lambda,\Dom}\big|_{\Sing_{\max}}$ is the projection on $\pi_{\max}\KK_\lambda$ according to the decomposition $\Sing_{\max}=\pi_{\max}\KK_{\lambda}\oplus \pi_{\max}\Dom$.
\end{lemma}

\begin{proof}
The map $\pi_{\max}\pi_{\KK_\lambda,\Dom}\big|_{\Sing_{\max}}$ is a projection. Indeed, in view of \eqref{PiKLambdaDFactors},
\begin{equation*}
\pi_{\max}\pi_{\KK_\lambda,\Dom}\pi_{\max}\pi_{\KK_\lambda,\Dom}\big|_{\Sing_{\max}}
=\pi_{\max}\pi_{\KK_\lambda,\Dom}\pi_{\KK_\lambda,\Dom}\big|_{\Sing_{\max}}
=\pi_{\max}\pi_{\KK_\lambda,\Dom}\big|_{\Sing_{\max}},
\end{equation*}
The operator $\pi_{\max}\pi_{\KK_\lambda,\Dom}\big|_{\Sing_{\max}}$
has kernel containing $\pi_{\max}\Dom$, since the latter space is contained in $\Dom$, and range contained in $\pi_{\max}\KK_\lambda$.
To complete the proof we only need to show that $\ker \pi_{\max}\pi_{\KK_\lambda,\Dom}\big|_{\Sing_{\max}}=\pi_{\max}\Dom$. Suppose that $u\in \ker \pi_{\max}\pi_{\KK_\lambda,\Dom}\big|_{\Sing_{\max}}$. Then $\phi=\pi_{\KK_\lambda,\Dom}u\in \KK_\lambda$ has the property that $\pi_{\max}\phi=0$. Thus $\phi\in \Dom_{\min}$. But since $A-\lambda$ is injective on $\Dom_{\min}$ (since $\lambda\in \bgres A$), $\phi=0$. That is, $u\in \ker \pi_{\KK_\lambda,\Dom}$. Since $u$ is already in $\Sing_{\max}$, this gives $u\in \Dom\cap \Sing_{\max}$. But the latter space is $\pi_{\max}\Dom$.
\end{proof}

In the course of the proof of Lemma \ref{ResAndKernelBundle} we showed that if $\lambda\in \bgres A$, then $\pi_{\max}|_{\KK_\lambda}:\KK_\lambda\to \Sing_{\max}$ is injective.  Thus, since the spaces $\KK_\lambda$ have dimension $d'=\Ind A_{\Dom_{\max}}$, we have a map
\begin{equation*}
\bgres A\ni \lambda \mapsto \pi_{\max}\KK_\lambda \in \Gr_{d'}(\Sing_{\max}).
\end{equation*}
Write $\KK_{\max}$ for this map, so $\KK_{\max}(\lambda)=\pi_{\max}\KK_\lambda$. If $\lambda_0\in \bgres A$, let $\phi_1,\dotsc,\phi_{d'}$ be a holomorphic frame of $\KK$, cf. Proposition \ref{BundleOfKernels}, near $\lambda_0$. Thus, in addition to independence, the maps $\lambda\mapsto \phi_j(\lambda)$ are holomorphic for $\lambda$ near $\lambda_0$. If $u_1,\dotsc, u_d$ is an orthonormal basis of $\Sing_{\max}$, then $\KK_{\max}(\lambda)$ is spanned by the vectors
\begin{equation*}
\pi_{\max}\phi_j(\lambda)=\sum_k (\phi_j(\lambda),u_k)u_k,
\end{equation*}
which depend holomorphically on $\lambda$. Thus $\KK_{\max}:\bgres A\to \Gr_{d'}(\Sing_{\max})$ is holomorphic.

If $\Dom\in \mathfrak G$, then Lemma \ref{ResAndKernelBundle} asserts that $\lambda\in \bgres A\cap \spec A_\Dom$ if and only if $\pi_{\max}\KK_\lambda\cap \pi_{\max}\Dom\ne 0$. Writing $\W=\pi_{\max}\Dom$, let
\begin{equation*}
\mathfrak V_\W=\set{\Vee\in\Gr_{d'}(\Sing_{\max}): \Vee\cap \W\ne 0}.
\end{equation*}
Then
\begin{equation*}
\lambda\in \bgres A\cap \spec A_\Dom\iff \KK_{\max}(\lambda)\in \mathfrak V_\W.
\end{equation*}

\begin{definition}\label{BadVariety}
For any nonnegative integer $d_0\leq d$ and $\W\in \Gr_{d_0}(\Sing_{\max})$ let
\begin{equation*}
\mathfrak V_\W=\set{\Vee\in\Gr_{d-d_0}(\Sing_{\max}): \Vee\cap \W\ne 0}.
\end{equation*}
If $\Dom\in \mathfrak D$, we write $\mathfrak V_{\Dom}$ for $\mathfrak V_{\pi_{\max}\Dom}$.
\end{definition}

Thus $\spec A_\Dom$ is the union of $\bgspec A$ and the pre-image of $\mathfrak V_{\Dom}$ under the map $\KK_{\max}$.

\begin{proposition}
The set $\mathfrak V_\W\subset \Gr_{d'}(\Sing_{\max})$ is a variety of (complex) codimension $1$. For each $\Dom\in \mathfrak G$,
\begin{equation*}
\spec A_{\Dom}=\bgspec A\cup \KK_{\max}^{-1}(\mathfrak V_\Dom).
\end{equation*}
This is a disjoint union.
\end{proposition}

\begin{proof}
We already showed the second statement. To prove the first statement, fix an ordered basis $\mathbf u=[u_1,\dotsc,u_d]$ for $\Sing_{\max}$. Pick some point  $\Vee_0\in \Gr_{d'}(\Sing_{\max})$ and let $\mathbf \Phi=[\phi_1,\dotsc,\phi_{d'}]$ be a holomorphic local section defined near $\Vee_0$ of the bundle of ordered bases of the canonical bundle over $\Gr_{d'}(\Sing_{\max})$. Thus
\begin{equation*}
\mathbf \Phi(\Vee)=\mathbf u\cdot Z(\Vee)
\end{equation*}
for some matrix $Z(\Vee)\in M^{d\times d'}(\C)$ depending holomorphically on $\Vee$. Let $\mathbf \Psi=[\psi_1,\dotsc,\psi_{d''}]$ be a basis of $\W$, so $\mathbf \Psi=\mathbf u\cdot W$ with $W\in M^{d\times d''}(\C)$. Then
\begin{equation*}
f(\Vee)=\det[Z(\Vee)|W]
\end{equation*}
is holomorphic in $\Vee$. Since $[\mathbf \Phi(\Vee),\mathbf \Psi]$ fails to be a basis of $\Sing_{\max}$ if and only if $\Vee\cap \W\ne 0$, $f(\Vee)$ vanishes if and only if $\Vee\cap \W\ne 0$. Thus $\mathfrak V_\W$ is a complex variety of codimension $1$.
\end{proof}

The norm of the factor $\pi_{\max}\pi_{\KK_\lambda,\Dom}\pi_{\max}$ in \eqref{ResolventOnDBis}, defined for $\lambda\in\bgres A\minus\spec A$, can be estimated in quite simple terms. Using Lemma \ref{ReductionOfPiToEmax}, the problem is generally to estimate, for any $\W\in \Gr_{d_0}(\Sing_{\max})$ and $\Vee\in \Gr_{d-d_0}(\Sing_{\max})\minus \mathfrak V_\W$, the norm of the projection
\begin{equation}\label{GeneralProjectionEmax}
\pi_{\Vee,\W}:\Sing_{\max}\to \Sing_{\max}
\end{equation}
on $\Vee$ according to the decomposition $\Sing_{\max}=\Vee\oplus \W$. We assume that the integer $d_0$ satisfies $0< d_0<d$.

Let then $\W\in \Gr_{d_0}(\Sing_{\max})$. Fix ordered orthonormal bases $\mathbf u=[u_1,\dotsc,u_d]$ for $\Sing_{\max}$ and $\mathbf \Psi=[\psi_1,\dotsc,\psi_{d_0}]$ for $\W$. Let $\Vee\in \Gr_{d-d_0}(\Sing_{\max})$ and let $\mathbf \Phi=[\phi_1,\dotsc,\phi_{d-d_0}]$ be an ordered orthonormal basis of $\Vee$. There are unique matrices $V\in M^{d\times (d-d_0)}(\C)$, $W\in M^{d\times d_0}(\C)$ such that
\begin{equation*}
\mathbf \Phi=\mathbf u \cdot V,\quad \mathbf \Psi=\mathbf u\cdot W
\end{equation*}
Define
\begin{equation*}
\delta(\Vee,\W)=|\det[V|W]|.
\end{equation*}
The columns of $V$, likewise the columns of $W$, form an orthonormal set of vectors in $\C^d$. We claim that $\delta(\Vee,\W)$ is independent of the choices of orthonormal bases $\mathbf \Phi$ and $\mathbf \Psi$. Indeed, if $\mathbf \Phi'$ and $\mathbf \Psi'$ are other choices of ordered orthonormal bases of, respectively, $\Vee$ and $\W$, then $\mathbf \Phi'=\mathbf \Phi\cdot U_1$, $\mathbf \Psi'=\mathbf \Psi\cdot U_2$ with unitary matrices $U_1$ and $U_2$. Thus $\mathbf \Phi'= \mathbf u \cdot VU_1$ and $\mathbf \Psi'= \mathbf u \cdot WU_2$. But
\begin{equation*}
[VU_1,WU_2]=
 [V|W]
 \begin{bmatrix}
 U_1 & 0 \\ 0 & U_2
 \end{bmatrix}
\end{equation*}
so $|\det[VU_1|WU_2]|=|\det[V|W]\det U_1\det U_2|=|\det[V|W]|$ since unitary matrices have determinant of modulus $1$. Thus we get a globally defined function
\begin{equation*}
\delta:\Gr_{d-d_0}(\Sing_{\max})\times \Gr_{d_0}(\Sing_{\max})\to \R.
\end{equation*}
This function is clearly continuous, and $\mathfrak V_\W$ is the set of zeros of $\Vee\mapsto \delta(\Vee,\W)$. Suppose $\Vee\notin \mathfrak V_\W$ and let $\pi_{\Vee,\W}:\Sing_{\max}\to \Vee$ be the projection \eqref{GeneralProjectionEmax} on $\Vee$. The basis $\mathbf u$ can be written in terms of the basis $[\mathbf \Phi,\mathbf \Psi]$, as
\begin{equation*}
\mathbf u=[\mathbf \Phi,\mathbf \Psi]\cdot Q
\end{equation*}
where $Q$ is the inverse of $P=[V|W]$. Let $\tilde P$ be the matrix of minors of $P$, so that $Q=(\det P)^{-1}\tilde P$. The entries of $\tilde P$ are polynomials of degree $d-1$ in the entries of $P$. Since the columns of the latter matrix are vectors in the unit sphere in $\C^d$, the entries $\tilde p^j_k$ of $\tilde P$ are bounded by a constant independent of $P$. If $u=\sum a^\ell u_\ell\in \Sing_{\max}$, then the two terms in
\begin{equation*}
u=\bigg[\frac{1}{\det P}\sum_{k=1}^{d-d_0} \phi_k \sum_{\ell=1}^d\tilde p^k_\ell a^\ell\bigg]  + \bigg[\frac{1}{\det P}\sum_{k=1}^{d_0} \psi_k \sum_{\ell=1}^d\tilde p^{k+d-d_0}_\ell a^\ell\bigg]
\end{equation*}
correspond to the decomposition $\Sing_{\max}=\Vee\oplus \W$. Thus
\begin{equation*}
\pi_{\Vee,\W} u = \frac{1}{\det P}\sum_{k=1}^{d-d_0} \phi_k \sum_{\ell=1}^d\tilde p^k_\ell a^\ell.
\end{equation*}
This gives:

\begin{lemma}\label{BoundsViaGrassmann}
Let $\W\in \Gr_{d_0}(\Sing_{\max})$ and let $\Vee\in \Gr_{d-d_0}(\Sing_{\max})\minus \mathfrak V_{\W}$. Then
\begin{equation*}
\|\pi_{\Vee,\W}\|\leq \frac{C}{\delta(\Vee,\W)}.
\end{equation*}
The constant $C$ is independent of $\Vee$.
\end{lemma}

The question arises as to whether there is $\Dom\in \mathfrak G$ such that $\spec A_\Dom$ is discrete. The following proposition shows that if there is one such domain, then the set of such domains is open and connected, and its complement is a set with empty interior.

\begin{proposition}\label{BadSpectralVariety}
The set
\begin{equation*}
\mathfrak V=\set{\Dom\in \mathfrak G : \spec A_\Dom=\C}
\end{equation*}
is a variety. Thus, since $\mathfrak G$ is connected, $\mathfrak V \ne \mathfrak G$ if and only if $\mathfrak V$ has empty interior.
\end{proposition}
\begin{proof}
We identify $\mathfrak G$ with $\Gr_{d''}(\Sing_{\max})$ using the map \eqref{GrassmannMap}. Let $\gamma\to \Gr_{d''}(\Sing_{\max})$ be the canonical vector bundle. This is a holomorphic vector bundle. Let $\Dom_0\in \mathfrak G$ and let $u_1,\dotsc,u_{d''}$ be a holomorphic frame for $\gamma$ in a neighborhood $U$ of $\Dom_0$. Thus, if $u_1^0,\dotsc,u_d^0$ is a basis of $\Sing_{\max}$, then
\begin{equation*}
u_j(\Dom)=\sum_{\ell=1}^d g^\ell_j(\Dom)u^0_\ell,\quad j=1,\dotsc,d''
\end{equation*}
with holomorphic functions $g^\ell_j:U\to\C$. Any $u\in \Dom\in U$ can be written uniquely as
\begin{equation*}
u=v+\sum_{j=1}^{d''} \alpha^ju_j(\Dom)
\end{equation*}
with $v\in \Dom_{\min}$. For $\Dom\in U$ define $F(\Dom):\Dom_0\to\Dom$ by
\begin{equation*}
F(\Dom)(v + \sum_{j=1}^{d''} \alpha^j u_j(\Dom_0))
=v+\sum_{j=1}^{d''} \alpha^j u_j(\Dom),\quad v\in \Dom_{\min}.
\end{equation*}
This operator is bijective, and continuous in the graph norm of $A$. The operators
\begin{equation*}
A(\Dom,\lambda)=(A-\lambda)\circ F(\Dom):\Dom_0\to x^{-m/2}L^2_b(M;E)
\end{equation*}
depend holomorphically on $(\Dom,\lambda)\in U\times\C$, and the invertibility of $A_{\Dom}-\lambda$ is equivalent to the invertibility of $A(\Dom,\lambda)$.

If $\Dom_0\notin \mathfrak V$, then there is $\lambda_0\notin \spec A_{\Dom_0}$, and therefore, there is a neighborhood $U'\subset U$ of $\Dom_0$ and $\eps>0$ such that $A(\Dom,\lambda)$ is invertible for $(\Dom,\lambda)\in U'\times B(\lambda_0,\eps)$, where $B(\lambda_0,\eps)$ is the open disc in $\C$ with center $\lambda_0$ and radius $\eps$. Thus $U'$ is disjoint from $\mathfrak V$, which proves that $\mathfrak V$ is closed.

Suppose that $\lambda_0\in \spec A_{\Dom_0}$, let $K=\ker(A_{\Dom_0} - \lambda_0)$, and let $K^\perp$ be the orthogonal of $K$ in $\Dom_0$. Let $R=(A-\lambda_0)(\Dom_0)$, let $R^\perp$ be the orthogonal in $x^{-m/2}L^2_b(M;E)$ of $R$, and let $\pi_R$ and $\pi_{R^\perp}$ be the respective orthogonal projections. Define
\begin{align*}
A_{11}&=\pi_{R^\perp}A(\Dom,\lambda)|_K  &
A_{12}&=\pi_{R^\perp}A(\Dom,\lambda)|_{K^\perp}\\
A_{21}&=\pi_{R}A(\Dom,\lambda)|_K &
A_{22}&=\pi_{R}A(\Dom,\lambda)|_{K^\perp}
\end{align*}
so that
\begin{equation*}
A(\Dom,\lambda)=
\begin{bmatrix}
A_{11} & A_{12}\\
A_{21} & A_{22}
\end{bmatrix}
:
\begin{matrix}
K\\\oplus\\
K^\perp
\end{matrix}
\to
\begin{matrix}
R^\perp\\\oplus\\
R
\end{matrix}.
\end{equation*}
The operators $A_{ij}$ are continuous as operators into their target spaces as subspaces of $x^{-m/2}L^2_b(M;E)$ when their domains are given the graph norm of $A$, and depend holomorphically on $(\Dom,\lambda)$ for $\Dom\in U$ and $\lambda$ close to $\lambda_0$. Since $A_{22}(\Dom_0,\lambda_0)$ is invertible, we can, perhaps after shrinking $U$, find $\eps>0$ such that $A_{22}(\Dom,\lambda)$ is invertible if $\Dom\in U\times B(\lambda_0,\eps)$. If $(\Dom,\lambda)\in U\times B(\lambda_0,\eps)$, the elements of the kernel of $A_\Dom-\lambda$ are in one-to-one correspondence with the elements in the kernel of
\begin{equation*}
\mathcal A=A_{11}-A_{12}A_{22}^{-1}A_{21}:K\to R^\perp
\end{equation*}
via the map
\begin{equation*}
\ker \mathcal A \ni u \mapsto u-A_{22}^{-1}A_{21}u\in \ker A(\Dom,\lambda)\cong \ker A_\Dom-\lambda.
\end{equation*}
Since $A_\Dom$ has index $0$, $K$ and $R^\perp$ have the same dimension. Picking bases of $K$ and $R^\perp$ we can define a determinant $f(\Dom,\lambda)$ for $\mathcal A$. Since $\mathcal A$ depends holomorphically on $(\Dom,\lambda)$, so does $f(\Dom,\lambda)$. The set
\begin{equation*}
\set{(\Dom,\lambda)\in U\times B(\lambda_0,\eps): f(\Dom,\lambda)=0}
\end{equation*}
is the intersection of
\begin{equation*}
\mathfrak {spec} A=\set{(\Dom,\lambda)\in \mathfrak G\times\C:
\lambda \in \spec A_\Dom}
\end{equation*}
with $U\times B(\lambda_0,\eps)$ (thus $\mathfrak {spec}A$ is a variety). Write $f$ as
\begin{equation*}
f(\Dom,\lambda)=\sum_{\ell=0}^\infty f_\ell(\Dom)(\lambda-\lambda_0)^\ell;
\end{equation*}
the functions $f_\ell$ are holomorphic in $U$.  If $\Dom\in U\cap \mathfrak V$, then $f(\Dom,\lambda)=0$ for all $\lambda\in B(\lambda_0,\eps)$, so $f_\ell(\Dom)=0$. And if this condition holds for $\Dom$, then $f(\Dom,\lambda)=0$. So $\mathfrak V\cap U$ is the set of common zeros of the functions $f_\ell:U\to\C$, and $\mathfrak V$ is a variety.
\end{proof}

The following gives examples where $\mathfrak V$ is not empty.

\begin{example}
Let $A=e^{-\im \rho x}D_x$ on the interval $[-1,1]$, with $\rho\in \C$, $\rho\ne 0$. This is a cone operator:
\begin{equation*}
 A=(1-x^2)^{-1} e^{-\im\rho x}(1-x^2)D_x
\end{equation*}
and $(1-x^2)$ vanishes simply at $x=\pm 1$. We consider this operator initially as an unbounded operator
\begin{equation*}
C^\infty_0(-1,1)\subset (1-x^2)^{-1/2}L^2_b(-1,1) \to (1-x^2)^{-1/2}L^2_b(-1,1)
\end{equation*}
with the measure $\m=(1-x^2)^{-1}dx$. The space $(1-x^2)^{-1/2}L^2_b(-1,1)$ is just $L^2(-1,1)$ with the measure $dx$, so the domains of the minimal and maximal extensions are, respectively, the standard Sobolev spaces $H^1_0[-1,1]$ and $H^1(-1,1)$. Since $H^1(-1,1)$ consists of continuous functions on $[-1,1]$, the elements in $\Dom_{\max}$ can be evaluated at $x=-1$ and at $x=1$. The Mellin transforms at either boundary of elements of $(1-x^2)^{-1/2}L^2_b(-1,1)$ are holomorphic in $\Im \sigma>1/2$, of course. To compute the conormal symbol of $P=e^{-\im\rho x}(1-x^2)D_x$ at $x=-1$ let $x_{\Left}=1+x$. Then
\begin{equation*}
P=(2-x_{\Left}) e^{-\im\rho (x_L-1)} x_{\Left} D_{x_{\Left}} = 2e^{-\im\rho (x_L-1)}x_{\Left}D_{x_{\Left}}-x_{\Left}e^{-\im\rho (x_L-1)}x_{\Left}D_{x_{\Left}}
\end{equation*}
so the conormal symbol of $P$ at $x=-1$ with respect to $x_{\Left}$ is $2\sigma e^{\im\rho}$, giving a simple pole at $\sigma=0$ for the inverse of the conormal symbol. If $u\in \Dom_{\max}$, its value at $x=-1$ is essentially the residue at $\sigma=0$ of the Mellin transform of $u$. Using $x_{\Right}=1-x$ as defining function for $\set {x=1}$ we get
\begin{equation*}
P=-(2-x_{\Right}) e^{-\im\rho (1-x_R)} x_{\Right} D_{x_{\Right}} = -2e^{-\im\rho (1-x_R)}x_{\Right}D_{x_{\Right}}+x_{\Right}e^{-\im\rho (1-x_R)} x_{\Right}D_{x_{\Right}}
\end{equation*}
and the conormal symbol at that boundary is $-2\sigma e^{-\im\rho}$. Since the only point in $\spec_b(P)$ is $0$, we deduce that $\Dom_{\min}=(1-x^2)^{1/2}H^1_b$ and that
\begin{equation*}
\Dom_{\min}=\set{ u\in \Dom_{\max}: u(-1)=u(1)=0}.
\end{equation*}
The operator $A$ with the minimal domain is injective. The formal adjoint of $A$ is $A^\star = e^{\im\overline\rho x}(D_x+\overline \rho)$, and the Hilbert space adjoint of $A_{\Dom_{\min}}$ is $A^\star$ with its maximal domain, $\Dom^\star_{\max}$. The latter contains the function $e^{-\im \overline \rho x}$, which spans the kernel of $A^\star$, so the index of $A_{\Dom_{\min}}$ is $-1$. This also gives that the index of $A_{\Dom_{\max}}$ is $+1$.

The domains on which $A$ has index $0$ are of the form
\begin{equation*}
\Dom_{\alpha_-,\alpha_+}=\set{ u\in \Dom_{\max} :  \alpha_-u(-1) + \alpha_+u(1)=0}
\end{equation*}
with $(\alpha_-, \alpha_+)\in \C^2\minus 0$. If $z\ne 0$ then $(z\alpha_-,z\alpha_+)$ determines the same domain as $(\alpha_-,\alpha_+)$, so $\mathfrak G$, the manifold of domains where $A$ has index $0$, is $\C\mathbb P^1=S^2$.

Fix some $(\alpha_-,\alpha_+)\in \C^2\minus 0$. The kernel of $(A-\lambda)$ on $\Dom_{\max}$ is spanned by
\begin{equation*}
h_\lambda(x)=e^{\lambda e^{\im\rho x}/\rho}.
\end{equation*}
The condition that $h_\lambda\in \Dom_{\alpha_-,\alpha_+}$ is
\begin{equation*}
\alpha_-e^{\lambda e^{-\im\rho}/\rho}+
\alpha_+e^{\lambda e^{ \im \rho }/\rho}=0,
\end{equation*}
equivalently
\begin{equation*}
\alpha_- +\alpha_+e^{2\im\lambda\rho^{-1}\sin\rho}
= \alpha_- + \alpha_+e^{\lambda [e^{\im \rho}-e^{-\im \rho }]/\rho}
=0.
\end{equation*}
Thus, if $\rho\in \pi\mathbb Z$ ($\rho\ne 0$), then $\alpha_-+\alpha_+=0$ implies $\spec A_{\Dom_{\alpha_-,\alpha_+}} = \C$ while $\alpha_-+\alpha_+\ne 0$ implies $\spec A_{\Dom_{\alpha_-,\alpha_+}} = \varnothing$. And if $\rho\notin \pi\mathbb Z$, then for any $(\alpha_+,\alpha_-)\in \C^2\minus 0$, the spectrum of $A_{\Dom_{\alpha_-,\alpha_+}}$ is discrete, and empty if either $\alpha_-$ or $\alpha_+=0$.
\end{example}

%%%%%%%%%%%%%%%%%%%%%%%%%%%%%%%%%%%%%%%%%%%%%%%%%%%%%%%%%%%
%%%%%%%%%%%%%%%%%%%%%%%%%%%%%%%%%%%%%%%%%%%%%%%%%%%%%%%%%%%
\section{Selfadjointness }
\label{sec-SelfAdjointness}

We now discuss the important case where $A$ is symmetric on $\Dom_{\min}$ from the perspective of Section \ref{sec-DomainsAndSpectra}. The selfadjoint extensions of such operators were studied by Lesch \cite{Le97}. Suppose $A$ is such a $c$-elliptic symmetric operator. Since
\begin{equation*}
\|(A-\lambda)u\|\geq |\Im\lambda|\|u\|\quad \text{ if }u\in \Dom_{\min},
\end{equation*}
$A_{\Dom_{\min}}-\lambda$ is injective when $\Im\lambda\ne 0$. Since $A$ is Fredholm and the Hilbert space adjoint of $A_{\Dom_{\min}}$ is $A$ with domain $\Dom_{\max}$, $A_{\Dom_{\max}}-\lambda$ is surjective if $\Im\lambda\ne 0$. Since the operators $A_{\Dom_{\min}}-\lambda$ are Fredholm and depend continuously on $\lambda$, the indices at $\lambda=\im$ and $\lambda=-\im$ are equal. So the deficiency indices are the same, and $A$ admits selfadjoint extensions. If $A_\Dom$ is one such extension, then $\spec A_\Dom\subset \R$, therefore $\bgspec A$ is a discrete subset of $\R$.

The Dirichlet form of a general cone operator $A$ is the sesquilinear form
\begin{equation}\label{DirichletPairing}
[u,v]_A=(Au,v)-(u,A^\star v),\quad u\in \Dom_{\max}(A),\ v\in \Dom_{\max}(A^\star).
\end{equation}
It has the property that
\begin{equation*}
[u,v]_A=[\pi_{\max} u,\pi_{\max} v]_A,
\end{equation*}
because $[\pi_{\max}u,\pi_{\min} v]_A=[\pi_{\min} u,\pi_{\min}v]_A=0$ for any $u$ and $v$. Moreover, the induced sesquilinear pairing
\begin{equation}\label{NonSingPairing}
\Sing_{\max}(A)\times\Sing_{\max}(A^\star)\to\C \text{ is nonsingular}
\end{equation}
(cf. Theorems 7.11 and 7.17 in \cite{GiMe01}). If $\Dom\in \mathfrak D(A)$, let $\mathcal J\Dom\in \mathfrak D(A^\star)$ be the annihilator of $\Dom$ with respect to the pairing \eqref{DirichletPairing}. Thus if $\Dom\in \mathfrak D$, then the Hilbert space adjoint of $A_\Dom$ is $A^\star_{\mathcal J\Dom}$. We will prove in a moment that the mapping $\mathcal J:\mathfrak D(A)\to\mathfrak D(A^\star)$ is real-analytic.
Let $\mathcal J^\star:\mathfrak D(A^\star)\to\mathfrak D(A)$ be the analogously defined map. Clearly $\mathcal J^\star\mathcal J$ is the identity. If $A$ is symmetric on $\Dom_{\min}$ and $\Dom\in \mathfrak G$, then $\mathcal J:\mathfrak G\to \mathfrak G$, $\mathcal J^\star=\mathcal J$,  and $A_\Dom$ is selfadjoint if and only if $\Dom$ is a fixed point of $\mathcal J$. Such domains will be called selfadjoint.

\begin{lemma}\label{DualityLemma}
Let $A$ be an arbitrary $c$-elliptic cone operator.
\begin{enumerate}[$(i)$]
\item If $u\in \Sing_{\max}(A)$, then $Au\in \Sing_{\max}(A^\star)$. The map
\begin{equation*}
\Sing_{\max}(A) \ni u \mapsto Au\in \Sing_{\max}(A^\star)
\end{equation*}
is an isometry.
\item If $u\in \Sing_{\max}(A)$ and $v\in \Dom_{\max}(A^\star)$, then
\begin{equation*}
(Au,v)_{A^\star}=[u,v]_A.
\end{equation*}
\item Let $\Dom\in \mathfrak D(A)$. If $u\in \Dom\cap\Sing_{\max}(A)$, then $Au$ is orthogonal to $\mathcal J \Dom$ with respect to the inner product defined by $A^\star$.
\end{enumerate}
\end{lemma}
\begin{proof}
To prove the first assertion in (i), suppose $u\in \Sing_{\max}(A)$. Then $Au\in\Dom_{\max}(A^\star)$. In addition, $A^\star Au=-u$, so $AA^\star (Au)=-Au$, that is, $Au\in \Sing_{\max}(A^\star)$. If $u$, $v\in \Sing_{\max}(A)$, then
\begin{equation*}
(Au,Av)_{A^\star}=(A^\star Au,A^\star Av)+(Au,Av)=(u,v)+(Au,Av)=(u,v)_A.
\end{equation*}
For part (ii), suppose $u$ and $v$ are as indicated. Then
\begin{equation*}
(Au,v)_{A^\star}=(A^\star A u, A^\star v)+(Au,v)=(-u,A^\star v)+(Au,v)=[u,v]_A.
\end{equation*}
For part (iii) we observe that if $u\in \Dom$ and $v\in \mathcal J\Dom$, then $[u,v]_A=0$, and use part (ii).
\end{proof}

\begin{proposition}
The mapping $\mathcal J:\mathfrak D(A)\to\mathfrak D(A^\star)$ is real-analytic.
\end{proposition}
\begin{proof}
Let $\Dom_0\in \mathfrak D(A)$, and let $\phi_1,\dotsc,\phi_d$ be an $A$-orthonormal basis of $\Sing_{\max}(A)$ whose first $d_0$ elements form a basis of $\Dom_0\cap \Sing_{\max}(A)$. Let $\psi_j=A\phi_j$, $j=1,\dotsc,d$. The $\psi_j$ form an orthonormal basis of $\Sing_{\max}(A^\star)$, by part (i) of Lemma \ref{DualityLemma}. Therefore, by part (ii) of the same lemma, \begin{equation}\label{DualityRelations}
[\phi_j,\psi_k]_A=\delta_{jk}.
\end{equation}
We deduce that
\begin{equation*}
\mathcal J\Dom_0=\LinSpan\set{\psi_{d_0+1},\dotsc,\psi_d}\oplus \Dom_{\min}(A^\star).
\end{equation*}
Write $\mathbf \Phi_1=[\phi_1,\dotsc,\phi_{d_0}]$, $\mathbf \Phi_2=[\phi_{d_0+1},\dotsc,\phi_d]$, and analogously $\mathbf \Psi_1$, $\mathbf \Psi_2$. The set
\begin{equation*}
\set{\LinSpan(\mathbf \Phi_1 + \mathbf \Phi_2\cdot Z): Z\in M^{(d-d_0)\times d_0}(\C)}
\end{equation*}
is a neighborhood of $\Dom_0$ in (a component of) $\mathfrak D(A)$, and the map
\begin{equation*}
Z\mapsto \LinSpan(\mathbf \Phi_1 + \mathbf \Phi_2\cdot Z)
\end{equation*}
is the inverse of a holomorphic chart. Likewise, parametrize the $(d-d_0)$-dimensional subspaces of $\Sing_{\max}(A^\star)$ in a neighborhood of $\mathcal J\Dom_0$ by
\begin{equation*}
W\mapsto \LinSpan(\mathbf \Psi_2 + \mathbf \Psi_1\cdot W)
\end{equation*}
with $W\in M^{d_0\times(d-d_0)}(\C)$. The condition that the vector space spanned by
\begin{equation*}
\psi_{k+d_0}+\sum_{j=1}^{d_0} W^j_k\psi_j,\quad k=1,\dotsc,d-d_0
\end{equation*}
is $[\cdot,\cdot]_A$-orthogonal to
\begin{equation*}
\phi_j+\sum_{k=1}^{d-d_0} Z^{k}_j\phi_{k+d_0},\quad j=1,\dotsc,d_0
\end{equation*}
is $Z^k_j+\overline W^j_k=0$ because of \eqref{DualityRelations}. Thus, in coordinates, $\mathcal J$ maps the space determined by $Z$ to the space determined by $W=-Z^*$. We conclude that $\mathcal J$ is real-analytic.
\end{proof}

Since $\mathcal J$ is real-analytic, its set of fixed points is a real-analytic variety. In fact:

\begin{proposition}
Let $A$ be symmetric on $\Dom_{\min}$. The set $\mathfrak {SA}$ of domains $\Dom\in \mathfrak G$ such that $A_\Dom$ is selfadjoint is a real-analytic (smooth) submanifold of $\mathfrak G$ of codimension $(d')^2$, $d'=\Ind A_{\Dom_{\max}}$.
\end{proposition}

\begin{proof}
If $\Ind A_{\Dom_{\max}}=0$ then also $\Ind A_{\Dom_{\min}}=0$, $\Dom_{\min}=\Dom_{\max}$, and $A_{\Dom_{\min}}$ is the only selfadjoint extension of $A$. Assume then that $\Ind A_{\Dom_{\max}}>0$ and pick a selfadjoint domain $\Dom_0$. Let $\phi_1,\dotsc,\phi_{d'}$ be an orthonormal basis of $\Dom_0\cap \Sing_{\max}$. Then $[\phi_j,\phi_k]_A=0$. By part (i) of Lemma \ref{DualityLemma}, $(A\phi_j,A\phi_k)_A=\delta_{jk}$. Thus by part (ii), $[\phi_j,A\phi_k]_A=(A\phi_j,A\phi_k)_A=\delta_{jk}$. Also by part (ii), $[A\phi_j,A\phi_k]_A=-(\phi_j,A\phi_j)_A$, which vanishes by part (iii). As above, write $\mathbf \Phi_1=[\phi_1,\dotsc,\phi_{d'}]$ and let $\mathbf \Phi_2=[A\phi_1,\dotsc,A\phi_{d'}]$. So a neighborhood $U\subset \mathfrak G$ of $\Dom_0$ is parametrized by the vector spaces associated with the bases $\mathbf \Phi_1+\mathbf \Phi_2\cdot Z$, $Z\in \mathfrak{gl}(\C,d')$. Writing the components of $\mathbf \Phi_1 + \mathbf \Phi_2\cdot Z$ as
\begin{equation*}
\phi_j+\sum_{k=1}^{d'} Z^k_j A\phi_k,\quad j=1,\dotsc,d'
\end{equation*}
we see that the selfadjoint domains in $U$ are those that satisfy
\begin{equation*}
Z^j_k-\overline Z^k_j=0
\end{equation*}
(\ie, $Z$ is a selfadjoint matrix). These equations represent $(d')^2$ real-analytic conditions.
\end{proof}

\begin{proposition}\label{RealPtsInSAspectrum}
Let $A$ be symmetric on $\Dom_{\min}$ and assume that $-\Ind A_{\Dom_{\min}}>0$. For any $\lambda\in \R$ there is $\Dom\in \mathfrak {SA}$ such that $\lambda\in \spec A_\Dom$.
\end{proposition}

\begin{proof}
If $\lambda$ belongs to $\bgspec A$, then $\lambda$ already belongs to the spectrum of any extension, selfadjoint or not, of $A$. Suppose $\lambda \in \bgres A$. If $u$, $v \in \KK_\lambda$, cf. Proposition \ref{BundleOfKernels}, then
\begin{equation*}
[u,v]_A = 2i \Im(\lambda) (u,v),
\end{equation*}
so if $\lambda$ is real, then the Dirichlet form of $A$ vanishes on $\KK_\lambda$. Since
\begin{equation*}
\dim \KK_\lambda = \Ind A_{\Dom_{\max}}=-\Ind A_{\Dom_{\min}},
\end{equation*}
and since $\pi_{\max}$ is injective on $\KK_\lambda$,
\begin{equation}\label{DLambda}
\Dom_\lambda=\KK_\lambda + \Dom_{\min}
\end{equation}
is an element of $\mathfrak G$ on which the Dirichlet form vanishes. Thus $A_\Dom$ is selfadjoint, and $\lambda\in \spec A_\Dom$.
\end{proof}

Note that there is no assumption on semiboundedness of $A$.

\begin{proposition}
Let $A$ be symmetric on $\Dom_{\min}$. Then
\begin{equation}\label{SAspecIsAll}
\bgspec A = \bigcap_{\Dom\in \mathfrak{SA}} \spec A_\Dom.
\end{equation}
\end{proposition}
\begin{proof}
If $\Ind A_{\Dom_{\min}}=0$, then $\mathfrak G=\set{\Dom_{\min}}$ and $\bgspec A=\spec A_{\Dom_{\min}}$. Since $A_{\Dom_{\min}}$ is already selfadjoint, \eqref{SAspecIsAll} is an identity.

Suppose then that $-\Ind A_{\Dom_{\min}}>0$. Denote the set on the right in \eqref{SAspecIsAll} by $S$. From the definition of $\bgspec A$ we get $\bgspec A\subset S$.

To prove the opposite inclusion suppose that $\lambda_0\in \bgres A$. If $\Im\lambda_0 \ne 0$, then $\lambda_0\notin S$, since $S\subset \R$. If $\lambda_0\in \R\cap \bgres A$, consider $A_{\Dom_{\lambda_0}}$, where $\Dom_{\lambda_0}$ is as in \eqref{DLambda}. From the proof of Proposition \ref{RealPtsInSAspectrum} we know that $A_{\Dom_{\lambda_0}}$ is selfadjoint. Since $A_{\Dom_{\lambda_0}}$ is Fredholm and $\spec A_{\Dom_{\lambda_0}}\ne \C$, this spectrum is discrete. Since $\KK_{\lambda_0}\subset \Dom_{\lambda_0}$, $\lambda_0\in \spec A_{\Dom_{\lambda_0}}$. We can therefore find a neighborhood $U\subset \bgres A$ of $\lambda_0$ with the property that $U\cap \spec A_{\Dom_{\lambda_0}}=\set{\lambda_0}$. We claim that if $\lambda\in U\minus\set{\lambda_0}$, then $\lambda_0\notin \spec A_{\Dom_\lambda}$. To see this, let $\lambda\in U$ and assume that $\lambda_0\in \spec A_{\Dom_{\lambda}}$. Then $\KK_{\lambda_0}\cap \Dom_\lambda\ne 0$. Thus there are $\phi\in \KK_{\lambda_0}$ with $\phi\ne 0$, and $\psi\in \KK_{\lambda}$, $v\in \Dom_{\min}$ such that $\phi=\psi+v$. The element $\phi-v\in \Dom_{\lambda_0}$ is equal to $\psi$, so $\Dom_{\lambda_0}\cap \KK_{\lambda}\ne 0$. Necessarily $\psi\ne 0$, since $\pi_{\max}\phi=\pi_{\max}\psi$ and $\phi\ne 0$. Thus $\lambda\in U\cap \spec A_{\Dom_{\lambda_0}}$, which implies $\lambda=\lambda_0$. It follows that if $\lambda\in U\cap\R\minus \lambda_0$, then $\Dom_\lambda\in \mathfrak {SA}$ and $\lambda_0\notin \spec A_{\Dom_{\lambda}}$, hence $\lambda_0\notin S$. Therefore $S\subset \bgspec A$.
\end{proof}

%%%%%%%%%%%%%%%%%%%%%%%%%%%%%%%%%%%%%%%%%%%%%%%%%%%%%%%%%%%
%%%%%%%%%%%%%%%%%%%%%%%%%%%%%%%%%%%%%%%%%%%%%%%%%%%%%%%%%%%
\section{The model operator}
\label{sec-ModelOperator}

In this section we focus on the spectra of closed extensions of the operator $A_\wedge$, cf \eqref{InvariantAWedge}. We continue to assume that the operator $A\in x^{-m}\Diff^m_b(M;E)$ is $c$-elliptic. We will usually write $\Dom_{\wedge,\min}$ for $\Dom_{\min}(A_\wedge)$ and $\Dom_{\wedge,\max}$ for $\Dom_{\max}(A_\wedge)$. Recall that the inner product on $\Dom_{\wedge,\max}$ is given by \eqref{AWedgeInner}. The nature of $\Dom_{\wedge,\min}$ was described in Proposition \ref{OnWedgeDomains}. We also noted there that $\Dom_{\wedge,\max}/\Dom_{\wedge,\min}$ is finite dimensional. Because of the finite dimensionality of this quotient, many of
the results concerning the closed extensions of $A$ find their
analogue in the situation at hand, despite the fact that neither of the operators
\begin{equation*}\label{MinimalAWedge}
A_\wedge:\Dom_{\wedge,\min}\subset x^{-m/2}L^2_b(Y^\wedge;E)\to x^{-m/2}L^2_b(Y^\wedge;E)
\end{equation*}
nor
\begin{equation*}\label{MaximalAWedge}
A_\wedge:\Dom_{\wedge,\max}\subset x^{-m/2}L^2_b(Y^\wedge;E)\to x^{-m/2}L^2_b(Y^\wedge;E)
\end{equation*}
needs to be Fredholm. On the other hand,
the homogeneity property
\begin{equation}\label{HomogeneityOfAwedgeBis}
 A_\wedge - \lambda =
 \varrho^m \kappa_\varrho (A_\wedge -\lambda/\varrho^m)\kappa_\varrho^{-1}
 \;\text{ for every } \varrho>0
\end{equation}
of $A_\wedge-\lambda$, $\lambda\in \C$, cf. \eqref{HomogeneityOfA}, not available in such simple form in the case of $A$, permits an essentially complete understanding of the spectra and resolvents for the closed extensions of $A_\wedge$.

We begin our analysis with:
\begin{definition}\label{bgSpecAndResWedge}
The background spectrum of $A_\wedge$ is the set
\begin{equation*}
\bgspec A_\wedge=\set{\lambda\in \C:\lambda\in \spec A_{\wedge,\Dom}\ \forall \Dom\in \mathfrak D_\wedge}.
\end{equation*}
The complement of this set, $\bgres A_\wedge$, is the background resolvent set.
\end{definition}

The analogue
\begin{equation*}
\bgres A_\wedge =\{\lambda\in \C: A_{\wedge,\Dom_{\min}}-\lambda \text{ is injective  and } A_{\wedge,\Dom_{\max}}-\lambda\text{ is surjective}\}
\end{equation*}
of Lemma \ref{Spec2} holds for $A_\wedge$ in place of $A$, with the same proof.

\begin{lemma}\label{GoodSectorWedge}
If $\lambda\in \bgres A_\wedge$ and $\Dom\in \mathfrak D_\wedge$, then $A_{\wedge,\Dom} -\lambda$ is Fredholm. The set $\bgres A_\wedge$ is a union of open sectors.
\end{lemma}
\begin{proof}
Let $\lambda\in \bgres A_\wedge$. Since $\Dom_{\wedge,\max}/\Dom_{\wedge,\min}$ is finite dimensional and $A_\wedge-\lambda$ is injective on $\Dom_{\wedge,\min}$, $A_{\wedge,\Dom_{\max}}-\lambda$ has finite dimensional kernel. Thus $A_{\wedge,\Dom_{\max}}-\lambda$ is Fredholm, and so is its restriction to any subspace of $\Dom_{\wedge,\max}$ of finite codimension.

Next, suppose that $\lambda_0\in \bgres A_\wedge$ and let $\lambda=\varrho^m\lambda_0$. Since $\kappa_\varrho$ is invertible and
\begin{equation*}
 A_{\wedge} - \lambda =
 \varrho^m{\kappa}_{\varrho}(A_{\wedge} -
 \lambda_0){\kappa}_{\varrho}^{-1},
\end{equation*}
$ A_{\wedge} - \lambda$ is injective on $\Dom_{\wedge,\min}$ and surjective on $\Dom_{\wedge,\max}$. Thus the ray $\set{r\lambda_0: r>0}$ is contained in $\bgres A_\wedge$. Since $\lambda_0\in \bgres A_\wedge$, $A_\wedge-\lambda_0$ admits a continuous left inverse $B(\lambda_0):x^{-m/2}L^2_b(Y^\wedge;E)\to \Dom_{\wedge,\min}$. Since the inclusion $\Dom_{\wedge,\min}\embed x^{-m/2}L^2_b(Y^\wedge;E)$ is continuous, the formula
\begin{equation*}
B(\lambda_0)(A_\wedge-\lambda)=I+(\lambda_0-\lambda)B(\lambda_0)
\end{equation*}
gives that $(A_{\wedge,\min}-\lambda)$ admits a left inverse if $\lambda$ is close to $\lambda_0$. Likewise $(A_{\wedge,\max}-\lambda)$ admits a right inverse if $\lambda$ is close to $\lambda_0$. So $\bgres A_\wedge$ is open. Therefore its connected components are open sectors.
\end{proof}

Label the connected components of $\bgres A_\wedge$ by $\open \Lambda_\alpha$, $\alpha\in \mathfrak I\subset \N$. Since the inclusion map $\Dom\embed x^{-m/2}L^2_b(Y^\wedge;E)$ is continuous for any $\Dom\in \mathfrak D_\wedge$,
\begin{equation*}
\open \Lambda_\alpha\ni \lambda\mapsto \Ind (A_{\wedge,\Dom}-\lambda)
\end{equation*}
is constant, and
\begin{equation*}
\Ind(A_{\wedge,\Dom}-\lambda)=\Ind(A_{\wedge,\Dom_{\min}}-\lambda)+\dim\Dom/\Dom_{\wedge,\min},\quad\lambda\in \bgres A_\wedge.
\end{equation*}
Let
\begin{equation*}
d'_\alpha=\Ind(A_{\wedge,\Dom_{\max}}-\lambda), \quad d''_\alpha=-\Ind(A_{\wedge,\Dom_{\min}}-\lambda),\quad \lambda\in \open \Lambda_\alpha
\end{equation*}
and let
\begin{equation*}
\mathfrak G_{\wedge,\alpha}=\set{\Dom\in \mathfrak D_\wedge:
\dim\Dom/\Dom_{\wedge,\min}=d''_\alpha},\quad \alpha\in \mathfrak I.
\end{equation*}
The elements of $\mathfrak G_{\wedge,\alpha}$ are thus the domains $\Dom$ for which $A_{\wedge,\Dom}-\lambda$ has index $0$ when $\lambda\in \open \Lambda_\alpha$. Write $\Sing_{\wedge,\max}$ for the orthogonal of $\Dom_{\wedge,\min}$ in $\Dom_{\wedge,\max}$. Using that
\begin{equation*}
\mathfrak D_{\wedge}\ni \Dom\mapsto \Dom\cap \Sing_{\wedge,\max}
\end{equation*}
is a bijection onto the set of finite dimensional subspaces of $\Sing_{\wedge,\max}$ we give each of the $\mathfrak G_{\wedge,\alpha}$ the structure of a complex manifold.

The proofs of the following lemma and proposition parallel the arguments in the proofs of Propositions \ref{Spec4} and \ref{BadSpectralVariety}, respectively.

\begin{lemma}
For every $\alpha\in\mathfrak I$ such that $\dim \mathfrak G_{\wedge,\alpha,0}>0$ and every $\lambda\in \open \Lambda_\alpha$ there is $\Dom_0 \in \mathfrak G_{\wedge,\alpha}$ such that $\lambda\in \spec A_{\wedge,\Dom_0}$.
\end{lemma}

\begin{proposition}\label{BadSpectralVarietyWedge}
For every $\alpha\in \mathfrak I$ the set
\begin{equation*}
\mathfrak V_\alpha=\set{\Dom\in \mathfrak G_{\wedge,\alpha}: \open \Lambda_\alpha\subset \spec A_{\wedge,\Dom}}
\end{equation*}
is a variety.
\end{proposition}

If $\Dom\in \mathfrak D_\wedge$, then $\kappa_\varrho\Dom$ is again an element of $\mathfrak D_\wedge$. Indeed,
\begin{equation*}
\Dom_{\wedge,\min}\subset \kappa_\varrho\Dom
\subset \Dom_{\wedge,\max}
\end{equation*}
since both $\Dom_{\wedge,\min}$ and $\Dom_{\wedge,\max}$ are $\kappa$-invariant. Define
\begin{equation*}
\pmb \kappa_\varrho:\mathfrak D_\wedge\to \mathfrak D_\wedge,\quad \pmb \kappa_\varrho(\Dom)=\kappa_\varrho \Dom.
\end{equation*}
Since $\kappa_\varrho\Dom_{\wedge,\min}=\Dom_{\wedge,\min}$,
\begin{equation*}
\pi_{\wedge,\max}\kappa_\varrho= \pi_{\wedge,\max}\kappa_\varrho\pi_{\wedge,\max},
\end{equation*}
and therefore the map
\begin{equation}\label{RestrictedKappa}
\R\ni \xi \mapsto \pi_{\wedge,\max}\kappa_{e^\xi}\big|_{\Sing_{\wedge,\max}}:\Sing_{\wedge,\max}\to\Sing_{\wedge,\max}
\end{equation}
is a (continuous) one-parameter group of isomorphisms of $\Sing_{\wedge,\max}$, necessarily given by exponentiation of its infinitesimal generator. So \eqref{RestrictedKappa} extends to a holomorphic action of $\C$ on $\Sing_{\wedge,\max}$. We will use the notation $\pmb\kappa_\varrho(\Vee)$ for $\pi_{\wedge,\max}\kappa_\varrho(\Vee)$ when $\Vee\subset \Sing_{\wedge,\max}$ is a subspace.

\begin{proposition}\label{KappaVectorField}
Let $d_0<d=\dim\Sing_{\wedge,\max}$ be a nonnegative integer. The map
\begin{equation*}
\R\times \Gr_{d_0}(\Sing_{\wedge,\max}) \ni (\xi,\Vee)\mapsto \pmb \kappa_{e^\xi}\Vee \in \Gr_{d_0}(\Sing_{\wedge,\max})
\end{equation*}
extends to a holomorphic map
\begin{equation*}
\pmb \kappa_{\exp}:\C\times \Gr_{d_0}(\Sing_{\wedge,\max})\to
\Gr_{d_0}(\Sing_{\wedge,\max})
\end{equation*}
with the property that
\begin{equation*}
\pmb \kappa_{\exp}(\zeta+\zeta',\Vee)=\pmb \kappa_{\exp}(\zeta,\pmb \kappa_{\exp}(\zeta',\Vee))
\end{equation*}
for all $\zeta$, $\zeta'\in \C$ and $\Vee\in \Gr_{d_0}(\Sing_{\wedge,\max})$. In particular, for each $\Vee\in \Gr_{d_0}(\Sing_{\wedge,\max})$, the curve
\begin{equation*}
\R \ni \xi\mapsto \pmb \kappa_{e^{\xi}}\Vee \in \Gr_{d_0}(\Sing_{\wedge,\max})
\end{equation*}
is real-analytic, and the infinitesimal generator of the group action $\pmb \kappa_{\exp}$ is the real part of a holomorphic vector field.
\end{proposition}

\begin{proof}
The proof is an elementary argument on Grassmannian varieties.
Let $\Vee_0\in \Gr_{d_0}(\Sing_{\wedge,\max})$ and pick a basis $\mathbf \Phi = [\phi_1,\dotsc,\phi_d]$ of $\Sing_{\wedge,\max}$ whose first $d_0$ elements form a basis of $\Vee_0$. Then $\pi_{\wedge,\max}\kappa_{e^\zeta}\big|_{\Sing_{\wedge,\max}}$ sends the basis $\mathbf \Phi$ to the basis $\mathbf \Phi\cdot \pmb \kappa(\zeta)$ whose $j$-th component is
\begin{equation*}
\sum_{k=1}^d \phi_k \pmb \kappa^k_j(\zeta);
\end{equation*}
the matrix $\pmb \kappa(\zeta)=[\kappa^k_j(\zeta)]$ depends holomorphically on $\zeta$. Let $\mathbf \Phi_1=[\phi_1,\dotsc,\phi_{d_0}]$, $\mathbf \Phi_2=[\phi_{d_0+1},\dotsc,\phi_d]$. If $Z\in M^{(d-d_0)\times d_0}(\C)$ is a $(d-d_0)\times d_0$ matrix with complex entries, then $\mathbf \Phi_2\cdot Z$ is defined, the entries of $\mathbf \Phi_1+\mathbf \Phi_2\cdot Z$ are independent, and
\begin{equation*}
\Vee(Z)=\LinSpan(\mathbf \Phi_1+\mathbf \Phi_2\cdot Z)\subset \Sing_{\wedge,\max}
\end{equation*}
defines an element of the Grassmannian $\Gr_{d_0}(\Sing_{\wedge,\max})$. For a fixed basis $\mathbf \Phi$ the collection of elements $\Vee(Z)$ is a neighborhood $U$ of $\Vee_0$ and $Z\mapsto \Vee(Z)$ is the inverse of a holomorphic chart of $\Gr_{d_0}(\Sing_{\wedge,\max})$. Write the $d \times d$ matrix $\pmb\kappa(\zeta)$ in block form,
\begin{equation*}
\pmb \kappa(\zeta)=
\begin{bmatrix}
\pmb \kappa^1_1(\zeta) & \pmb \kappa^1_2(\zeta)\\
\pmb \kappa^2_1(\zeta) & \pmb \kappa^2_2(\zeta)
\end{bmatrix},
\end{equation*}
with $\kappa^1_1(\zeta)\in M^{d_0\times d_0}(\C)$. With this notation,  $\pi_{\wedge,\max}\kappa_{e^\zeta}\big|_{\Sing_{\max}}$ maps the components of $\mathbf \Phi_1 + \mathbf \Phi_2\cdot Z$ to the components of
\begin{equation*}
\mathbf \Phi_1\cdot \big(\pmb \kappa^1_1(\zeta)+\pmb \kappa^1_2(\zeta)Z\big)+\mathbf \Phi_2\cdot \big(\pmb \kappa^2_1(\zeta) +\pmb \kappa^2_2(\zeta)Z\big).
\end{equation*}
If $Z$ belongs to a bounded set in $M^{(d-d_0)\times d_0}(\C)$ then for $\zeta$ small enough the matrix $\pmb \kappa^1_1(\zeta)+\pmb \kappa^1_2(\zeta)Z$ is invertible, since $\pmb \kappa(0)=\id$, and we get from
\begin{equation*}
\{\mathbf \Phi_1 +\mathbf \Phi_2\cdot \big(\pmb \kappa^2_1(\zeta) +\pmb \kappa^2_2(\zeta)Z\big) \big(\pmb \kappa^1_1(\zeta)+\pmb \kappa^1_2(\zeta)Z\big)^{-1}\}\cdot \big(\pmb \kappa^1_1(\zeta)+\pmb \kappa^1_2(\zeta)Z\big)
\end{equation*}
that $\kappa_{e^\zeta}$ maps the point in $U$ of coordinates $Z$ to the point in $U$ of coordinates $\big(\pmb \kappa^2_1(\zeta) +\pmb \kappa^2_2(\zeta)Z\big) \big(\pmb \kappa^1_1(\zeta)+\pmb \kappa^1_2(\zeta)Z\big)^{-1}$. The latter is a holomorphic function of $\zeta$ and $Z$.
\end{proof}

If $\Dom\in \mathfrak D_\wedge$, then $\Dom=\pi_{\wedge,\max}\Dom\oplus \Dom_{\wedge,\min}$. Therefore
\begin{equation*}
\kappa_{e^\xi}\Dom=(\pi_{\wedge,\max}\kappa_{e^\xi}\pi_{\wedge,\max}\Dom) \oplus \Dom_{\wedge,\min}
\end{equation*}
for real $\xi$. For $\zeta\in \C$, define
\begin{equation*}
\kappa_{e^\zeta}\Dom=(\pi_{\wedge,\max}\kappa_{e^\zeta}\pi_{\wedge,\max}\Dom) \oplus \Dom_{\wedge,\min}.
\end{equation*}
The $\pmb \kappa_{e^\xi}: \Gr_{d_0}(\Sing_{\wedge,\max})\to \Gr_{d_0}(\Sing_{\wedge,\max})$ with $\xi\in \R$ form a one-parameter group of biholomorphisms. Let $\mathcal T_\wedge$ be the infinitesimal generator. The points where $\mathcal T_\wedge$ vanishes are the fixed points of $\pmb \kappa_{e^\xi}$. The vector field $\mathcal T_\wedge$ is the real part of a holomorphic vector field $\mathcal T_\wedge'$ (a holomorphic section  of $T^{1,0}\mathfrak D_\wedge$). Since $\mathcal T_\wedge$ vanishes at a point if and only if $\mathcal T_\wedge'$ vanishes at that point, we have that the set of fixed points of $\pmb \kappa_{e^\xi}$ in each $\Gr_{d_0}(\Sing_{\wedge,\max})$ is an analytic variety.

\begin{corollary}
The set of $\kappa$-invariant domains in $\mathfrak D_\wedge$ is an analytic variety.
\end{corollary}

Thus the set of $\kappa$-invariant domains is a small set.

\begin{remark}
By Lemma 5.12 of \cite{GiMe01}, a subspace of $\Sing_{\wedge,\max}$ is $\kappa$-invariant if and only if it is a direct sum of subspaces $\tilde \Sing_j\subset \Sing_{\sigma_j}(A_\wedge)$, each of which is itself $\kappa$-invariant. The set of $\kappa$-invariant subspaces of $\Sing_{\sigma_j}(A_\wedge)$ of a given dimension needs not be a discrete subset of the corresponding Grassmannian.
\end{remark}

Again as in Section \ref{sec-DomainsAndSpectra}, let
\begin{equation*}
\KK_{\wedge,\lambda}=\ker (A_{\wedge,\Dom_{\max}}-\lambda),\quad \lambda \in \bgres A_\wedge.
\end{equation*}
The proof of Proposition \ref{BundleOfKernels} gives that the $\KK_{\wedge,\lambda}$ are the fibers of a Hermitian holomorphic vector bundle
\begin{equation}\label{BundleOfKWedgeKernels}
\KK_\wedge\to \bgres A_\wedge
\end{equation}
over $\bgres A_\wedge$; the rank of $\KK_{\wedge}|_{\open \Lambda_\alpha}$ is the number $d'_\alpha$, which may change with $\alpha$.

\begin{lemma}\label{KappaOnKwedgeLambda}
Let $\lambda\in \bgres A_\wedge$. The map $\kappa_\varrho$ sends $\KK_{\wedge,\lambda}$ to $\KK_{\wedge,\varrho^m\lambda}$, and so gives a vector bundle morphism $\KK_\wedge\to\KK_\wedge$.
\end{lemma}

\begin{proof}
Writing \eqref{HomogeneityOfAwedgeBis} in the form
\begin{equation}\label{HomogeneityOfAwedgeBis2}
 A_\wedge - \lambda = \varrho^{-m} \kappa_\varrho^{-1}
 (A_\wedge -\varrho^m\lambda)\kappa_\varrho
 \;\text{ for every } \varrho>0,
\end{equation}
and letting each member of this identity
act on  $\phi\in\KK_{\wedge,\lambda}$, we see that $\kappa_\varrho
\phi\in\KK_{\wedge,\varrho^m\lambda}$.
\end{proof}

Lemma \ref{ResAndKernelBundle} has a word by word translation to the situation at hand and if $\Dom \in \mathfrak G_{\wedge,\alpha}$, then
\begin{equation}\label{ResAndWedgeKernelBundle}
\lambda\in \resolv A_{\wedge,\Dom}\cap \open\Lambda_\alpha\iff \KK_{\wedge,\lambda}\cap \Dom=0.
\end{equation}
For such $\lambda$,
\begin{equation}\label{KWedgeoplusDom}
\KK_{\wedge,\lambda}\oplus  \Dom=\Dom_{\wedge,\max}.
\end{equation}
Let
\begin{equation*}
\KK_{\wedge,\max}(\lambda)=\pi_{\wedge,\max}\KK_{\wedge,\lambda}.
\end{equation*}
Then
\begin{equation*}
\open\Lambda_\alpha\ni\lambda\to \KK_{\wedge,\max}(\lambda)\in \Gr_{d_\alpha'}(\Sing_{\wedge,\max})
\end{equation*}
is holomorphic.

Suppose that $\lambda_0\in \open \Lambda_\alpha$ and let $\Gamma=\set{r\lambda_0:r>0}$ be the ray through $\lambda_0$. In view of Lemma \ref{KappaOnKwedgeLambda}, the set $\spec A_{\wedge,\Dom}\cap \Gamma$ will not contain points $\lambda$ with $|\lambda|$ large if and only if $\kappa_\varrho \KK_{\wedge,\lambda_0}\cap \Dom=0$ for $\varrho$ large. With the notation introduced in Definition \ref{BadVariety} (of course with $\Sing_{\max}$ replaced by $\Sing_{\wedge,\max}$), this will happen if and only if $\pmb\kappa_\varrho \KK_{\wedge,\max}(\lambda_0) \notin \mathfrak V_{\pi_{\wedge,\max}\Dom}$ for large $\varrho$. Since $\mathfrak V_{\pi_{\wedge,\max}\Dom}\subset \Gr_{d_\alpha'}(\Sing_{\wedge,\max})$ is of complex codimension $1$ and
$\varrho\mapsto \pmb\kappa_\varrho\KK_{\wedge,\max}(\lambda_0)$ is a real curve, these curves generically do not intersect $\mathfrak V_{\pi_{\wedge,\max}\Dom}$. However, it can happen that $\pmb\kappa_\varrho \KK_{\wedge,\max}(\lambda_0)\in \mathfrak V_{\pi_{\wedge,\max}\Dom}$ for all $\varrho$, for instance if $\KK_{\wedge,\lambda_0}\cap \Dom$ contains a nontrivial $\kappa$-invariant subspace. It can also happen that $\pmb\kappa_\varrho \KK_{\wedge,\max}(\lambda_0)\in \mathfrak V_{\pi_{\wedge,\max}\Dom}$ infinitely often. For example, suppose that $\Sing_{\wedge,\max}$ is two-dimensional and that the infinitesimal generator of the action $\pmb \kappa_\varrho$ has two distinct eigenvalues $\im\sigma_1$ and $\im\sigma_2$ with $\Im\sigma_1=\Im\sigma_2$. Let $u_1$, $u_2$ be eigenvectors for these eigenvalues. If $a_1u_1+a_2u_2$ is a basis element for $\pi_{\wedge,\max}\Dom$ and $\KK_{\wedge,\max}(\lambda_0)$ is spanned by the same vector, then $\pmb \kappa_\varrho \KK_{\wedge,\max}(\lambda_0)$ is spanned by $a_1\varrho^{\im\sigma_1}u_1+a_2\varrho^{\im\sigma_2}u_2$, and $\pmb \kappa_\varrho\KK_{\wedge,\max}(\lambda_0)\cap \pi_{\wedge,\max}\Dom\ne 0$ whenever $\varrho=e^{2\pi k/(\Re\sigma_2-\Re\sigma_1)}$ with $k\in \mathbb Z$.

We will show that the spaces $\KK_{\wedge,\lambda}$, $\lambda\in \open\Lambda_\alpha$, can be obtained directly from a single space $\KK_{\wedge,\lambda_0}$, $\lambda_0\in \open\Lambda_\alpha$ via the action of $\kappa$ and $B_{\wedge,\min}(\lambda)$, the left inverse of $A_{\wedge,\Dom_{\min}}-\lambda$ with kernel equal to the orthogonal of $\RR_{\wedge,\lambda}=\rg (A_{\wedge,\Dom_{\min}}-\lambda)$. The family $B_{\wedge,\min}(\lambda)$ depends smoothly on $\lambda\in \bgres A_\wedge$, cf. Section \ref{sec-DomainsAndSpectra}.

Fix some sector $\open \Lambda_\alpha$ and for the sake of simplicity let $\lambda_0\in \open \Lambda_\alpha$ lie in the axis of symmetry $\Gamma_\alpha$ of $\open \Lambda_\alpha$. So
\begin{equation*}
\open \Lambda_\alpha=\set{\lambda:\,|\arg(\lambda/\lambda_0)|<\theta_\alpha}
\end{equation*}
where $\arg$ is the principal branch of the argument function on $\C\minus\overline \R_-$. Let $\log$ be the principal branch of the logarithm on the same set. Then
\begin{equation*}
\mathfrak P_{\max}(\lambda)\phi = \pi_{\wedge,\max}\kappa_{e^{\log(\lambda/\lambda_0)/m}} \pi_{\wedge,\max}\phi
\end{equation*}
is well defined for $\phi\in \KK_{\wedge,\lambda_0}$ and is holomorphic in $\lambda$ for $\lambda\notin -\Gamma_\alpha$. Thus we have a map
\begin{equation*}
\mathfrak P_{\max}(\lambda):\KK_{\wedge,\lambda_0}\to \Sing_{\wedge,\max}
\end{equation*}
depending holomorphically on $\lambda$ for $\lambda\notin -\Gamma_\alpha$. In general, if $\phi\in \KK_{\wedge,\lambda}$, then
\begin{equation*}
(A_\wedge-\lambda)\pi_{\wedge,\min}\phi=-(A_\wedge-\lambda)\pi_{\wedge,\max}\phi.
\end{equation*}
Thus the right hand side belongs to the range $\RR_{\wedge,\lambda}$ of $A_{\wedge,\Dom_{\min}}-\lambda$, and
\begin{equation*}
\pi_{\wedge,\min}\phi=-B_{\wedge,\min}(\lambda)(A_\wedge-\lambda)\pi_{\wedge,\max}\phi
\end{equation*}
Conversely, if $u\in \Sing_{\wedge,\max}$ and $(A_\wedge-\lambda)u\in \RR_{\wedge,\lambda}$, then
\begin{equation*}
u - B_{\wedge,\min}(\lambda)(A_\wedge-\lambda)u\in \KK_{\wedge,\lambda}.
\end{equation*}
Define
\begin{equation*}
\mathfrak P_{\min}(\lambda):\KK_{\wedge,\lambda_0}\to\Dom_{\wedge,\min},\quad \lambda\in \open \Lambda_\alpha
\end{equation*}
by
\begin{equation*}
\mathfrak  P_{\min}(\lambda)\phi = -B_{\wedge,\min}(\lambda)(A_\wedge-\lambda)\mathfrak P_{\max}(\lambda)\phi,\quad \phi\in \KK_{\wedge,\lambda_0}.
\end{equation*}
Let also
\begin{equation}\label{HolomorphicExtensions}
\mathfrak P(\lambda)=\mathfrak P_{\max}(\lambda)+\mathfrak P_{\min}(\lambda):\KK_{\wedge,\lambda_0}\to\Dom_{\wedge,\max}
\end{equation}
The operators $B_{\wedge,\min}(\lambda)$ depend smoothly, but not holomorphically, on $\lambda$ (unless $\Ind (A_{\wedge,\Dom_{\min}}-\lambda)=0$ for $\lambda\in \open \Lambda_\alpha$). So it is not obvious that $\mathfrak  P_{\min}(\lambda)$ depends holomorphically on $\lambda$.

\begin{proposition}\label{HolomorphicExtensionsSplittingProp}
The map $\mathfrak P_{\min}(\lambda) : \KK_{\wedge,\lambda_0} \to \Dom_{\wedge,\min}$ depends holomorphically on $\lambda\in \open \Lambda_\alpha$, and
\begin{equation}\label{HolomorphicExtensionsSplitting}
\mathfrak P(\lambda)\phi\in \KK_{\wedge,\lambda}\quad \text{for }\lambda\in \open \Lambda_\alpha \text{ and } \phi\in \KK_{\wedge,\lambda_0}.
\end{equation}
\end{proposition}

\begin{proof}
Since $\lambda_0\in \bgres A_\wedge$, there is $\Dom \in \mathfrak G_{\wedge,\alpha}$ such that $\lambda_0\in \resolv A_{\wedge,\Dom}$, so $\resolv A_{\wedge,\Dom}\cap \open \Lambda_\alpha\ne \varnothing$. Let then  $B_{\wedge,\Dom}(\lambda)$ denote the resolvent of $A_{\wedge,\Dom}$ on $\resolv A_{\wedge,\Dom}\cap \open \Lambda_\alpha$. In particular $B_{\wedge,\Dom}(\lambda)$ depends  holomorphically on $\lambda\in \resolv A_{\wedge,\Dom}\cap \open \Lambda_\alpha$. Let $\phi\in \KK_{\wedge,\lambda_0}$ and define
\begin{equation*}
\mathfrak  P_{\min,\Dom}(\lambda)\phi = -\pi_{\wedge,\min}B_{\wedge,\Dom}(\lambda)(A_\wedge-\lambda)\mathfrak P_{\max}(\lambda)\phi
\end{equation*}
for $\lambda\in \resolv A_{\wedge,\Dom}\cap \open \Lambda_\alpha$. Thus $\mathfrak  P_{\min,\Dom}(\lambda)$ is holomorphic on the set where we defined it. Note that $\pi_{\wedge,\min}B_{\wedge,\Dom}(\lambda)$ is a holomorphic left inverse for $A_{\wedge,\Dom_{\min}}-\lambda$ when $\lambda\in \resolv A_{\wedge,\Dom} \cap \open \Lambda_\alpha$.

If $\phi\in \KK_{\wedge,\lambda_0}$ and $\varrho\in \R_+$ then $\kappa_\varrho\phi \in \KK_{\wedge,\varrho^m\lambda_0}$ and
\begin{equation*}
\mathfrak P_{\max}(\varrho^m\lambda_0)\phi = \pi_{\wedge,\max}\kappa_\varrho\phi.
\end{equation*}
Thus if $\varrho^m\lambda_0\in \resolv A_{\wedge,\Dom}$ then
\begin{align*}
\pi_{\wedge,\min}\kappa_\varrho\phi &=-\pi_{\wedge,\min}B_{\wedge,\Dom}(\varrho^m\lambda_0) (A_\wedge-\varrho^m\lambda_0)\pi_{\wedge,\max}\kappa_\varrho \phi\\
&=-\pi_{\wedge,\min}B_{\wedge,\Dom}(\varrho^m\lambda_0)(A_\wedge-\varrho^m\lambda_0)\mathfrak P_{\max}(\varrho^m\lambda_0)\phi.
\end{align*}
Consequently,
\begin{equation*}
\mathfrak P_{\max}(\varrho^m\lambda_0)\phi+\mathfrak P_{\min,\Dom}(\varrho^m\lambda_0)\phi=\kappa_{\varrho}\phi\in \KK_{\wedge,\varrho^m\lambda_0}.
\end{equation*}
This implies that the equation
\begin{equation*}
(A_\wedge-\lambda)[\mathfrak P_{\max}(\lambda)\phi+\mathfrak P_{\min,\Dom}(\lambda)\phi]=0
\end{equation*}
is satisfied when $\lambda\in \Gamma_\alpha\cap \resolv A_{\wedge,\Dom}$. By unique continuation, it is satisfied for any $\lambda\in \open \Lambda_\alpha\cap \resolv A_{\wedge,\Dom}$. Thus
\begin{equation}\label{HolomorphicExtensionsPrelim}
\mathfrak P_{\max}(\lambda)\phi+\mathfrak P_{\min,\Dom}(\lambda)\phi
\in \KK_{\wedge,\lambda}\ \text{ if }\lambda\in \resolv A_{\wedge,\Dom}\cap \open \Lambda_\alpha.
\end{equation}
For such $\lambda$ we therefore have
\begin{equation*}
(A_\wedge-\lambda)\mathfrak P_{\max}(\lambda)\phi = -(A_\wedge-\lambda)\mathfrak P_{\min,\Dom}(\lambda)\phi
\end{equation*}
so
\begin{equation}
\mathfrak P_{\min,\Dom}(\lambda)\phi=-B_{\wedge,\min}(\lambda)(A_\wedge-\lambda)\mathfrak P_{\max}(\lambda)\phi
\end{equation}
that is,
\begin{equation}\label{HolomorphicExtensionsHolCont}
\mathfrak P_{\min,\Dom}(\lambda)\phi=\mathfrak P_{\min}(\lambda)\phi.
\end{equation}
Replacing this in \eqref{HolomorphicExtensionsPrelim} we see that formula \eqref{HolomorphicExtensionsSplitting} holds where $\mathfrak P_{\min,\Dom}(\lambda)\phi$ is defined.

Since the left hand side of \eqref{HolomorphicExtensionsHolCont} is holomorphic where defined, so is the right hand side. Since the right hand side is continuous on $\open \Lambda_\alpha$, the singularities of the left hand side, {\ie} the elements of the discrete set $\spec A_{\wedge,\Dom}\cap\open \Lambda_\alpha$, are removable. Thus $\mathfrak P_{\min}(\lambda)\phi$ is holomorphic for $\lambda\in \open \Lambda_\alpha$ and by continuity \eqref{HolomorphicExtensionsSplitting} holds.
\end{proof}

The sets
\begin{equation*}
\mathcal L_\alpha=\set{\KK_{\wedge,\max}(\lambda) : \lambda \in \open\Lambda_\alpha}\subset \Gr_{d'_\alpha}(\Sing_{\wedge,\max})
\end{equation*}
play an important role, particularly their intersection with the varieties $\mathfrak V_{\pi_{\wedge,\max}\Dom}$, $\Dom\in \mathfrak G_{\wedge,\alpha}$. These sets are invariant under the action of $\pmb \kappa_{e^\xi}$ for $\xi$ real. If $\Vee\in\Gr_{d_0}(\Sing_{\wedge,\max})$, then $\C\ni\zeta\mapsto \pmb \kappa_{e^\zeta}\Vee\in \Gr_{d_0}(\Sing_{\wedge,\max})$ is a holomorphic map, cf. Proposition \ref{KappaVectorField}. The generator of the one-parameter group $(\xi,\Vee)\mapsto \kappa_{e^\xi}\Vee$, the vector field $\mathcal T_\wedge$, is the real part of a holomorphic vector field $\mathcal T_\wedge'$, cf. the paragraph following the proof of Proposition \ref{KappaVectorField}, which at $\Vee$ is the image of the Cauchy-Riemann vector field at $\zeta=0$, $\partial_\zeta|_0$,  under the differential of the map $\zeta\mapsto \pmb \kappa_{e^\zeta}\Vee$. If $\Vee$ is not $\kappa$-invariant, \ie, $\mathcal T_\wedge\ne 0$ at $\Vee$, then also the imaginary part of $\mathcal T_\wedge'$ is different from $0$ at $\Vee$; thus $\zeta\mapsto \pmb \kappa_{e^\zeta}\Vee$ is a local embedding near $\zeta=0$ if $\mathcal T_\wedge\ne 0$ at $\Vee$. As a consequence we get that the real and imaginary parts of $\mathcal T_\wedge'$ commute at the noninvariant points. Since the set of invariant points is closed with empty interior, the real and imaginary parts of $\mathcal T_\wedge'$ commute everywhere. We can view the images of the maps $\C\ni \zeta\mapsto \pmb \kappa_{e^\zeta}\Vee$ as a point (if $\Vee$ is invariant) or as an integral manifold of the involutive Frobenius distribution generated by $\Re\mathcal T_\wedge'$, $\Im\mathcal T_\wedge'$ on $\Gr_{d_0}(\Sing_{\wedge,\max})\minus \set{\Vee\in \Gr_{d_0}:\Vee\text{ is } \kappa \text{-invariant}}$.

\begin{theorem}
The set $\mathcal L_\alpha$ is contained in one orbit of $\mathcal T_\wedge'$.
\end{theorem}

In fact, the set $\mathcal L_\alpha$ is identical to the set
\begin{equation*}
\set{ \kappa_{e^{\log(\lambda/\lambda_0)/m}}\pi_{\wedge,\max}\KK_{\wedge,\lambda_0}: \lambda \in \open\Lambda_\alpha},
\end{equation*}
a subset of the orbit of $\mathcal T_\wedge'$ containing $\pi_{\wedge,\max}\KK_{\wedge,\lambda_0}$. Thus, if $\dim_\R\Gr_{d'_\alpha}(\Sing_{\wedge,\max})>2$, then $\mathcal L_\alpha$ is in principle a small set (nevertheless it could be dense).

The following lemma completes our description of the vector bundle $\KK_{\wedge}$ introduced in \eqref{BundleOfKWedgeKernels}.

\begin{lemma}
If $\phi\in \KK_{\wedge,\lambda_0}$, then
\begin{equation}\label{HolomorphicExtensionsHomogeneity}
\kappa_\varrho\mathfrak P(\lambda)\phi=\mathfrak P(\varrho^m\lambda)(\phi), \quad \varrho\in \R_+.
\end{equation}
\end{lemma}

\begin{proof}
Write $(\lambda/\lambda_0)^{1/m}=e^{[\log(|\lambda/\lambda_0|)+\im\arg(\lambda/\lambda_0)]/m}$, $\lambda\notin -\Gamma_\alpha$. For real $\varrho$ and $\lambda\notin-\Gamma_\alpha$ we have
\begin{equation*}
\pi_{\wedge,\max}\kappa_\varrho \pi_{\wedge,\max}\kappa_{(\lambda/\lambda_0)^{1/m}}\pi_{\wedge,\max} =
\pi_{\wedge,\max}\kappa_{(\varrho^m\lambda/\lambda_0)^{1/m}}\pi_{\wedge,\max}
\end{equation*}
Thus $\pi_{\wedge,\max}\kappa_{\varrho}\mathfrak P_{\max}(\lambda)\phi=\mathfrak P_{\max}(\varrho^m\lambda)\phi$. But if $\lambda\in \open\Lambda_\alpha$, then $\kappa_\varrho\mathfrak P(\lambda)\phi\in \KK_{\wedge,\varrho^m\lambda}$ and  $\mathfrak P_{\min}(\varrho^m\lambda)\phi$ is the unique element of $\Dom_{\wedge,\min}$ such that
\begin{equation*}
\mathfrak P_{\max}(\varrho^m\lambda)\phi+\mathfrak P_{\min}(\varrho^m\lambda)\phi \in \KK_{\wedge,\varrho^m\lambda}
\end{equation*}
so we have
\begin{equation*}
\kappa_\varrho \mathfrak P(\lambda)\phi=\mathfrak P_{\max}(\varrho^m\lambda)\phi+\mathfrak P_{\min}(\varrho^m\lambda)\phi.
\end{equation*}
This is \eqref{HolomorphicExtensionsHomogeneity}.
\end{proof}

Suppose  $\phi\in \KK_{\wedge,\lambda_0}$. If $\phi(\lambda)=\mathfrak P(\lambda)\phi$ vanishes at some $\lambda_1$, then $\phi(\lambda)$ vanishes along the ray through $\lambda_1$. Therefore it vanishes identically, since $\phi(\lambda)$ is holomorphic. Thus, if we pick a basis $\set{\phi_j}$ of $\KK_{\wedge,\lambda_0}$, then the sections $\mathfrak P(\lambda)\phi_j$ form a frame over $\open \Lambda_\alpha$ for the bundle $\KK_\wedge$.

%%%%%%%%%%%%%%%%%%%%%%%%%%%%%%%%%%%%%%%%%%%%%%%%%%%%%%%%%%%
%%%%%%%%%%%%%%%%%%%%%%%%%%%%%%%%%%%%%%%%%%%%%%%%%%%%%%%%%%%
\section{Resolvents for the model operator}
\label{sec-ResolventsOnYWedge}

We now turn our attention to determining the existence of sectors of minimal growth for extensions of $A_\wedge$.

If $\Dom\in \mathfrak G_{\wedge,\alpha}$, we write $B_{\wedge,\Dom}(\lambda)$ for the inverse of $A_{\wedge,\Dom}-\lambda$, $\lambda\in \resolv A_{\wedge,\Dom}$. If $\lambda\in \C\minus 0$, then $\hat\lambda=\lambda/|\lambda|$. By a closed sector we shall mean a set of the form
\begin{equation*}
\Lambda = \set{z\in\C : z=re^{i\theta} \text{ for } r\geq 0,\
 \theta\in \R, \ |\theta-\theta_0|\leq  a}.
\end{equation*}
If $R>0$ and $\Lambda$ is a closed sector, then  $\Lambda_R=\set{\lambda\in \Lambda : |\lambda|\geq R}$. Let $\Dom\in \mathfrak G_{\wedge,\alpha}$. Let $\Lambda$ be a closed sector with $\Lambda\minus 0 \subset \open\Lambda_\alpha$. Then $\Lambda$ is called a sector of minimal growth for $A_{\wedge,\Dom}$ if there is $R>0$ such that $A_{\wedge,\Dom}-\lambda$ is invertible if $\lambda\in \Lambda_R$, and either of the equivalent estimates
\begin{equation}\label{ResolventEstimatesWedge}
 \|B_{\wedge,\Dom}(\lambda)\|_{\L(x^{-m/2}L^2_b)}
 \leq C/|\lambda|, \quad
 \|B_{\wedge,\Dom}(\lambda)\|_{\L(x^{-m/2}L^2_b,
 \Dom_{\wedge,\max})} \leq C
\end{equation}
holds for some $C>0$ when $\lambda \in \Lambda_R$.

The following lemma is immediate, in view of \eqref{HomogeneityOfAwedgeBis} and the fact that $\kappa_\varrho$ is an isometry on $x^{-m/2}L^2_b(Y^\wedge;E)$.

\begin{lemma}\label{GeneralKappaDomWedge}
If $\Dom\in \mathfrak G_{\wedge,\alpha}$, then $\resolv A_{\wedge,\kappa_\varrho^{-1}\Dom}= \varrho^{-m}\resolv A_{\wedge,\Dom}$. If $\resolv A_{\wedge,\Dom}\ne \varnothing$, then
\begin{equation}\label{TransfromationRuleWedge}
B_{\wedge,\kappa_\varrho^{-1}\Dom}(\lambda)= \varrho^m \kappa_\varrho^{-1} B_{\wedge,\Dom}(\varrho^m\lambda)\kappa_{\varrho},\quad \varrho>0.
\end{equation}
Thus, if the closed sector $\Lambda$, $\Lambda\minus 0\subset \open\Lambda_\alpha$, is a sector of minimal growth for $A_{\wedge,\Dom}$, then $\Lambda$ is a sector of minimal growth for $A_{\wedge,\kappa_\varrho^{-1}\Dom}$.
\end{lemma}

In fact, if the first estimate in \eqref{ResolventEstimatesWedge} holds when $\lambda\in \Lambda_R$, then
\begin{equation*}
\|B_{\wedge,\kappa_\varrho^{-1}\Dom}(\lambda)\|_{\L(x^{-m/2}L^2_b)} \leq C/|\lambda|,\quad \lambda\in \Lambda_{R/\varrho^m}
\end{equation*}
with the same constant $C$.

The simplest domains are those that are $\kappa$-invariant:

\begin{proposition}\label{ResOnInvWedge}
Suppose $\Dom\in \mathfrak G_{\wedge,\alpha}$ is $\kappa$-invariant. Then either $\open \Lambda_\alpha\cap \resolv A_{\wedge,\Dom}=\varnothing$, or $\open \Lambda_\alpha\subset \resolv A_{\wedge,\Dom}$.  In the latter case the resolvent $B_{\wedge,\Dom}(\lambda)$ of $A_{\wedge,\Dom}$ satisfies
\begin{equation}\label{MinimalGrowthForInvariant}
\|B_{\wedge,\Dom}(\lambda)\|_{\L(x^{-m/2}L^2_b)}\leq C/|\lambda|
\end{equation}
for some $C>0$ when $\lambda\in \Lambda\minus 0$, $\Lambda$ a closed sector with $\Lambda\minus 0\subset \open \Lambda_\alpha$.
\end{proposition}
\begin{proof}
Suppose that $\lambda_0\in \open \Lambda_\alpha\cap\spec A_{\wedge,\Dom}$. The homogeneity property \eqref{HomogeneityOfAwedgeBis} and the assumption that $\Dom$ is $\kappa$-invariant  give that $\lambda/\varrho^m\in \spec A_{\wedge,\Dom}$ for every $\varrho>0$. Thus $\spec A_{\wedge,\Dom}\cap \open \Lambda_\alpha$ is not discrete.
On the other hand, if $\Dom\in \mathfrak G_{\wedge,\alpha}\minus \mathfrak V_\alpha$, cf. Proposition~\ref{BadSpectralVarietyWedge}, then $\spec A_{\wedge,\Dom}\cap \open \Lambda_\alpha$ is a discrete closed subset of $\open \Lambda_\alpha$. In particular, for every ray
\begin{equation}\label{TheRay}
 \Gamma=\set{z\in\C : z=r e^{i\theta_0} \text{ for } r>0}
\end{equation}
contained in $\open \Lambda_\alpha$, $\Gamma\cap \spec A_{\wedge,\Dom}$ is closed and discrete.
Thus, if $\Dom$ is $\kappa$-invariant, then $\open \Lambda_\alpha\cap\spec A_{\wedge,\Dom}\ne \varnothing$ implies $\open \Lambda_\alpha\subset \spec A_{\wedge,\Dom}$.

Suppose $\Lambda$ is a closed sector with $\Lambda\minus 0 \subset \resolv A_{\wedge,\Dom}$. Since $\Dom$ is $\kappa$-invariant, \eqref{TransfromationRuleWedge} reads
\begin{equation*}
B_{\wedge,\Dom}(\lambda)= \varrho^m \kappa_\varrho^{-1} B_{\wedge,\Dom}(\varrho^m\lambda)\kappa_{\varrho}.
\end{equation*}
Setting $\varrho=|\lambda|^{-1/m}$  gives
\begin{equation*}
B_{\wedge,\Dom}(\lambda)= |\lambda|^{-1} \kappa_{|\lambda|^{1/m}} B_{\wedge,\Dom}(\hat\lambda)\kappa_{|\lambda|^{1/m}}^{-1}.
\end{equation*}
For $\hat\lambda\in \Lambda$ ($|\hat\lambda|=1$) we have a uniform estimate for $\|B_{\wedge,\Dom}(\hat\lambda)\|_{\L(x^{-m/2}L^2_b)}$, and \eqref{MinimalGrowthForInvariant} follows immediately, since $\kappa_\varrho$ is an isometry on $x^{-m/2}L^2_b(Y^\wedge;E)$.
\end{proof}

If the domain $\Dom\in \mathfrak G_{\wedge,\alpha}$ is not $\kappa$-invariant, the existence of a ray or sector of minimal growth for $B_{\wedge,\Dom}(\lambda)$ is more complicated:

\begin{theorem}\label{NecAndSuffWedge}
Let $\Dom\in \mathfrak G_{\wedge,\alpha}$, let $\Lambda$ be a closed sector with $\Lambda\minus0\subset \open\Lambda_\alpha$. Then $\Lambda$ is a sector of minimal growth for $A_{\wedge,\Dom}$ if and only if there are $C$, $R>0$ such that $\Lambda_R\subset \resolv A_{\wedge,\Dom}$ and
\begin{equation}\label{NecAndSuffWedgeCond}
\|\pi_{\wedge,\max}\pi_{\KK_{\wedge,\hat\lambda},\kappa_{|\lambda|^{1/m}}^{-1}\Dom}\big|_{\Sing_{\wedge,\max}}\|_{\L(\Dom_{\wedge,\max})}\leq C,\quad \lambda\in \Lambda_R.
\end{equation}
\end{theorem}

If $\Dom$ is $\kappa$-invariant, then
\begin{equation*}
\pi_{\wedge,\max}\pi_{\KK_{\wedge,\hat\lambda}, \kappa_{|\lambda|^{1/m}}^{-1}\Dom}\big|_{\Sing_{\wedge,\max}} = \pi_{\wedge,\max}\pi_{\KK_{\wedge,\hat\lambda}, \Dom}\big|_{\Sing_{\wedge,\max}},
\end{equation*}
and the theorem reduces to the trivial situation of Proposition \ref{ResOnInvWedge}.

The proof of the theorem requires some preparation. Define $B_{\wedge,\max}(\lambda)$ for $\lambda\in \bgres A_\wedge$ as the right inverse of
\begin{equation*}
A_\wedge-\lambda:\Dom_{\wedge,\max}\subset x^{-m/2}L^2_b(Y^\wedge;E) \to x^{-m/2}L^2_b(Y^\wedge;E)
\end{equation*}
with range in $\KK_{\wedge,\lambda}^\perp$. Thus $B_{\wedge,\max}(\lambda)$ has the virtue of being the right inverse of $A_{\wedge,\max}-\lambda$ with the smallest operator norm. It depends smoothly on $\lambda\in \bgres A_\wedge$; this is proved in the same way as the corresponding statement for $B_{\max}(\lambda)$ in Section~\ref{sec-DomainsAndSpectra}. Recall that $B_{\wedge,\min}(\lambda)$ is the left inverse of $A_{\wedge,\Dom_{\min}}-\lambda$ with kernel equal to the orthogonal of $\RR_{\wedge,\lambda}=\rg (A_{\wedge,\Dom_{\min}}-\lambda)$.
For $\lambda\in \resolv A_{\wedge,\Dom}$ let
\begin{equation*}
\pi_{\KK_{\wedge,\lambda},\Dom}:\Dom_{\wedge,\max}\to \Dom_{\wedge,\max}
\end{equation*}
be the projection on $\KK_{\wedge,\lambda}$ according to the decomposition $\Dom_{\wedge,\max}=\KK_{\wedge,\lambda}\oplus \Dom$, cf. \eqref{KWedgeoplusDom}. Then the resolvent of $A_{\wedge,\Dom}$ is
\begin{equation}\label{FormulaForBWedgeDomFirst}
B_{\wedge,\Dom}(\lambda)= B_{\wedge,\max}(\lambda)-\big(I- B_{\wedge,\min}(\lambda)(A_\wedge-\lambda)\big)\pi_{\KK_{\wedge,\lambda},\Dom} B_{\wedge,\max}(\lambda),
\end{equation}
cf. \eqref{ResolventOnDviaMinMaxTmp}. We will take advantage of this formula by using the group action $\kappa$. We begin with estimates for $B_{\wedge,\min}(\lambda)$ and $B_{\wedge,\max}(\lambda)$.

\begin{lemma}
The operator $B_{\wedge,\min}(\lambda)$ is $\kappa$-homogeneous of degree $-m$,
\begin{equation}\label{FormulaForBWedgeMin}
B_{\wedge,\min}(\lambda) = \varrho^{-m}\kappa_\varrho B_{\wedge,\min}(\lambda/\varrho^m)\kappa_\varrho^{-1}.
\end{equation}
Therefore, if $\Lambda$ is a closed sector with $\Lambda\minus 0\subset \bgres A_\wedge$, then
\begin{equation}\label{EsimateForBWedgeMin}
\|B_{\wedge,\min}(\lambda)\|_{\L(x^{-m/2}L^2_b)}\leq C/|\lambda|
\end{equation}
for some $C>0$ when $\lambda\in \Lambda\minus 0$.
\end{lemma}
\begin{proof}
Let $B'_{\wedge,\min}(\lambda)=\varrho^{-m}\kappa_\varrho B_{\wedge,\min}(\lambda/\varrho^m)\kappa_\varrho^{-1}$. The operator $B'_{\wedge,\min}(\lambda)$ maps into $\Dom_{\wedge,\min}$ because the latter space is $\kappa$-invariant. Using the $\kappa$-homogeneity of $A_\wedge$ one verifies that the operator $B'_{\wedge,\min}(\lambda)$ is a left inverse for $A_{\wedge,\min}-\lambda$. Also because of the $\kappa$-invariance of $\Dom_{\wedge,\min}$ and the $\kappa$-homogeneity of $A_\wedge-\lambda$, $\RR_{\lambda} = \kappa_\varrho \RR_{\lambda/\varrho^m}$. The kernel of $B'_{\wedge,\min}(\lambda)$ is $\kappa_{\varrho}\RR_{\lambda/\varrho^m}^\perp$. Since $\kappa_\varrho$ is an isometry on $x^{-m/2}L^2_b(Y^\wedge;E)$, $\kappa_\varrho$ preserves the orthogonality of the decomposition $\RR_{\lambda/\varrho^m}\oplus \RR_{\lambda/\varrho^m}^\perp$. Hence, the kernel of $B'_{\wedge,\min}(\lambda)$ is orthogonal to $\RR_{\lambda}$. Thus $B'_{\wedge,\min}(\lambda) = B_{\wedge,\min}(\lambda)$, and \eqref{FormulaForBWedgeMin} holds.

The estimate in \eqref{EsimateForBWedgeMin} follows from setting $\varrho^m=|\lambda|$ in \eqref{FormulaForBWedgeMin}.
\end{proof}

The operator family $B_{\wedge,\max}(\lambda)$ is not $\kappa$-homogeneous. Nevertheless its norm satisfies good estimates.

\begin{proposition}\label{PropOnFormulaForBWedgeMax}
Let $\pi_{\KK_{\wedge,\lambda}}:\Dom_{\wedge,\max}\to \Dom_{\wedge,\max}$ be the orthogonal projection on $\KK_{\wedge,\lambda}$.  Regard the finite dimensional space $\KK_{\wedge,\lambda}$ as a subspace of $x^{-m/2}L^2_b(Y^\wedge;E)$ and let
\begin{equation*}
\mathfrak p_{\KK_{\wedge,\lambda}}:x^{-m/2}L^2_b(Y^\wedge;E)\to x^{-m/2}L^2_b(Y^\wedge;E)
\end{equation*}
be the orthogonal projection on $\KK_{\wedge,\lambda}$. Then
\begin{equation}\label{FormulaBWedgeMax}
B_{\wedge,\max}(\lambda)=
|\lambda|^{-1} \kappa_{|\lambda|^{1/m}}\big(\id -  \frac{1-|\lambda|^2}{1+|\lambda|^2} \mathfrak p_{\KK_{\wedge,\hat\lambda}}\big) B_{\wedge,\max}(\hat\lambda)\kappa_{|\lambda|^{1/m}}^{-1}.
\end{equation}
Therefore, if $\Lambda$ is a closed sector such that $\Lambda\minus 0\subset \bgres A_\wedge$, then
\begin{equation}\label{EstimateForBWedgeMax}
\|B_{\wedge,\max}(\lambda)\|_{\L(x^{-m/2}L^2_b)} \leq C/|\lambda|
\end{equation}
for some $C>0$ when $\lambda\in \Lambda\minus 0$.
\end{proposition}

The proof will require:
\begin{lemma}\label{LemmaOnProjectionsWedge}
For any $\lambda\in \bgres A_\wedge$ and $\varrho>0$,
\begin{equation}\label{RelationBetweenPisWedge}
\kappa_\varrho^{-1} \pi_{\KK_{\wedge,\varrho^m\lambda}}\kappa_\varrho =
\frac{1+|\lambda|^2}{1+|\varrho^m\lambda|^2}\varrho^{2m} \pi_{\KK_{\wedge,\lambda}} + \frac{1-\varrho^{2m}}{1+|\varrho^m\lambda|^2}\mathfrak p_{\KK_{\wedge,\lambda}}.
\end{equation}
\end{lemma}
We will prove the lemma later.

\begin{proof}[Proof of Proposition \ref{PropOnFormulaForBWedgeMax}]
Suppose $f\in x^{-m/2}L^2_b(Y^\wedge;E)$ and let $u=B_{\wedge,\max}(\lambda)f$. Then $(A_\wedge-\varrho^m\lambda)\kappa_\varrho u=\varrho^m\kappa_\varrho(A_\wedge-\lambda)u= \varrho^m\kappa_\varrho f$, and consequently
\begin{equation*}
B_{\wedge,\max}(\varrho^m\lambda)\varrho^m\kappa_\varrho f
=\kappa_\varrho u -\pi_{\KK_{\wedge,\varrho^m\lambda}} \kappa_\varrho u.
\end{equation*}
This gives the formula
\begin{equation*}
\varrho^m B_{\wedge,\max}(\varrho^m\lambda)\kappa_\varrho f
=\kappa_\varrho B_{\wedge,\max}(\lambda)f -\pi_{\KK_{\wedge,\varrho^m\lambda}} \kappa_\varrho B_{\wedge,\max}(\lambda)f
\end{equation*}
which, in view of \eqref{RelationBetweenPisWedge} and the fact that the range of $B_{\wedge,\max}(\lambda)$ is orthogonal to $\KK_{\wedge,\lambda}$, reduces to
\begin{equation*}
\varrho^m B_{\wedge,\max}(\varrho^m\lambda)\kappa_\varrho f
=\kappa_\varrho B_{\wedge,\max}(\lambda)f -\frac{1-\varrho^{2m}}{1+|\varrho^m\lambda|^2}\kappa_\varrho\mathfrak p_{\KK_{\wedge,\lambda}} B_{\wedge,\max}(\lambda)f.
\end{equation*}
Thus
\begin{equation*}
B_{\wedge,\max}(\varrho^m\lambda)
=\varrho^{-m} \kappa_\varrho \big(\id - \frac{1-\varrho^{2m}}{1+|\varrho^m\lambda|^2}\mathfrak p_{\KK_{\wedge,\lambda}}\big)  B_{\wedge,\max}(\lambda)\kappa_\varrho^{-1}.
\end{equation*}
The formula \eqref{FormulaBWedgeMax} is obtained from this by replacing $\varrho^m$ by $|\lambda|$ and $\lambda$ by $\hat\lambda$. The estimate \eqref{EstimateForBWedgeMax} is evident given the formula \eqref{FormulaBWedgeMax}.
\end{proof}

The operator
\begin{equation*}
B_{\wedge,\max}^h(\lambda)=|\lambda|^{-1} \kappa_{|\lambda|^{1/m}} B_{\wedge,\max}(\hat\lambda)\kappa_{|\lambda|^{1/m}}^{-1},\quad \lambda\in \bgres A_\wedge
\end{equation*}
is a $\kappa$-homogeneous right inverse of $A_{\wedge,\max}-\lambda$ of degree $-m$ that coincides with $B_{\wedge,\max}(\lambda)$ when $|\lambda|=1$. For any closed sector $\Lambda$ with $\Lambda\minus 0\subset \bgres A_\wedge$ there is $C$ such that
\begin{equation*}
\|B_{\wedge,\max}^h(\lambda)\|_{\L(x^{-m/2}L^2_b)}\leq C/|\lambda|,\quad \lambda\in \Lambda\minus 0,
\end{equation*}
and for any closed sector $\Lambda$ as above and $R>0$,
\begin{equation*}
\|B_{\wedge,\max}(\lambda)-B_{\wedge,\max}^h(\lambda)\|_{\L(x^{-m/2}L^2_b)}\leq C/|\lambda|
\end{equation*}
for $\lambda\in \Lambda_R$. So in some estimates below, it makes little difference whether the correction term involving $\mathfrak p_{\KK_{\wedge,\hat\lambda}}$ is present or not. However, we will keep on using $B_{\wedge,\max}(\lambda)$ instead of $B_{\wedge,\max}^h(\lambda)$, as the former family is in some sense more natural than the latter.

\begin{proof}[Proof of Lemma \ref{LemmaOnProjectionsWedge}]
Let $\phi_1,\dotsc,\phi_{d_\alpha'}$ be an $A_\wedge$-orthonormal basis of
$\KK_{\wedge,\lambda}$. Then
\begin{equation*}
 \delta_{jk}=(\phi_j,\phi_k)_{A_\wedge} = (1+|\lambda|^2)(\phi_j,\phi_k).
\end{equation*}
In particular, the $\sqrt{1+|\lambda|^2}\phi_j\in \KK_{\wedge,\lambda}$ are
orthonormal in $x^{-m/2}L^2_b(Y^\wedge;E)$. On the other hand, using that
$\kappa_\varrho$ is an isometry on $x^{-m/2}L^2_b(Y^\wedge;E)$,
\begin{multline*}
(\kappa_\varrho\phi_j,\kappa_\varrho\phi_k)_{A_\wedge}
= \varrho^{2m}(\phi_j,\phi_k)_{A_\wedge} + (1-\varrho^{2m})(\phi_j,\phi_k)\\
= \Big(\varrho^{2m} + \frac{1-\varrho^{2m}}{1+|\lambda|^2}\Big)\delta_{jk}
= \frac{1+|\varrho^m \lambda|^2}{1+|\lambda|^2}\,\delta_{jk}.
\end{multline*}
Thus the $\sqrt{(1+|\lambda|^2)/(1+|\varrho^m \lambda|^2)} \kappa_\varrho
\phi_j\in \KK_{\wedge,\varrho^m\lambda}$ are $A_\wedge$-orthonormal, and if
$u\in \KK_{\wedge,\lambda}$, then
\begin{align*}
\pi_{\KK_{\wedge,\varrho^m\lambda}} \kappa_\varrho u
&= \frac{1+|\lambda|^2}{1+|\varrho^m \lambda|^2}
\sum_j(\kappa_\varrho u,\kappa_\varrho\phi_j)_{A_\wedge}\kappa_\varrho\phi_j\\
&= \frac{1+|\lambda|^2}{1+|\varrho^m \lambda|^2}
\sum_j\big[\varrho^{2m}(u,\phi_j)_{A_\wedge} + (1-\varrho^{2m})(u,\phi_j)
\big]\kappa_\varrho\phi_j \\
&= \frac{1+|\lambda|^2}{1+|\varrho^m \lambda|^2}\kappa_\varrho
\sum_j\big[\varrho^{2m}(u,\phi_j)_{A_\wedge}\phi_j
+ (1-\varrho^{2m})(u,\phi_j)\phi_j \big] \\
&= \kappa_{\varrho}\big(\frac{1+|\lambda|^2}{1+|\varrho^m\lambda|^2}\varrho^{2m}
\pi_{\KK_{\wedge,\lambda}}u
+ \frac{1-\varrho^{2m}}{1+|\varrho^m \lambda|^2}
 \sum_j (1+|\lambda|^2)(u,\phi_j)\phi_j\big)\\
&= \kappa_{\varrho}\big( \frac{1+|\lambda|^2}{1+|\varrho^m\lambda|^2}\varrho^{2m}
 \pi_{\KK_{\wedge,\lambda}}u
+ \frac{1-\varrho^{2m}}{1+|\varrho^m \lambda|^2}
\mathfrak p_{\KK_{\wedge,\lambda}}u\big).
\end{align*}
This gives the formula in the statement of the lemma.
\end{proof}

Note that if $u\in \Dom_{\wedge,\max}$, then
\begin{equation}\label{NormOfFrakP}
\|\mathfrak p_{\KK_{\wedge,\lambda}} u\|\leq \|u\|\leq \|u\|_{A_\wedge}, \quad \|\mathfrak p_{\KK_{\wedge,\lambda}} u\|_{A_\wedge} \leq {\sqrt{1+|\lambda|^2}}\|u\|_{A_\wedge}.
\end{equation}

\begin{lemma}\label{GoodLemmaWedge}
Let $\Lambda$ be some closed sector, let $R>0$, and let
\begin{equation*}
P(\lambda):x^{-m/2}L^2_b(Y^\wedge;E)\to \Dom_{\wedge,\max}
\end{equation*}
be a family of operators defined for $\lambda\in\Lambda_R$. Then
\begin{equation}\label{AbstractREstimatesWedge}
\|P(\lambda)\|_{\L(x^{-m/2}L^2_b)}
 \leq C/|\lambda|\quad{and}\quad
\|P(\lambda)\|_{\L(x^{-m/2}L^2_b,
 \Dom_{\wedge,\max})} \leq C
\end{equation}
hold for some $C>0$ and all $\lambda\in \Lambda_R$ if and only if
\begin{equation}\label{AbstractREstimatesWedgeKappa}
\|\kappa_{|\lambda|^{1/m}}^{-1}P(\lambda)\|_{\L(x^{-m/2}L^2_b,
 \Dom_{\wedge,\max})}\leq C/|\lambda|
\end{equation}
holds for some $C>0$ and all $\lambda\in \Lambda_R$.
\end{lemma}
\begin{proof}
Using that $A_\wedge \kappa_{|\lambda|^{1/m}}^{-1} P(\lambda)=|\lambda|^{-1} \kappa_{|\lambda|^{1/m}}^{-1} A_\wedge  P(\lambda)$, and that $\kappa_{|\lambda|^{1/m}}^{-1}$ is an isometry in $x^{-m/2}L^2_b(Y^\wedge;E)$, we obtain
\begin{align*}
\|\kappa_{|\lambda|^{1/m}}^{-1}P(\lambda)f\|_{A_\wedge }^2
&= \|A_\wedge \kappa_{|\lambda|^{1/m}}^{-1} P(\lambda)f\|^2 + \|\kappa_{|\lambda|^{1/m}}^{-1}P(\lambda)f\|^2\\
&= |\lambda|^{-2}\|\kappa_{|\lambda|^{1/m}}^{-1}A_\wedge P(\lambda)f\|^2+\|\kappa_{|\lambda|^{1/m}}^{-1}P(\lambda)f\|^2\\
&=
|\lambda|^{-2}\|A_\wedge P(\lambda)f\|^2+\|P(\lambda)f\|^2
\end{align*}
if $f\in x^{-m/2}L^2_b(Y^\wedge;E)$. Thus \eqref{AbstractREstimatesWedgeKappa} follows from \eqref{AbstractREstimatesWedge}.

Assume now that \eqref{AbstractREstimatesWedgeKappa} holds and let $f\in x^{-m/2}L^2_b(Y^\wedge;E)$. Then
\begin{equation*}
\|P(\lambda)f\|=\|\kappa_{|\lambda|^{1/m}}^{-1}P(\lambda)f\|\leq \|\kappa_{|\lambda|^{1/m}}^{-1}P(\lambda)f\|_{A_\wedge}
\end{equation*}
gives the first estimate in \eqref{AbstractREstimatesWedge}. To obtain the second, write $\|P(\lambda)f\|_{A_\wedge}^2$ as
\begin{equation*}
\|A_\wedge P(\lambda)f\|^2+\|P(\lambda)f\|^2 =\|\kappa_{|\lambda|^{1/m}}^{-1}A_\wedge P(\lambda)f\|^2+\|\kappa_{|\lambda|^{1/m}}^{-1}P(\lambda)f\|^2
\end{equation*}
and use the $\kappa$-homogeneity of $A_\wedge $ to conclude that
\begin{align*}
\|P(\lambda)f\|_{A_\wedge}^2
&= |\lambda|^2\| A_\wedge \kappa_{|\lambda|^{1/m}}^{-1} P(\lambda)f\|^2 + \|\kappa_{|\lambda|^{1/m}}^{-1}P(\lambda)f\|^2\\
&\leq(|\lambda|^2 +1)\|\kappa_{|\lambda|^{1/m}}^{-1} P(\lambda)f\|_{A_\wedge}^2.
\end{align*}
The second estimate in \eqref{AbstractREstimatesWedge} follows from this.
\end{proof}

\begin{corollary}
Let $\Dom\in \mathfrak G_{\wedge,\alpha}$, let $\Lambda$ be a closed sector. Then $\Lambda$ is a sector of minimal growth for $A_{\wedge,\Dom}$ if and only if there are $C$, $R>0$ such that
\begin{equation}\label{EstimateforBWedgeDom}
\|\kappa_{|\lambda|^{1/m}}^{-1}B_{\wedge,\Dom}(\lambda)\|_{\L(x^{-m/2}L^2_b,
 \Dom_{\wedge,\max})}\leq C/|\lambda|,\quad \lambda\in \Lambda_R.
\end{equation}
\end{corollary}

Both $B_{\wedge,\min}(\lambda)$ and $B_{\wedge,\max}(\lambda)$ satisfy
\eqref{AbstractREstimatesWedge}, therefore \eqref{AbstractREstimatesWedgeKappa} for any closed sector $\Lambda$ with $\Lambda\minus 0\subset \open \Lambda_\alpha$. In the case of $B_{\wedge,\min}(\lambda)$, the first of the estimates in \eqref{AbstractREstimatesWedge} is \eqref{EsimateForBWedgeMin}. To prove the second we note that
\begin{equation*}
A_\wedge B_{\wedge,\min}(\lambda)=\pi_{\RR_{\wedge,\lambda}} +\lambda B_{\wedge,\min}(\lambda)
\end{equation*}
where $\pi_{\RR_{\wedge,\lambda}}:x^{-m/2}L^2_b(Y^\wedge;E)\to x^{-m/2}L^2_b(Y^\wedge;E)$ is the orthogonal projection on $\RR_{\wedge,\lambda}$. The norm of this operator is $1$, and $\|\lambda B_{\wedge,\min}(\lambda)\|_{\L(x^{-m/2}L^2_b)}$ is bounded independently of $\lambda$ when $\lambda\in \Lambda$ and $|\lambda|$ is large. The argument for $B_{\wedge,\max}(\lambda)$ is analogous, using \eqref{EstimateForBWedgeMax} and the fact that this operator is a right inverse for $A_\wedge-\lambda$.

\begin{proof}[Proof of Theorem \ref{NecAndSuffWedge}]
We will prove that \eqref{NecAndSuffWedgeCond} is equivalent to \eqref{EstimateforBWedgeDom}. Recalling the formula \eqref{FormulaForBWedgeDomFirst} for $B_{\wedge,\Dom}(\lambda)$ and that $B_{\wedge,\max}(\lambda)$ satisfies \eqref{AbstractREstimatesWedgeKappa}, we see that $B_{\wedge,\Dom}(\lambda)$ satisfies \eqref{EstimateforBWedgeDom} if and only if
\begin{equation*}
\|\kappa_{|\lambda|^{1/m}}^{-1}\big(I- B_{\wedge,\min}(\lambda)(A_\wedge-\lambda)\big)\pi_{\KK_{\wedge,\lambda},\Dom} B_{\wedge,\max}(\lambda)
\|_{\L(x^{-m/2}L^2_b,\Dom_{\max})} \leq C/|\lambda|
\end{equation*}
for $\lambda\in \Lambda$, $|\lambda|$ large. We have
\begin{multline*}
\kappa_{|\lambda|^{1/m}}^{-1}\big(I- B_{\wedge,\min}(\lambda)(A_\wedge-\lambda)\big)\pi_{\KK_{\wedge,\lambda},\Dom} B_{\wedge,\max}(\lambda)
=\\
\big(I- B_{\wedge,\min}(\hat \lambda)(A_\wedge-\hat \lambda)\big)\kappa_{|\lambda|^{1/m}}^{-1}\pi_{\KK_{\wedge,\lambda},\Dom} \kappa_{|\lambda|^{1/m}} \kappa_{|\lambda|^{1/m}}^{-1} B_{\wedge,\max}(\lambda)
\end{multline*}
Evidently $\KK_{\wedge,\lambda}\cap \Dom=0$ if and only if $\kappa_{|\lambda|^{1/m}}^{-1}\KK_{\wedge,\lambda}\cap \kappa_{|\lambda|^{1/m}}^{-1}\Dom=0$. By Lemma \ref{KappaOnKwedgeLambda},  $\kappa_{|\lambda|^{1/m}}^{-1}\KK_{\wedge,\lambda} = \KK_{\wedge,\hat\lambda}$, and it is not hard to see that
\begin{equation*}
\kappa_{|\lambda|^{1/m}}^{-1}\pi_{\KK_{\wedge,\lambda},\Dom} \kappa_{|\lambda|^{1/m}} = \pi_{\KK_{\wedge,\hat\lambda}, \kappa_{|\lambda|^{1/m}}^{-1}\Dom}.
\end{equation*}
Using that $I- B_{\wedge,\min}(\lambda)(A_\wedge-\lambda)$ and $\pi_{\KK_{\wedge,\lambda},\Dom}$ both vanish on $\Dom_{\wedge,\min}$ regardless of $\lambda$ and $\Dom$, we arrive at the conclusion that $B_{\wedge,\Dom}(\lambda)$ satisfies \eqref{EstimateforBWedgeDom} if and only if the norm of
\begin{equation}\label{BareWedge}
\big(I- B_{\wedge,\min}(\hat\lambda)(A_\wedge-\hat\lambda)\big)\pi_{\wedge,\max}\, \pi_{\KK_{\wedge,\hat\lambda},\kappa_{|\lambda|^{1/m}}^{-1}\Dom}\, \pi_{\wedge,\max} \kappa_{|\lambda|^{1/m}}^{-1} B_{\wedge,\max}(\lambda)
\end{equation}
as an operator $x^{-m/2}L^2_b(Y^\wedge;E)\to\Dom_{\wedge,\max}$ is bounded by $C/|\lambda|$ for some $C$ if $\lambda\in \Lambda$, $|\lambda|$ large. By Lemma \ref{GoodLemmaWedge},
\begin{equation*}
\|\kappa_{|\lambda|^{1/m}}^{-1} B_{\wedge,\max}(\lambda)
\|_{\L(x^{-m/2}L^2_b,\Dom_{\wedge,\max})}\leq C/|\lambda|
\end{equation*}
for $\lambda\in \Lambda$, $|\lambda|$ large. Evidently
\begin{equation*}
\|I- B_{\wedge,\min}(\hat\lambda)(A_\wedge-\hat\lambda)\|_{\L(\Dom_{\wedge,\max})}
\end{equation*}
is bounded independently of $\lambda$, $\lambda\in \Lambda\minus 0$. Thus if \eqref{NecAndSuffWedgeCond} holds, then the norm of the operator \eqref{BareWedge} is bounded by $C/|\lambda|$ for some $C$ when $\lambda\in \Lambda$, $|\lambda|$ large.

Conversely, suppose that the norm of the operator \eqref{BareWedge} is bounded by $C/|\lambda|$ for some $C$ when $\lambda\in \Lambda$, $|\lambda|$ large. Composing with $\pi_{\wedge,\max}$ on the left we get that the norm of
\begin{equation*}
\pi_{\wedge,\max}\, \pi_{\KK_{\wedge,\hat\lambda},\kappa_{|\lambda|^{1/m}}^{-1}\Dom}\, \pi_{\wedge,\max} \kappa_{|\lambda|^{1/m}}^{-1} B_{\wedge,\max}(\lambda)
\end{equation*}
as an operator $x^{-m/2}L^2_b(Y^\wedge;E)\to\Sing_{\wedge,\max}$ satisfies the same estimate. Using the formula \eqref{FormulaBWedgeMax} for $B_{\wedge,\max}(\lambda)$ we get
\begin{multline*}
\pi_{\wedge,\max}\,\pi_{\KK_{\wedge,\hat\lambda},\kappa_{|\lambda|^{1/m}}^{-1}\Dom}\kappa_{|\lambda|^{1/m}}^{-1}  B_{\wedge,\max}(\lambda) =\\
|\lambda|^{-1} \pi_{\wedge,\max}\,\pi_{\KK_{\wedge,\hat\lambda},\kappa_{|\lambda|^{1/m}}^{-1}\Dom} \big(\id -  \frac{1-|\lambda|^2}{1+|\lambda|^2} \mathfrak p_{{\KK_\wedge,\hat\lambda}}\big) B_{\wedge,\max}(\hat\lambda)\kappa_{|\lambda|^{1/m}}^{-1}.
\end{multline*}
We dismiss the factor $\kappa_{|\lambda|^{1/m}}^{-1}$ at the end of the last formula, since this is an isometry on $x^{-m/2}L^2_b(Y^\wedge;E)$. Since $\mathfrak p_{\KK_{\wedge,\hat\lambda}}$ has range in $\KK_{\wedge,\hat\lambda}$,
\begin{multline*}
|\lambda|^{-1} \pi_{\wedge,\max}\,\pi_{\KK_{\wedge,\hat\lambda},\kappa_{|\lambda|^{1/m}}^{-1}\Dom} \bigg(\frac{1-|\lambda|^2}{1+|\lambda|^2}\bigg) \mathfrak p_{\KK_{\wedge,\hat\lambda}} B_{\wedge,\max}(\hat\lambda)\\
=
|\lambda|^{-1} \frac{1-|\lambda|^2}{1+|\lambda|^2}\pi_{\wedge,\max}  \mathfrak p_{\KK_{\wedge,\hat\lambda}} B_{\wedge,\max}(\hat\lambda).
\end{multline*}
This operator $x^{-m/2}L^2_b(Y^\wedge;E)\to\Dom_{\wedge,\max}$ evidently has norm $O(|\lambda|^{-1})$ if $\lambda\in \Lambda$, $|\lambda|\to\infty$ (cf. \eqref{NormOfFrakP}). We conclude that if the norm of \eqref{BareWedge} is bounded as indicated, then the norm of
\begin{equation}\label{BareWedge1}
\pi_{\wedge,\max}\,\pi_{\KK_{\wedge,\hat\lambda},\kappa_{|\lambda|^{1/m}}^{-1}\Dom}B_{\wedge,\max}(\hat\lambda)
\end{equation}
is bounded by a constant when $\lambda\in \Lambda$, $|\lambda|$ large. The operator $A_\wedge-\hat\lambda$ on $\Dom_{\wedge,\max}$ satisfies $\|A_\wedge-\hat\lambda\|_{\L(\Dom_{\wedge,\max},x^{-m/2}L^2_b)}\leq 1$. So composing the operator \eqref{BareWedge1} with $A_\wedge-\hat\lambda$ on the right we get that the norm of
\begin{equation*}
\pi_{\wedge,\max}\,\pi_{\KK_{\wedge,\hat\lambda},\kappa_{|\lambda|^{1/m}}^{-1}\Dom} (\id -\pi_{\KK_{\wedge,\hat \lambda}})=
\pi_{\wedge,\max}\,\pi_{\KK_{\wedge,\hat\lambda},\kappa_{|\lambda|^{1/m}}^{-1}\Dom}\,\pi_{\wedge,\max}-\pi_{\wedge,\max}\,\pi_{\KK_{\wedge,\hat \lambda}}
\end{equation*}
satisfies the same estimate. Since $\|\pi_{\wedge,\max}\pi_{\KK_{\wedge,\hat \lambda}}\|_{\L(\Dom_{\wedge,\max})}\leq 1$, and using that  $\pi_{\KK_{\wedge,\hat\lambda},\Dom}=\pi_{\KK_{\wedge,\hat\lambda},\Dom}\pi_{\wedge,\max}$, we obtain that if $\Lambda$ is a sector of minimal growth for $A_{\wedge,\Dom}$, then
\begin{equation*}
\|\pi_{\wedge,\max}\,\pi_{\KK_{\wedge,\hat\lambda},\kappa_{|\lambda|^{1/m}}^{-1}\Dom}\,\pi_{\wedge,\max}\|_{\L(\Dom_{\wedge,\max})}
\end{equation*}
is bounded for $\lambda\in \Lambda$, $|\lambda|$ large. This completes the proof of the theorem.
\end{proof}

Let $\KK_{\wedge,\max}(\lambda) = \pi_{\wedge,\max}\KK_{\wedge,\lambda}$. Let $\Dom\in \mathfrak G_{\wedge,\alpha}$, let $\lambda_0\in \open\Lambda_\alpha$ be such that $|\lambda_0|=1$, and let $R>0$. The condition that
\begin{equation}\label{weakRayCondition}
\varrho^m\lambda_0\in\resolv A_{\wedge,\Dom} \text{ for } \varrho\geq R
\end{equation}
is equivalent to the statement that $\KK_{\wedge,\varrho^m\lambda_0}\cap \Dom=0$ for $\varrho\geq R$, which in turn is equivalent to the condition that $\KK_{\wedge,\lambda_0}\cap \kappa_\varrho^{-1}\Dom=0$ for $\varrho\geq R$. Thus, since $\KK_{\wedge,\lambda_0}\cap \kappa_\varrho^{-1}\Dom=0$ if and only if $\pi_{\wedge,\max}\KK_{\wedge,\lambda_0}\cap \pi_{\wedge,\max}\kappa_\varrho^{-1}\Dom=0$, the condition in \eqref{weakRayCondition} is equivalent to the statement that the curve $\gamma$ defined by
\begin{equation}\label{CurveOfDomains}
[R\neutral],\neutral(\infty)\ni \varrho\mapsto \gamma(\varrho)=\pi_{\wedge,\max}\kappa_\varrho^{-1}\Dom \in \Gr_{d_\alpha''}(\Sing_{\wedge,\max})
\end{equation}
does not intersect the variety $\mathfrak V_{\KK_{\wedge,\max}(\lambda_0)}$ introduced in Definition \ref{BadVariety} (with $\Sing_{\wedge,\max}$ in place of $\Sing_{\max})$. With the proof of Lemma \ref{ReductionOfPiToEmax}, $\pi_{\wedge,\max}\,\pi_{\KK_{\wedge,\lambda_0},\kappa_\varrho^{-1}\Dom}\big|_{\Sing_{\wedge,\max}}$ is the projection on $\KK_{\wedge,\max}(\lambda_0)$ according to the decomposition
\begin{equation*}
\KK_{\wedge,\max}(\lambda_0)\oplus \pi_{\wedge,\max}\kappa_\varrho^{-1} \Dom = \Sing_{\wedge,\max}.
\end{equation*}
Thus if there is a neighborhood $U$ of $\mathfrak V_{\KK_{\wedge,\max}(\lambda_0)}$ in $\Gr_{d_\alpha''}(\Sing_{\wedge,\max})$ such that $\gamma(\varrho)\notin U$ if $\varrho$ is sufficiently large, then Lemma \ref{BoundsViaGrassmann} gives that
\begin{equation*}
\|\pi_{\KK_{\wedge,\max}(\lambda_0),\gamma(\varrho)}\|\leq \frac{C}{\delta(\KK_{\wedge,\max}(\lambda_0),\gamma(\varrho))}
\end{equation*}
is bounded as $\varrho\to\infty$. Therefore the necessary condition of Theorem \ref{NecAndSuffWedge} is satisfied, and we get:

\begin{theorem}
Let $\lambda_0\in \bgres A_\wedge$ belong to $\open\Lambda_\alpha$. Let $\Dom\in \mathfrak G_{\wedge,\alpha}$ and suppose that there is a neighborhood $U\subset \Gr_{d_\alpha''}(\Sing_{\wedge,\max})$ of $\mathfrak V_{\KK_{\wedge,\max}(\lambda_0)}$ such that $\pi_{\wedge,\max}\kappa_{\varrho}^{-1}\Dom\notin U$ for all sufficiently large $\varrho$. Then there is a closed sector $\Lambda$ containing $\lambda_0$ which is a sector of minimal growth for $A_{\wedge,\Dom}$.
\end{theorem}

%%%%%%%%%%%%%%%%%%%%%%%%%%%%%%%%%%%%%%%%%%%%%%%%%%%%%%%%%%%
%%%%%%%%%%%%%%%%%%%%%%%%%%%%%%%%%%%%%%%%%%%%%%%%%%%%%%%%%%%

\section{Resolvents}
\label{sec-ResolventsOnM}

We will now prove the analogue of Theorem \ref{NecAndSuffWedge} for $A\in x^{-m}\Diff^m_b(M;E)$.

Define $A_\varrho=\varrho^{-m}\kappa_\varrho^{-1}A\kappa_\varrho$. Then $A_\varrho$ is $c$-elliptic, since $A$ is assumed to be $c$-elliptic, cf. \eqref{csymArho}. Using that $A_\varrho-A$ belongs to $x^{-m+1}\Diff^m_b(M;E)$, \cite[Proposition 4.1(1)]{GiMe01} gives the first formula in
\begin{equation*}
\Dom_{\min}(A_\varrho)=\Dom_{\min}(A),\quad \Dom_{\max}(A_\varrho)=\kappa_\varrho^{-1}\Dom_{\max}(A).
\end{equation*}
The second is obtained using that $\kappa_\varrho$ preserves $C_0^\infty(\open M;E)$ and $x^{-m/2}L^2_b(M;E)$. We will write $\Dom_{\varrho,\max}$ instead of $\Dom_{\max}(A_\varrho)$, and $\Dom_{\min}$ instead of $\Dom_{\min}(A_\varrho)$. Generally we prepend the symbol $\varrho$ to subindices of objects associated with $A_\varrho$ originally associated with $A$. In particular,
\begin{equation*}
\Sing_{\varrho,\max}=\Dom_{\min}^\perp
\end{equation*}
with the orthogonal computed in $\Dom_{\varrho,\max}$ using the inner product
\begin{equation*}
(u,v)_{A_\varrho}=(A_\varrho u,A_\varrho v)+(u,v)
\end{equation*}
of $\Dom_{\varrho,\max}$, and $\pi_{\varrho,\max}: \Dom_{\varrho,\max}\to \Dom_{\varrho,\max}$ is the orthogonal projection on $\Sing_{\varrho,\max}$. It is not hard to verify that
\begin{equation}\label{EmaxAsKernelRho}
\Sing_{\varrho,\max}=\kappa_\varrho^{-1}[\ker (A^\star A+\varrho^{2m})\cap \Dom_{\max}],
\end{equation}
cf. Lemma \ref{EmaxAsKernel}.

Using
\begin{equation}\label{variantHomogeneityOfA}
\varrho^{-m}\kappa_\varrho^{-1}(A-\varrho^m\lambda)\kappa_\varrho = A_\varrho-\lambda
\end{equation}
we see that
\begin{equation*}
\bgres A_\varrho = \varrho^{-m}\bgres A.
\end{equation*}
For $\lambda\in \bgres A_\varrho$ let $\KK_{\varrho,\lambda}=\ker (A_{\varrho,\Dom_{\max}}-\lambda)$. Then $\KK_{\varrho,\lambda/\varrho^m}=\kappa_\varrho^{-1}\KK_{\lambda}$. If $\Dom\in \mathfrak G$, then $\kappa_\varrho^{-1}\Dom \in \mathfrak G_\varrho$, and if $\lambda\in \resolv A_\Dom$, then $\lambda/\varrho^m\in \resolv A_{\varrho,\kappa_\varrho^{-1}\Dom}$. It is easy to verify that $\Dom_{\varrho,\max}=\KK_{\varrho,\lambda/\varrho^m}\oplus \kappa_{\varrho}^{-1}\Dom$ and that
\begin{equation}\label{ConjOfProjLambdaDom}
\kappa_\varrho^{-1}\pi_{\KK_\lambda,\Dom}\kappa_\varrho = \pi_{\KK_{\varrho,\lambda/\varrho^m},\kappa_\varrho^{-1}\Dom}.
\end{equation}

Let $B_{\varrho,\min}(\lambda) = \varrho^m\kappa_{\varrho}^{-1} B_{\min}(\varrho^m\lambda) \kappa_\varrho$. This is a left inverse of $A_{\varrho,\Dom_{\min}}-\lambda$. The operator $B_{\varrho,\min}(\lambda)$ has range in $\Dom_{\min}$ since this subspace is $\kappa$-invariant and the range of $B_{\min}(\varrho^m\lambda)$ is $\Dom_{\min}$.

\begin{theorem}\label{NecAndSuff}
Let $\Dom\in \mathfrak G$ and let $\Lambda$ be a closed sector. Then $\Lambda$ is a sector of minimal growth for $A_\Dom$ if and only if there are positive constants $C$, $R$ such that $\Lambda_R\subset \resolv A_\Dom$,
\begin{equation}\label{MinimalGrowthMinMax}
\|B_{\min}(\lambda)\|_{\L(x^{-m/2}L^2_b)}\leq C/|\lambda|,\quad \|B_{\max}(\lambda)\|_{\L(x^{-m/2}L^2_b)}\leq C/|\lambda|,
\end{equation}
and
\begin{equation}\label{NecAndSuffCond}
\|\pi_{|\lambda|^{1/m},\max}\pi_{\KK_{|\lambda|^{1/m},\hat\lambda}, \kappa_{|\lambda|^{1/m}}^{-1}\Dom} \big|_{\Sing_{|\lambda|^{1/m},\max}}\|_{\L(\Dom_{|\lambda|^{1/m},\max})}\leq C,\quad \lambda\in \Lambda_R.
\end{equation}
\end{theorem}

The proof requires a number of analogues of results obtained in the previous section. Their proofs parallel those in that section.

\begin{lemma}\label{GoodLemma}
Let $\Lambda$ be some closed sector, let $R>0$, and let
\begin{equation*}
P(\lambda):x^{-m/2}L^2_b(M;E)\to \Dom_{\max}
\end{equation*}
be a family of operators defined for $\lambda\in\Lambda_R$. Then
\begin{equation}\label{AbstractREstimates}
\|P(\lambda)\|_{\L(x^{-m/2}L^2_b)}
 \leq C/|\lambda|\quad{and}\quad
\|P(\lambda)\|_{\L(x^{-m/2}L^2_b,
 \Dom_{\max})} \leq C
\end{equation}
hold for some $C>0$ and all $\lambda\in \Lambda_R$ if and only if
\begin{equation}\label{AbstractREstimatesKappa}
\|\kappa_{|\lambda|^{1/m}}^{-1}P(\lambda)\|_{\L(x^{-m/2}L^2_b,
 \Dom_{|\lambda|^{1/m},\max})}\leq C/|\lambda|
\end{equation}
holds for some $C>0$ and all $\lambda\in \Lambda_R$.
\end{lemma}
\begin{proof}
Using that $A_{|\lambda|^{1/m}} \kappa_{|\lambda|^{1/m}}^{-1} P(\lambda)=|\lambda|^{-1} \kappa_{|\lambda|^{1/m}}^{-1} A P(\lambda)$, and that $\kappa_{|\lambda|^{1/m}}^{-1}$ is an isometry in $x^{-m/2}L^2_b(M;E)$, we obtain
\begin{align*}
\|\kappa_{|\lambda|^{1/m}}^{-1}P(\lambda)f\|_{A_{|\lambda|^{1/m}}}^2
&=\|A_{|\lambda|^{1/m}} \kappa_{|\lambda|^{1/m}}^{-1}P(\lambda)f\|^2+\|\kappa_{|\lambda|^{1/m}}^{-1}P(\lambda)f\|^2\\
&=|\lambda|^{-2}\|\kappa_{|\lambda|^{1/m}}^{-1}A_{|\lambda|^{1/m}} P(\lambda)f\|^2+\|\kappa_{|\lambda|^{1/m}}^{-1}P(\lambda)f\|^2\\
&= |\lambda|^{-2}\|A_{|\lambda|^{1/m}} P(\lambda)f\|^2+\|P(\lambda)f\|^2
\end{align*}
if $f\in x^{-m/2}L^2_b(M;E)$. Thus \eqref{AbstractREstimatesKappa} follows from \eqref{AbstractREstimates}.

Assume now that \eqref{AbstractREstimatesKappa} holds and let $f\in x^{-m/2}L^2_b(M;E)$. Then
\begin{equation*}
\|P(\lambda)f\|=\|\kappa_{|\lambda|^{1/m}}^{-1}P(\lambda)f\|\leq \|\kappa_{|\lambda|^{1/m}}^{-1}P(\lambda)f\|_{A_{|\lambda|^{1/m}}}
\end{equation*}
gives the first estimate in \eqref{AbstractREstimates}. To obtain the second, write $\|P(\lambda)f\|_A^2$ as
\begin{equation*}
\|A P(\lambda)f\|^2+\|P(\lambda)f\|^2 =\|\kappa_{|\lambda|^{1/m}}^{-1}A P(\lambda)f\|^2+\|\kappa_{|\lambda|^{1/m}}^{-1}P(\lambda)f\|^2
\end{equation*}
and use the definition of $A_{|\lambda|^{1/m}}$ to conclude that
\begin{multline*}
\|P(\lambda)f\|_A^2 = |\lambda|^2 \|A_{|\lambda|^{1/m}} \kappa_{|\lambda|^{1/m}}^{-1}P(\lambda)f\|^2+\|\kappa_{|\lambda|^{1/m}}^{-1}P(\lambda)f\|^2\\
\leq(|\lambda|^2 +1)\|\kappa_{|\lambda|^{1/m}}^{-1} P(\lambda)f\|_{A_{|\lambda|^{1/m}}}^2.
\end{multline*}
The second estimate in \eqref{AbstractREstimates} follows from this.
\end{proof}

\begin{corollary}
Let $\Dom\in \mathfrak G$, let $\Lambda$ be a closed sector. Then $\Lambda$ is a sector of minimal growth for $A_\Dom$ if and only if there are positive constants $C$, $R$ such that $\Lambda_R\subset \resolv A_\Dom$ and
\begin{equation}\label{EstimateForBDom}
\|\kappa_{|\lambda|^{1/m}}^{-1}B_{\Dom}(\lambda)\|_{\L(x^{-m/2}L^2_b,
 \Dom_{|\lambda|^{1/m},\max})}\leq C/|\lambda|,\quad \lambda\in \Lambda_R.
\end{equation}
\end{corollary}

\begin{proof}[Proof of Theorem \ref{NecAndSuff}: Sufficiency of the condition]
We will show that \eqref{MinimalGrowthMinMax} and \eqref{NecAndSuffCond} imply \eqref{EstimateForBDom}. Since $B_{\min}(\lambda)$ and $B_{\max}(\lambda)$ satisfy the estimate in \eqref{MinimalGrowthMinMax}, and since these estimates imply for each of them the second estimate in \eqref{AbstractREstimates}, we obtain that $\kappa_{|\lambda|^{1/m}}^{-1}B_{\min}(\lambda)$ and $\kappa_{|\lambda|^{1/m}}^{-1}B_{\max}(\lambda)$ both satisfy \eqref{AbstractREstimatesKappa}. In particular, to prove \eqref{EstimateForBDom} we only need to prove that for some $C$,
\begin{equation*}
\|\kappa_{|\lambda|^{1/m}}^{-1}\big(B_{\Dom}(\lambda)-B_{\max}(\lambda)\big)\|_{\L(x^{-m/2}L^2_b,
 \Dom_{|\lambda|^{1/m},\max})}\leq C/|\lambda|,\quad \lambda\in \Lambda_R.
\end{equation*}
Writing $B_\Dom(\lambda)$ as in \eqref{ResolventOnDBis} we get
\begin{equation*}
\kappa_{|\lambda|^{1/m}}^{-1}\big(B_{\max}(\lambda)-B_{\Dom}(\lambda)\big)= \kappa_{|\lambda|^{1/m}}^{-1} \big(\id-B_{\min}(\lambda)(A-\lambda)\big) \pi_{\KK_{\lambda},\Dom}B_{\max}(\lambda).
\end{equation*}
We rewrite the right hand side as
\begin{equation}\label{Intermezzo}
\big(\id - B_{|\lambda|^{1/m},\min}(\hat\lambda)(A_{|\lambda|^{1/m}} - \hat \lambda)\big) \kappa_{|\lambda|^{1/m}}^{-1} \,\pi_{\KK_{\lambda},\Dom}\,  \kappa_{|\lambda|^{1/m}} \kappa_{|\lambda|^{1/m}}^{-1} B_{\max}(\lambda).
\end{equation}
Using that $\id - B_{|\lambda|^{1/m},\min}(\hat\lambda) (A_{|\lambda|^{1/m}} - \hat \lambda)$ vanishes on $\Dom_{\min}$, the identity \eqref{ConjOfProjLambdaDom}, and that $\Dom_{\min}\subset \kappa_{|\lambda|^{1/m}}^{-1}\Dom$ we replace the factor $\kappa_{|\lambda|^{1/m}}^{-1} \pi_{\KK_{\lambda},\Dom}\kappa_{|\lambda|^{1/m}}$ in \eqref{Intermezzo} by
\begin{equation*}
\pi_{|\lambda|^{1/m},\max}\,\pi_{\KK_{|\lambda|^{1/m},\hat\lambda},\kappa_{|\lambda|^{1/m}}^{-1}\Dom}\, \pi_{|\lambda|^{1/m},\max}.
\end{equation*}
By hypothesis the norms of these operators $\Dom_{|\lambda|^{1/m},\max}\to \Sing_{|\lambda|^{1/m},\max}$ are uniformly bounded when $\lambda\in \Lambda_R$. It is easy to verify that the norm of
\begin{equation*}
\id - B_{|\lambda|^{1/m},\min}(\hat\lambda)(A_{|\lambda|^{1/m}} - \hat \lambda):\Dom_{|\lambda|^{1/m},\max}\to \Dom_{|\lambda|^{1/m},\max}
\end{equation*}
is bounded independently of $\lambda$, $\lambda\in \Lambda_R$. Finally, as already discussed,
\begin{equation*}
\|\kappa_{|\lambda|^{1/m}}^{-1}B_{\max}(\lambda)\|_{\L(x^{-m/2}L^2_b,
 \Dom_{|\lambda|^{1/m},\max})}\leq C/|\lambda|
\end{equation*}
holds for some $C>0$ and all $\lambda\in \Lambda_R$. Altogether these estimates give \eqref{EstimateForBDom}.
\end{proof}

To prove the necessity of the condition in Theorem \ref{NecAndSuff} we will need two lemmas.

\begin{lemma}\label{MinimalGrwothBMinBMax}
Suppose that $\Dom\in \mathfrak G$ and that the closed sector $\Lambda$ is a sector of minimal growth for $A_\Dom$. Then there are positive constants $R$ and $C$ such that
\begin{equation*}
\|B_{\min}(\lambda)\|_{\L(x^{-m/2}L^2_b)}\leq C/|\lambda|,\quad \|B_{\max}(\lambda)\|_{\L(x^{-m/2}L^2_b)}\leq C/|\lambda|
\end{equation*}
for $\lambda\in \Lambda_R$.
\end{lemma}
This is a direct consequence of the formulas
\begin{equation*}
B_{\min}(\lambda)=B_{\Dom}(\lambda)\pi_{\RR_\lambda},\quad
B_{\max}(\lambda)=B_{\Dom}(\lambda)-\pi_{\KK_\lambda}B_{\Dom}(\lambda),
\end{equation*}
cf. \eqref{CorrectionFormulaMax} and \eqref{CorrectionFormulaMin} valid for $\lambda\in \resolv A_\Dom$.

\begin{lemma}\label{LemmaOnProjections}
Let $\mathfrak p_{\KK_{\varrho,\lambda}}:x^{-m/2}L^2_b(M;E)\to x^{-m/2}L^2_b(M;E)$ be the orthogonal projection on  $\KK_{\varrho,\lambda}$ regarded as a subspace of $x^{-m/2}L^2_b(M;E)$. Then
\begin{equation}\label{RelationBetweenPis}
\kappa_\varrho^{-1} \pi_{\KK_{\varrho^m\lambda}}\kappa_\varrho =
\frac{1+|\lambda|^2}{1+|\varrho^m\lambda|^2}\varrho^{2m} \pi_{\KK_{\varrho,\lambda}} + \frac{1-\varrho^{2m}}{1+|\varrho^m\lambda|^2}\mathfrak p_{\KK_{\varrho,\lambda}}.
\end{equation}
Moreover, $\|\mathfrak p_{\KK_{\varrho,\lambda}}\|_{\L{(\Dom_{\varrho,\max})}}\leq \sqrt {1 + |\lambda|^2}$.
\end{lemma}

\begin{proof}
The proof of \eqref{RelationBetweenPis} parallels that of Lemma \ref{LemmaOnProjectionsWedge}. Let $\phi_1,\dotsc,\phi_{d'}$ be an $A_\varrho$-orthonormal basis of $\KK_{\varrho,\lambda}=\kappa_\varrho^{-1}\KK_{\varrho^m\lambda}$. Then
\begin{equation*}
 \delta_{jk}=(\phi_j,\phi_k)_{A_\varrho} = (1+|\lambda|^2)(\phi_j,\phi_k).
\end{equation*}
In particular, the $\sqrt{1+|\lambda|^2}\phi_j\in \KK_{\varrho,\lambda}$ are orthonormal in $x^{-m/2}L^2_b(M;E)$. On the other hand, using that $\kappa_\varrho$ is an isometry on $x^{-m/2}L^2_b(M;E)$,
\begin{equation*}
(\kappa_\varrho\phi_j,\kappa_\varrho\phi_k)_{A}
= \varrho^{2m}(\phi_j,\phi_k)_{A_\varrho} + (1-\varrho^{2m})(\phi_j,\phi_k)
= \frac{1+|\varrho^m \lambda|^2}{1+|\lambda|^2}\,\delta_{jk}.
\end{equation*}
This gives an $A$-orthonormal basis of $\KK_{\varrho^m\lambda}$, and if $u\in \KK_{\varrho,\lambda}$, then
\begin{align*}
\pi_{\KK_{\varrho^m\lambda}} \kappa_\varrho u
&= \frac{1+|\lambda|^2}{1+|\varrho^m \lambda|^2}
\sum_j(\kappa_\varrho u,\kappa_\varrho\phi_j)_{A}\kappa_\varrho\phi_j\\
&= \frac{1+|\lambda|^2}{1+|\varrho^m \lambda|^2}
\sum_j\big[\varrho^{2m}(u,\phi_j)_{A_\varrho} + (1-\varrho^{2m})(u,\phi_j)
\big]\kappa_\varrho\phi_j \\
&= \kappa_\varrho
\big( \frac{1+|\lambda|^2}{1+|\varrho^m\lambda|^2}\varrho^{2m}
\pi_{\KK_{\varrho,\lambda}}u + \frac{1-\varrho^{2m}}{1+|\varrho^m \lambda|^2}
\mathfrak p_{\KK_{\varrho,\lambda}}u\big).
\end{align*}
Thus \eqref{RelationBetweenPis} follows. The estimate of the norm of $\mathfrak p_{\KK_{\varrho,\lambda}}$ is elementary.
\end{proof}

\begin{proof}[Proof of Theorem \ref{NecAndSuff}: Necessity of the condition]
Suppose that $\Lambda$ is a sector of minimal growth for $A_\Dom$. By Lemma \ref{MinimalGrwothBMinBMax}, \eqref{MinimalGrowthMinMax} holds. In particular there are $C$, $R$ such that the operator
\begin{equation*}
\kappa_{|\lambda|^{1/m}}^{-1}(B_{\max}(\lambda)-B_{\Dom}(\lambda))
= \kappa_{|\lambda|^{1/m}} ^{-1}\big(\id-B_{\min}(\lambda)(A-\lambda)\big) \pi_{\KK_\lambda,\Dom}B_{\max}(\lambda)
\end{equation*}
as an element of $\L(x^{-m/2}L^2_b,\Dom_{|\lambda|^{1/m},\max})$
has norm bounded by $C/|\lambda|$ if $\lambda\in \Lambda_R$. Composing with $\pi_{|\lambda|^{1/m},\max}$ on the left, and using that $\kappa_{|\lambda|^{1/m}}^{-1}$ preserves $\Dom_{\min}$, we conclude that
\begin{multline*}
\pi_{|\lambda|^{1/m},\max}\,\kappa_{|\lambda|^{1/m}} ^{-1}\big(\id-B_{\min}(\lambda)(A-\lambda)\big) \pi_{\KK_\lambda,\Dom}B_{\max}(\lambda) \\ = \pi_{|\lambda|^{1/m},\max}\,\kappa_{|\lambda|^{1/m}} ^{-1}\, \pi_{\KK_\lambda,\Dom}B_{\max}(\lambda)
\end{multline*}
satisfies the same estimate. The operator
\begin{equation*}
(A-\lambda)\kappa_{|\lambda|^{1/m}} =  |\lambda|\kappa_{|\lambda|^{1/m}}(A_{|\lambda|^{1/m}}-\hat\lambda),
\end{equation*}
as an element of $\L(\Dom_{|\lambda|^{1/m},\max},x^{-m/2}L^2_b)$, has norm bounded by $2|\lambda|$, $\lambda \ne 0$. Thus the operator
\begin{multline*}
\pi_{|\lambda|^{1/m},\max}\,\kappa_{|\lambda|^{1/m}} ^{-1}\, \pi_{\KK_\lambda,\Dom}B_{\max}(\lambda)(A-\lambda)\kappa_{|\lambda|^{1/m}}\\
=\pi_{|\lambda|^{1/m},\max}\,\kappa_{|\lambda|^{1/m}} ^{-1}\, \pi_{\KK_\lambda,\Dom}(I-\pi_{\KK_{\lambda}})\kappa_{|\lambda|^{1/m}},
\end{multline*}
as an element of $\L(\Dom_{|\lambda|^{1/m},\max})$, has norm bounded by a constant independent of $\lambda\in \Lambda_R$. Since $\pi_{\KK_\lambda,\Dom}\pi_{\KK_{\lambda}}
= \pi_{\KK_{\lambda}}$,
\begin{equation*}
\kappa_{|\lambda|^{1/m}} ^{-1} \pi_{\KK_\lambda,\Dom}\pi_{\KK_{\lambda}}\kappa_{|\lambda|^{1/m}}
=
\frac{2|\lambda|^2}{1+|\lambda|^2} \pi_{\KK_{|\lambda|^{1/m},\hat\lambda}} + \frac{1 - |\lambda|^2}{1+|\lambda|^2}\mathfrak p_{\KK_{|\lambda|^{1/m},\hat\lambda}}
\end{equation*}
using \eqref{RelationBetweenPis}. Thus
\begin{equation*}
\|\pi_{|\lambda|^{1/m},\max}\kappa_{|\lambda|^{1/m}} ^{-1} \pi_{\KK_\lambda,\Dom}\pi_{\KK_{\lambda}}\kappa_{|\lambda|^{1/m}}\|_{\L(\Dom_{|\lambda|^{1/m},\max})}\leq C,\quad \lambda\in \Lambda_R,
\end{equation*}
for some $C$ and consequently also
\begin{equation*}
\|\pi_{|\lambda|^{1/m},\max}\kappa_{|\lambda|^{1/m}} ^{-1} \pi_{\KK_\lambda,\Dom}\kappa_{|\lambda|^{1/m}}\|_{\L(\Dom_{|\lambda|^{1/m}, \max})}\leq C,\quad \lambda\in \Lambda_R,
\end{equation*}
for some other $C$. Using \eqref{ConjOfProjLambdaDom} we conclude that in particular
\begin{equation*}
\|\pi_{|\lambda|^{1/m},\max} \pi_{\KK_{|\lambda|^{1/m},\hat\lambda},\kappa_{|\lambda|^{1/m}}^{-1}\Dom}\big|_{\Sing_{|\lambda|^{1/m},\max}}\|_{\L(\Dom_{|\lambda|^{1/m},\max})}\leq C,\quad \lambda\in \Lambda_R,
\end{equation*}
This completes the proof of the necessity of the condition.
\end{proof}

%%%%%%%%%%%%%%%%%%%%%%%%%%%%%%%%%%%%%%%%%%%%%%%%%%%%%%%%%%%
%%%%%%%%%%%%%%%%%%%%%%%%%%%%%%%%%%%%%%%%%%%%%%%%%%%%%%%%%%%

\end{document}